%% file: manuscript.tex
\documentclass[review,onefignum,onetabnum]{siamart171218}

\input{ex_shared}

\ifpdf
\hypersetup{
	pdftitle={A new approach to data assimilation initialization problems with sparse data using multiple cost functions},
	pdfauthor={D. J. Albers, G. Hripcsak, L. Mamykina, M. Sirlanci, E. Tabak}
}
\fi

\usepackage{amsxtra}
\usepackage{amsfonts}
\usepackage{amssymb, verbatim}
\usepackage{bbm}
\usepackage[utf8]{inputenc}
\usepackage{graphicx}
\usepackage{multicol}
\usepackage{amsmath}
\usepackage{subcaption}
\usepackage{multimedia}
\usepackage{multirow}
\usepackage{xcolor}

\newlength{\arrayrulewidthOriginal}

\graphicspath{ {./figures/} }

\begin{document}
	
	\maketitle
	
	\begin{abstract}
		This article develops a novel data assimilation methodology, addressing challenges that are common in real-world settings, such as severe sparsity of observations, lack of reliable models, and non-stationarity of the system dynamics. These challenges often cause identifiability issues and can confound model parameter initialization, both of which can  lead to estimated models with unrealistic qualitative dynamics and induce deeper parameter estimation errors. The proposed methodology’s objective function is constructed as a sum of components, each serving a different purpose: enforcing point-wise and distribution-wise agreement between data and model output, enforcing agreement of variables and parameters with a model provided, and penalizing unrealistic rapid parameter changes, unless they are due to external drivers or interventions. This methodology was motivated by, developed and evaluated in the context of estimating blood glucose levels in different medical settings. Both simulated and real data are used to evaluate the methodology from different perspectives, such as its ability to estimate unmeasured variables, its ability to reproduce the correct qualitative blood glucose dynamics, how it manages known non-stationarity, and how it performs when given a range of dense and severely sparse data. The results show that a multicomponent cost function can balance the minimization of point-wise errors with global properties, robustly preserving correct qualitative dynamics and managing data sparsity.
	\end{abstract}
	
	

\section{Introduction} 

This article proposes a methodology for data assimilation in complex systems displaying at least some of the following characteristics:
\begin{enumerate}

\item {\bf Lack of reliable models.} Most biomedical systems are only partially understood, and certainly rarely understood to the point where a well-defined model provides an accurate surrogate for the underlying biological processes. Even in situations where those processes can be captured by models to some degree --think general circulation models for weather forecasting-- typically a large number of model parameters remain undetermined.

\item {\bf Latent and emerging bulk variables.} The phase-space of complex systems is high dimensional, while only a handful of variables are observed systematically. The value of the remaining, \emph{latent} variables can only be inferred indirectly from the observations, if at all. Since the full set of latent variables is often too large and hard to pinpoint, it is useful to consider instead \emph{emerging bulk variables}, a small set of parameters representing a larger set of unmeasured variables and processes that affect the dynamics in a coherent way.

\item {\bf Sparse observations.} The system may be observed not only partially, but also at a sparse set of times, insufficient to even marginally resolve some of its dynamically significant time-scales. The timing of these observations may be structured (for instance by a medical protocol), event-driven or purely random.

\item {\bf Prior knowledge.} Even though a detailed model of the underlying processes may be unavailable, one often has some prior information on the system's dynamics, such as the expected presence of oscillations and their typical frequency and amplitude, or the difference in behavior between driven and unforced scenarios.

\item {\bf Non-stationarity.} The system's underlying dynamics may evolve over time, due either to external drivers or to un-modeled components of the system. Non-stationarity may manifest itself through a slow time modulation of the model's parameters.

\end{enumerate}

The authors have encountered problems with these characteristics in biomedical and computational physiology, specifically computing parameter initializations for data assimilation using sparse data \cite{ALBERS2023104477,albers_plos_comp_bio_DA_I,cenkf,houlihan} with the goal of estimating and forecasting endocrine function related to blood glucose regulation.  Contexts included critical care in the intensive care unit (ICU) \cite{melike2019simple,ALBERS2023104477,WANG2023104547,melike_g_control,vanherpe06}, type-2 diabetes \cite{albers_plos_comp_bio_DA_I,jamia_da,levine_albers_plos_comp_bio_DA_II} and outpatient or \emph{in the wild} settings \cite{miller2020learninginsulinglucosedynamicswild,wang2023learningabsorptionratesglucoseinsulin,pmlr-v235-zou24b}. The origin of these problems, summarized in subsection \ref{sec:motivation}, can be traced to sources that include model rigidity, data sparsity, non-stationarity and the presence of errors in the recorded measurement times. Another pervasive source of error is the use of incomplete loss functions, such as least square error, which under sparse data can completely misspecify the system's true dynamics (replacing for instance an unresolved oscillatory signal by its mean value.)  \emph{Importantly, while we have faced these problems in biomedicine, the sources of these problems and DA initialization more broadly generalize to a wide range of other fields} such as atmospheric physics\cite{aos_da_initialization,CoupledDataAssimilationandEnsembleInitializationwithApplicationtoMultiyearENSOPrediction}.

\medskip

The approach that we propose to address problems of this type uses a stochastic model for the evolution of the system's state, not necessarily accurate or fully based on field-specific knowledge, but rather, flexible and tailored to the qualitative nature of the expected dynamics, oscillatory in our study case. The model has  observable variables $\{x^j\}$ directly relatable to the observations $\{y^j\}$ at times $\{t^j\}$ of the system's state $\mathcal{Y}^{t_j}$, which represents the system as a whole in an abstract manner, latent or emerging bulk variables $\{z^j\}$ that are not observed, and parameters $\alpha$. Examples for $x$ and $z$ are the blood stream's glucose, and the insulin content and associated biophysical processes respectively. We propose to estimate the model's parameters and the system's current state, and forecast its evolution under candidate interventions, through a methodology that includes the following ingredients:
\begin{enumerate}
\item A quantification of the agreement between the $\{x^j\}$ and $\{y^j\}$ that is not only point but also distribution-wise. This addresses, for instance, the pitfall described above, of wrongly approximating an oscillatory signal by its mean: a constant value may be a model's local minimizer of the point-wise squared distance to an oscillatory signal, but not in a distributional sense.

\item A slow modulation over time of said distribution and of  the model's parameters $\alpha$, to capture non-stationarity. In other words, we allow both the distribution and the parameter values to flex (vary) within the estimation window to some degree.

\item A three staged methodology that first initializes the model's parameters $\alpha$ and the observable variables $\{x^j\}$, then maximizes an objective function only over the latent variables $\{z^j\}$, and only then performs a full-fledged maximization of the objective function over all variables and parameters.  The reasons for such a staged approach are that, without sensible model parameter values, it is hopeless to try to estimate the latent variables $z$ and, without sensible values for $z$, any use of the model to estimate $x$ and $\alpha$ is equally doomed. Thus the first estimate for the $x$ is based only on their closeness to the observations $y$ (point and distribution-wise), the first estimate for the model parameters $\alpha$ depends only on the observable variables $x$, not the latent $z$, and the first estimation for $z$ is based on the dynamical model using those temporarily estimated values for $x$ and $\alpha$.

\item A way to handle interventions not included in the model and potentially not fully measured. No bio-physical model can capture all possible interventions, which in our motivating example may range from regular interventions, such as food intake and administered insulin, to more individual ones, such as drug administration or unexpected external developments. For those interventions not explicitly accounted for in the model, we allow the model's state and parameters to behave discontinuously at the intervention times, where the jumps' permissible amplitude depends only on an estimate of the intervention's intensity.

\end{enumerate}

\subsection{Motivating example}
\label{sec:motivation}

This paper addresses two highly related problems encountered when estimating an imperfect model with data from a system that generates complex dynamics, is non-stationary, and is brutally under-measured \cite{ALBERS2023104477,houlihan}.  To demonstrate how standard methods can fail in such realistic settings, consider an application, estimating the states and parameters of a glucose-insulin model given sparse data taken during a stay at  an intensive care unit (ICU) \cite{ALBERS2023104477}. In this setting data includes continuously measured tube-administered nutrition, point-wise blood glucose measurements collected according to clinical protocols \cite{george_eval,parkes_grid}, usually about once an hour, \emph{administered} insulin and its type, but never plasma or interstitial insulin measurements.  In the ICU, where nutrition administration is roughly constant because of its administration via a tube feed, blood glucose and insulin oscillate out of phase on the order of tens of minutes  \cite{sturis_91}.  Given this situation, the model we estimate is the Ultradian model \cite{sturis_91,keenerII,albers_plos_comp_bio_DA_I,ALBERS2023104477,WANG2023104547} that includes blood glucose, blood insulin, remote (interstitial) insulin, and a three state delay approximated with the linear-chain-trick \cite{dde_lct}. We estimate the model with a version of an ensemble Kalman filter (EnKF), initializing the parameters at the nominal parameters for the model \cite{sturis_91,albers_plos_comp_bio_DA_I}. Because blood or remote insulin is never measured, the model is not identifiable, but can still be estimated with some accuracy under many circumstances \cite{albers_plos_comp_bio_DA_I,ALBERS2023104477,WANG2023104547,jamia_da}. We focus on this example for a few reasons. \emph{First,} this situation poses a real world problem \cite{real_world_evidence,33_hyperglycemia_review,32_glucose_control_ICU,ALBERS2023104477,WANG2023104547} and is  highly representative of many biomedical and health care settings as well as any situation where the model is low-dimensional, the system is non-stationary, and data are sparely measured.  \emph{Second}, within this real world context, the forecasting and estimation needs are diverse and dependent on the decision-needs. For example, there are situations where: \emph{(a)} only parameter estimates are important but next-step glucose prediction is not, \emph{(b)} next-step glucose prediction is important but accurate parameter estimates are not, \emph{(c)} only accurate point-wise glucose forecasting is important, and \emph{(d)} accurate parameter estimates and distributional glucose forecasting are important but accurate point-wise glucose forecasting is not. This final case is likely to be the most common situation in applied settings. \emph{Third}, we have a qualitative understanding of the underlying dynamics of the data-generating system, including both glucose-insulin system and the measurement processes and protocols \cite{jamia_phys_ehr,jamia_hcpmodel13}.  And \emph{fourth}, the measurements of this system are sparse in time compared to the complexity of the underlying dynamics but the measurement times may not be entirely random. 

Operationally, we assume a situation where: 
\begin{enumerate}
\item[a.] we know the generating dynamics correspond to a driven and damped system that elicits a nondescript noisy oscillatory orbit;
\item[b.] the measurements are taken infrequently enough that the spectral composition of the orbit is difficult or impossible to resolve from the data alone;
\item[c.] the reporting of the measurement times has nontrivial errors.
\end{enumerate}
In such a setting it may be difficult or impossible to find a set of initial conditions and parameters that  synchronize the model to these data with an oscillatory solution. For example, if the model can exhibit a fixed point, a plausible error minimizing parameter set includes parameters for which the model remains at a fixed point whose mean is the mean of these data. However, because we know the generating system has oscillatory dynamics, we know that such solutions to the inverse problem \textbf{must be wrong}. This is because even if the estimated parameters minimize the least squared error, they produce qualitatively wrong dynamics. The root of this problem lies in how we compute and conceive of what it means for a model estimate to be \emph{wrong.} Given a model that is capable of producing dynamics that are qualitatively similar to known generating dynamics, model parameters that elicit qualitative dynamics that are substantially different from the qualitative dynamics of the generating system are not an accurate representation of the generating system. When the generating dynamics are complex and sparsely measured, it is frequently encountered that least squared error minimizing parameters produce qualitatively wrong generating dynamics.

\begin{figure}[t!]
	\centering
	\begin{subfigure}{0.5\textwidth}
		\includegraphics[width=\linewidth]{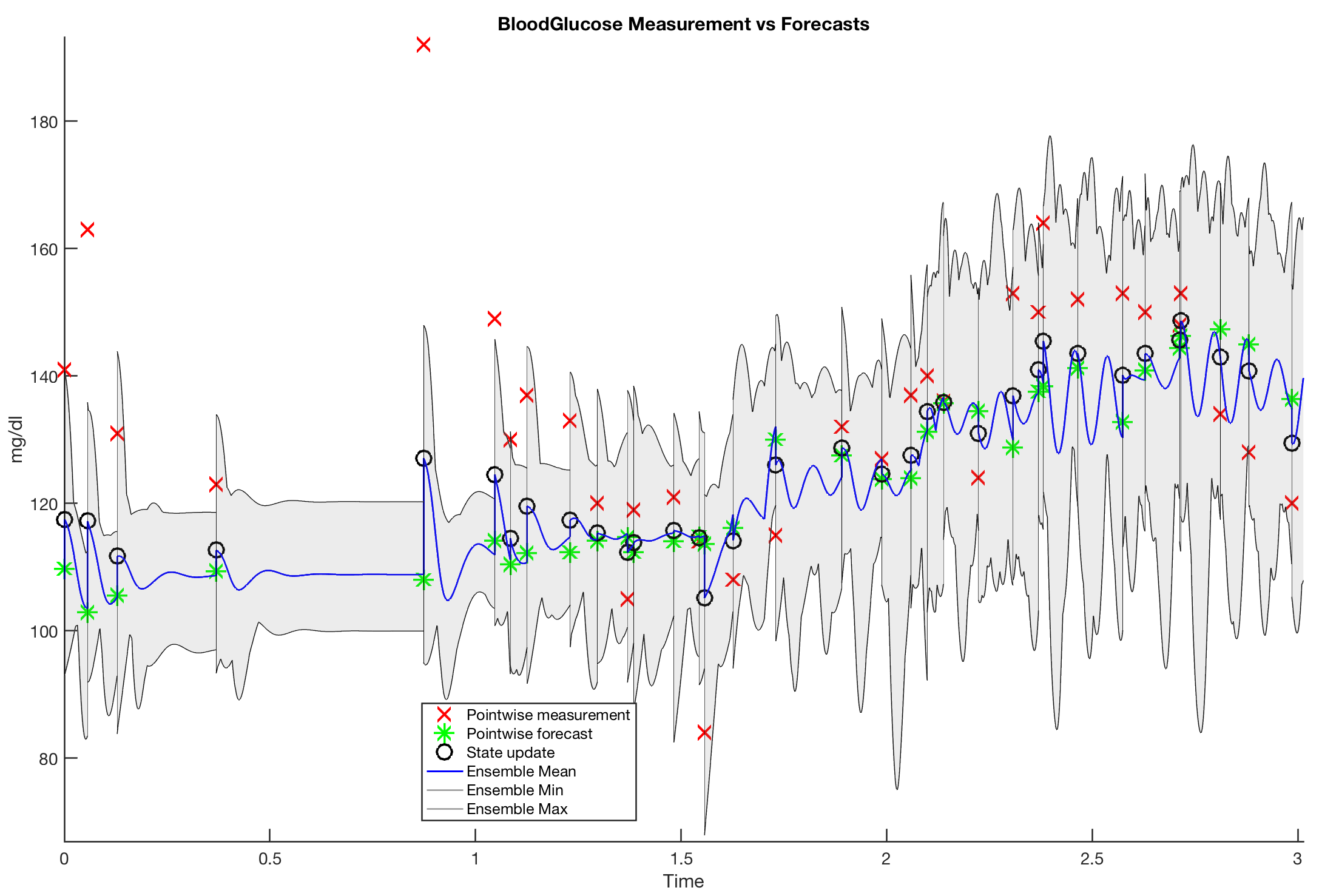}
	\caption{}
		\label{fig:motivation0a} 
	\end{subfigure}%
	\begin{subfigure}{0.5\textwidth}
		\includegraphics[width=\linewidth]{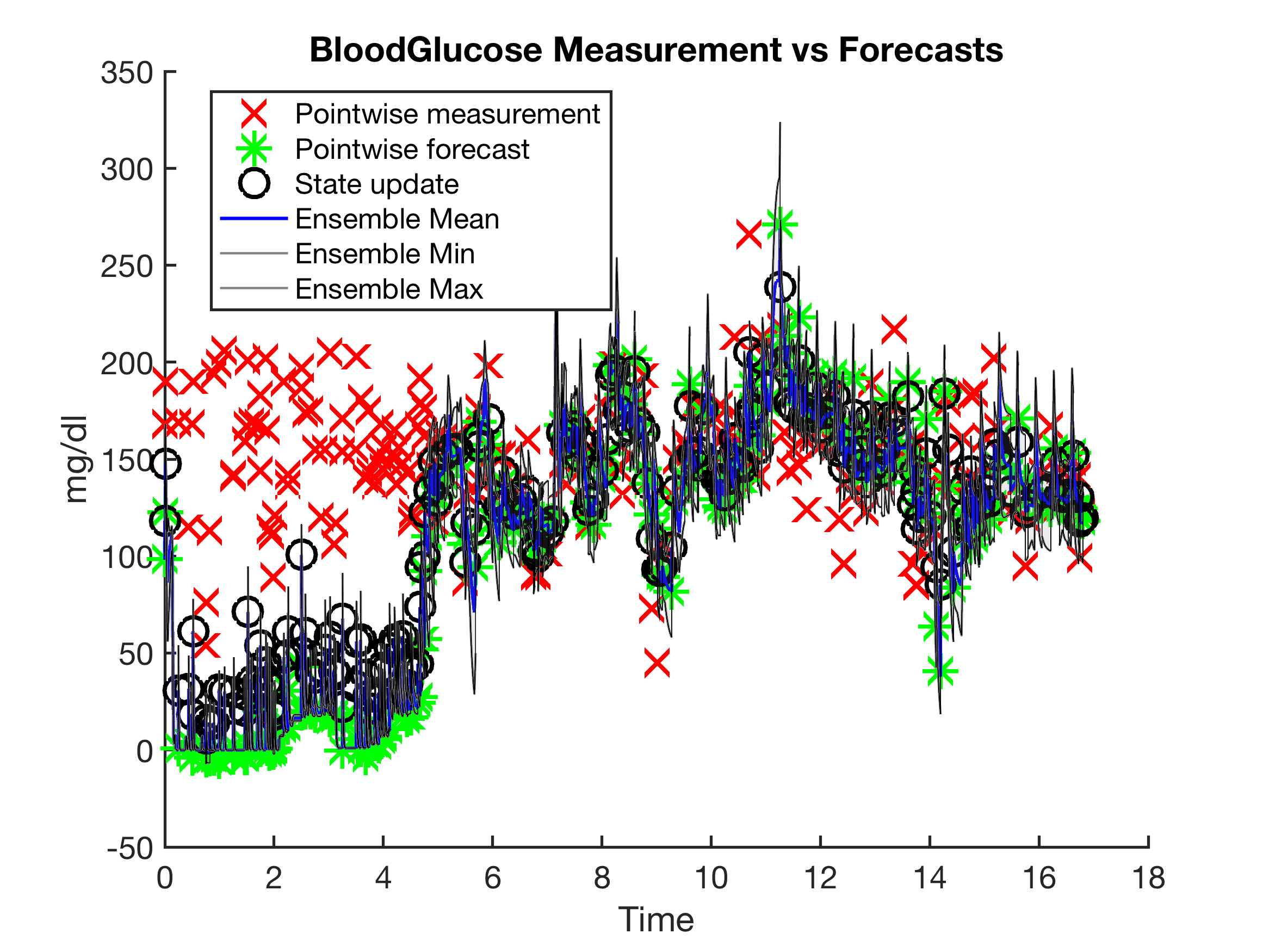}
		\caption{}
		\label{fig:motivation0b} 
	\end{subfigure}

	\caption{Estimating and forecasting  glucose trajectories for two patients in the ICU. The plot on the left shows accurate estimation and forecasting almost immediately ($<1$ day, $<6$ data points ) while the plot on the right shows poor estimation and prediction until about day $5$ ($\sim 125$ data points). Sources leading to model estimation accuracy include a complex interplay between data sparsity leading to model identifiability problems (Figs. \ref{fig:motivation1}-\ref{fig:motivation2}) and model initialization.}
	\label{fig:motivation0} 
\end{figure}

\begin{figure}[t!]
	\centering
	\begin{subfigure}{0.5\textwidth}
		\includegraphics[width=\linewidth]{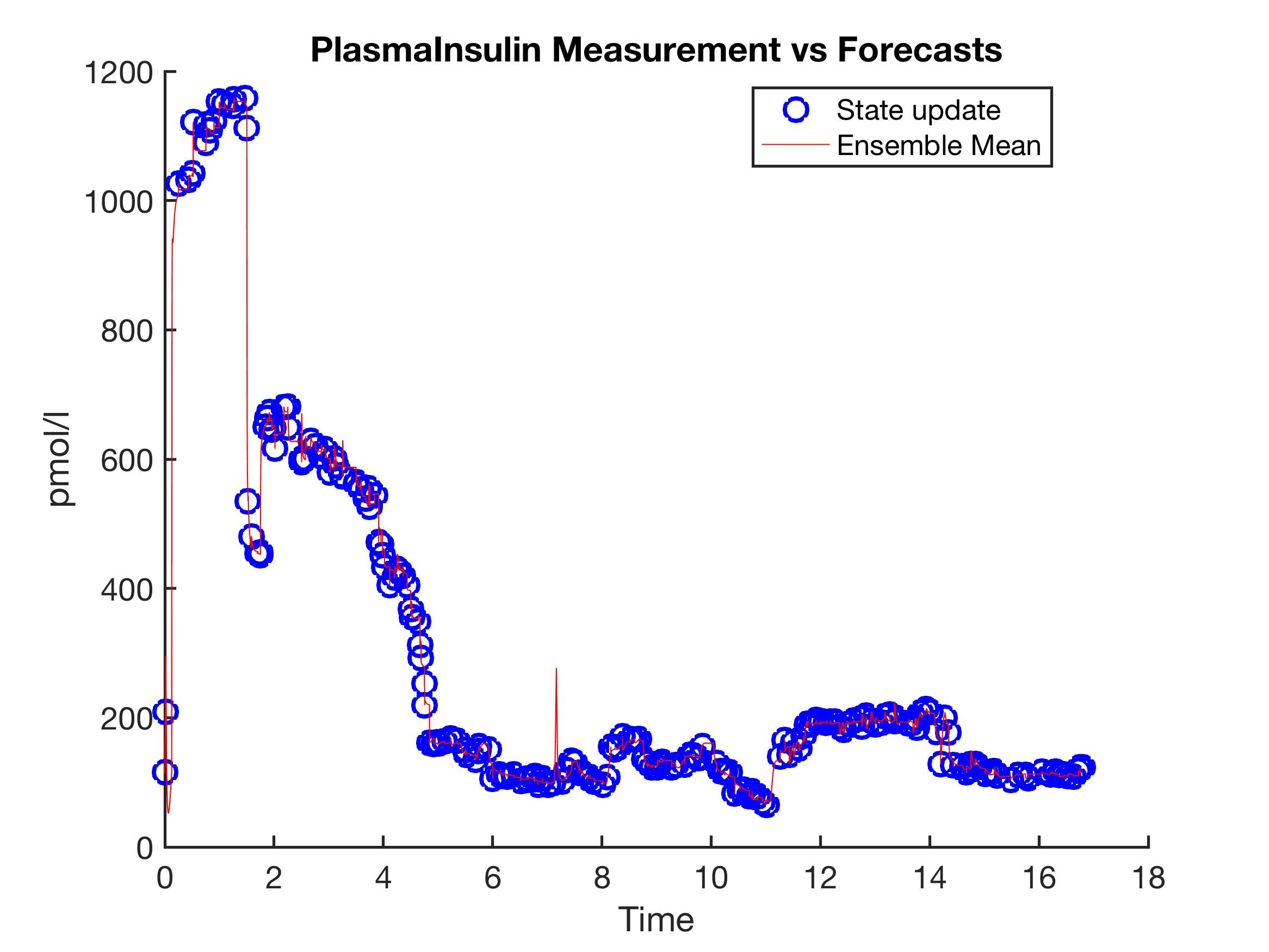}
		\label{fig:motivation1a} 
	\end{subfigure}%
	\begin{subfigure}{0.5\textwidth}
		\includegraphics[width=\linewidth]{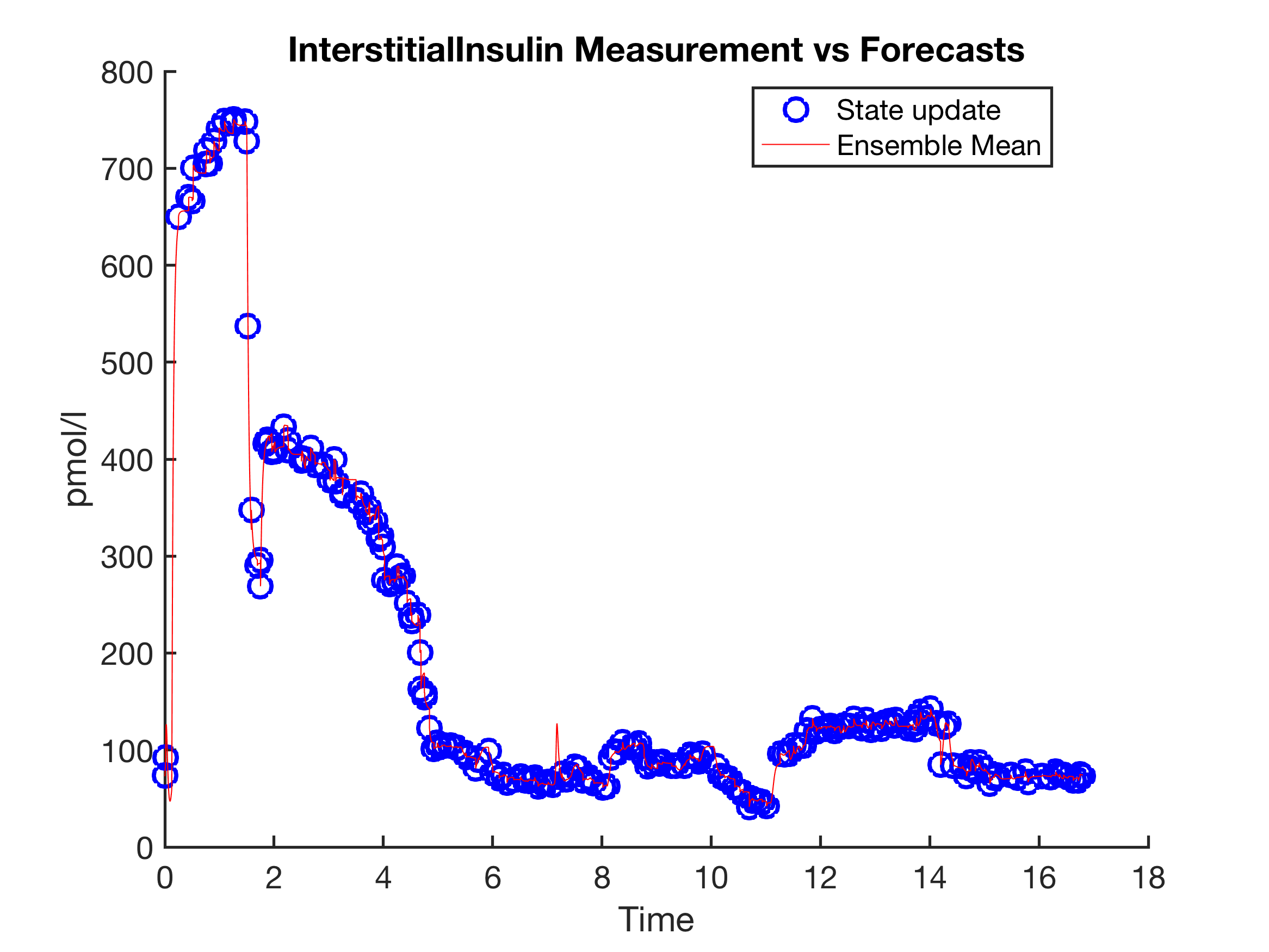}
		\label{fig:motivation1b} 
	\end{subfigure}
    
    \caption{In the ICU, insulin, one of the states that defines the glucose-insulin system and should be in the range of 25-400 picomoles per liter (pmol/l), is never measured. This can lead to model estimation and initialization problems, as seen in Fig. \ref{fig:motivation0}. Here we see the estimated interstitial and plasma insulin levels that are driving the forecasting errors seen in Fig. \ref{fig:motivation0b}. Note that after about $5$ days the model does eventually entrain to the patient and the insulin estimates take physiologically plausible values.}
    \label{fig:motivation1} 
\end{figure}

To demonstrate the problem, begin with the case where patient $426$'s data are estimated with the ensemble Kalman filter (EnKF), shown in Fig. \ref{fig:motivation0}a. Here we can see that the model estimation converges within about $1-1.5$ days, the ensemble mean oscillates as we might expect, the uncertainty in the estimate is reasonable, and these data are reasonably well predicted by the model-based forecast \cite{ALBERS2023104477}. This is an example of what successful estimation and forecasting looks given such data. A less successful estimation attempt is shown in Fig. \ref{fig:motivation0}b where we estimate patient $593$'s data. In this case it takes $\sim 4.5$ days for the model to entrain to the patient. Initially the EnKF estimates produce blood glucose values far below what would represent a living human, and the glycemic dynamics are not oscillatory. Fig. \ref{fig:motivation1} reveals one source of the problem, the unmeasured insulin states (blood and remote insulin) take values far above what is possible until about day $4.5$, the point where the model has converged such that the ensemble mean is tracking the glycemic dynamics; clearly the model was not initialized with accurate parameters. These two problems have two naively obvious solutions.


\begin{figure}[t!]
 	\centering
 	\begin{subfigure}{0.5\textwidth}
        \includegraphics[width=0.9\linewidth]{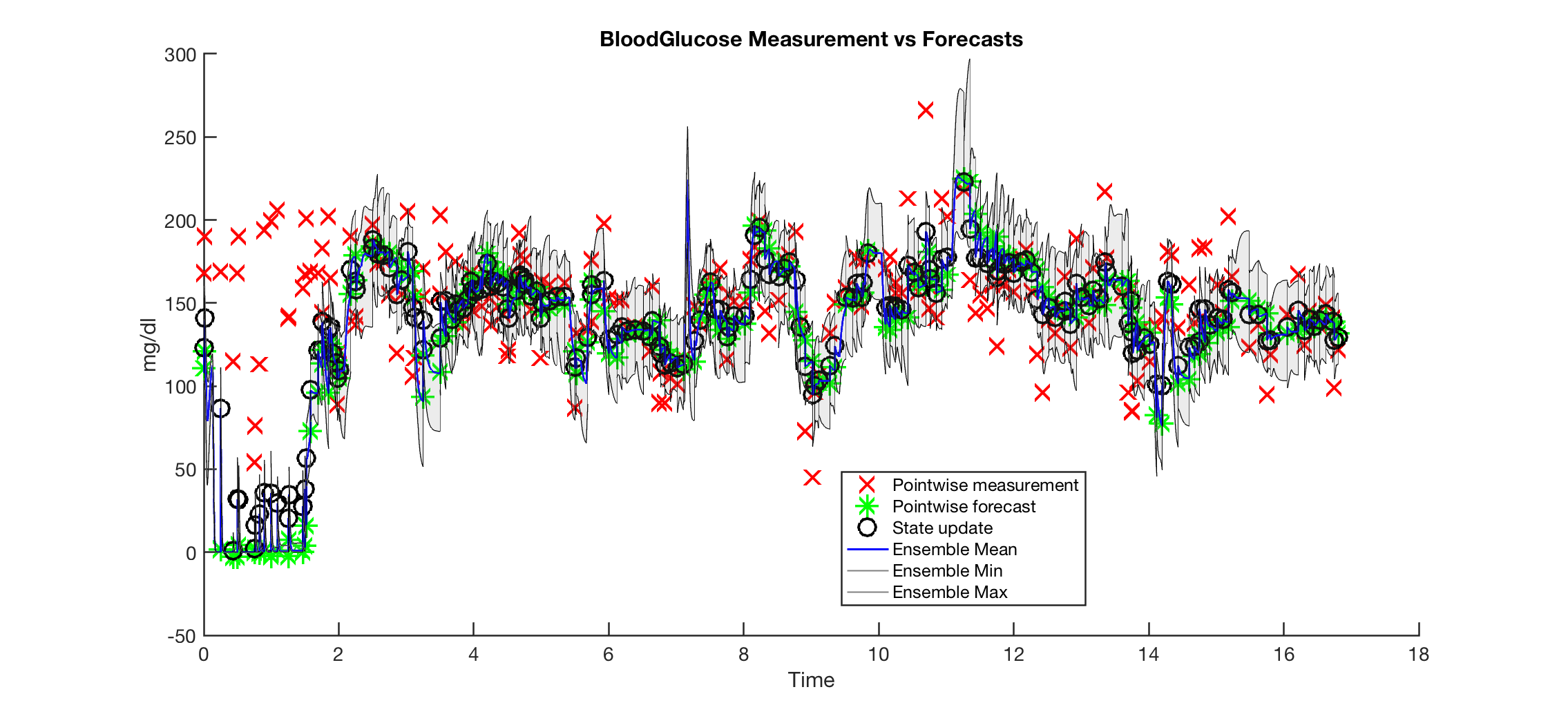}
       \caption{}
        \label{fig:motivation2a}
    \end{subfigure}%
	\begin{subfigure}{0.5\textwidth}
		\includegraphics[width=0.9\linewidth]{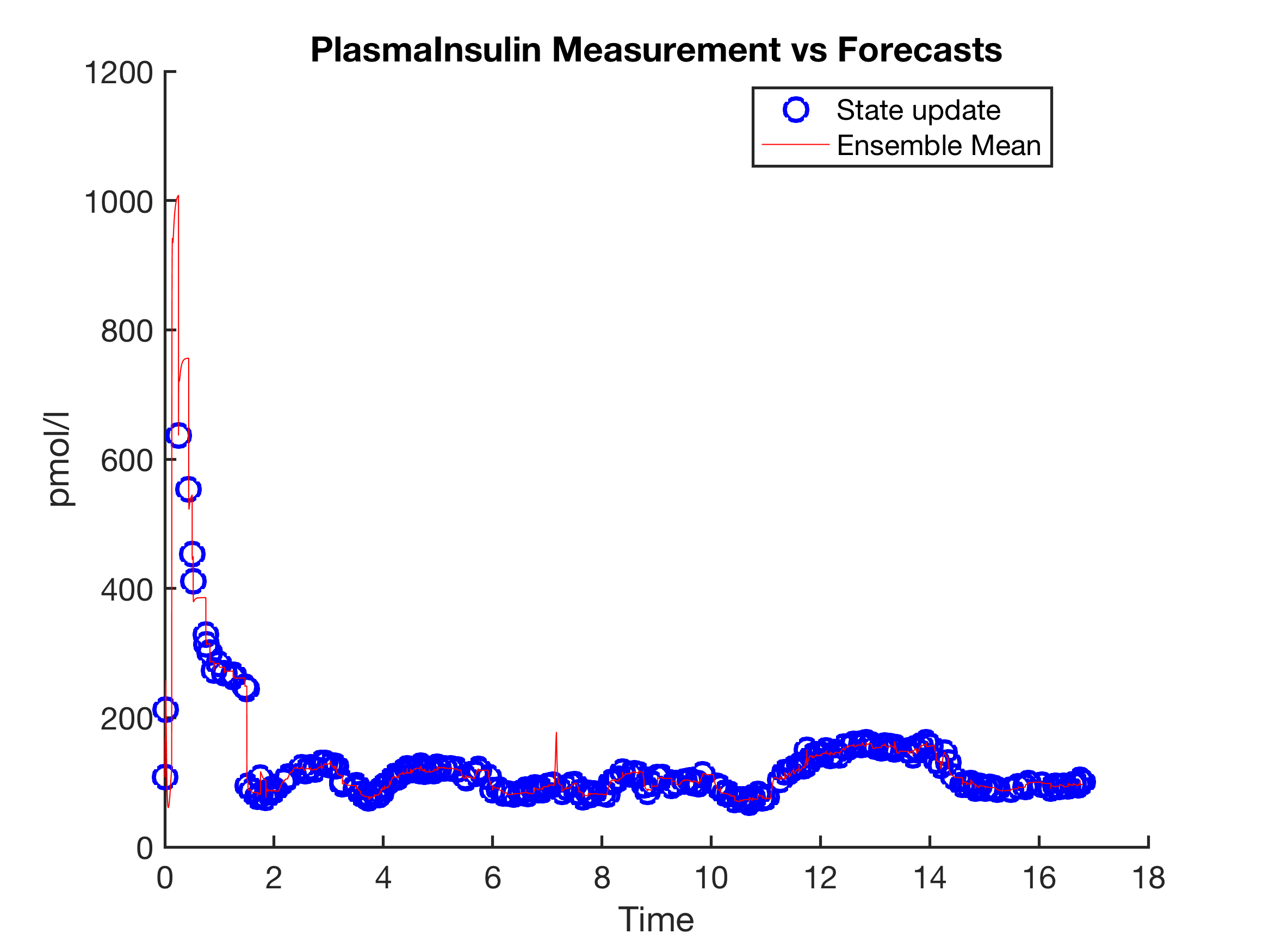}
		\caption{}
		\label{fig:motivation2b}
	\end{subfigure}

	\begin{subfigure}{0.5\textwidth}
		\includegraphics[width=0.9\linewidth]{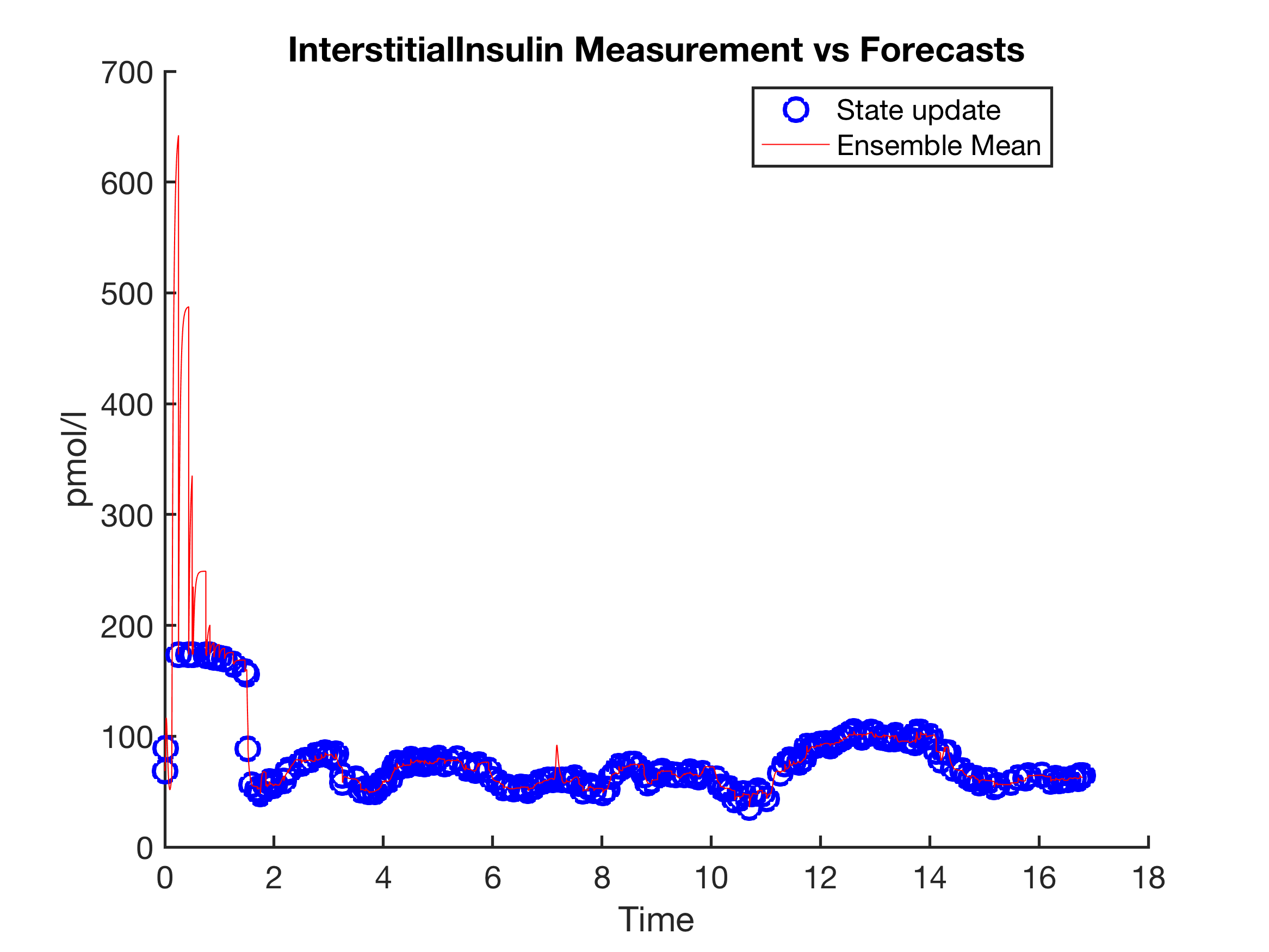}
		\caption{}
		\label{fig:motivation2f}
	\end{subfigure}%
	\begin{subfigure}{0.5\textwidth}
		\includegraphics[width=0.9\linewidth]{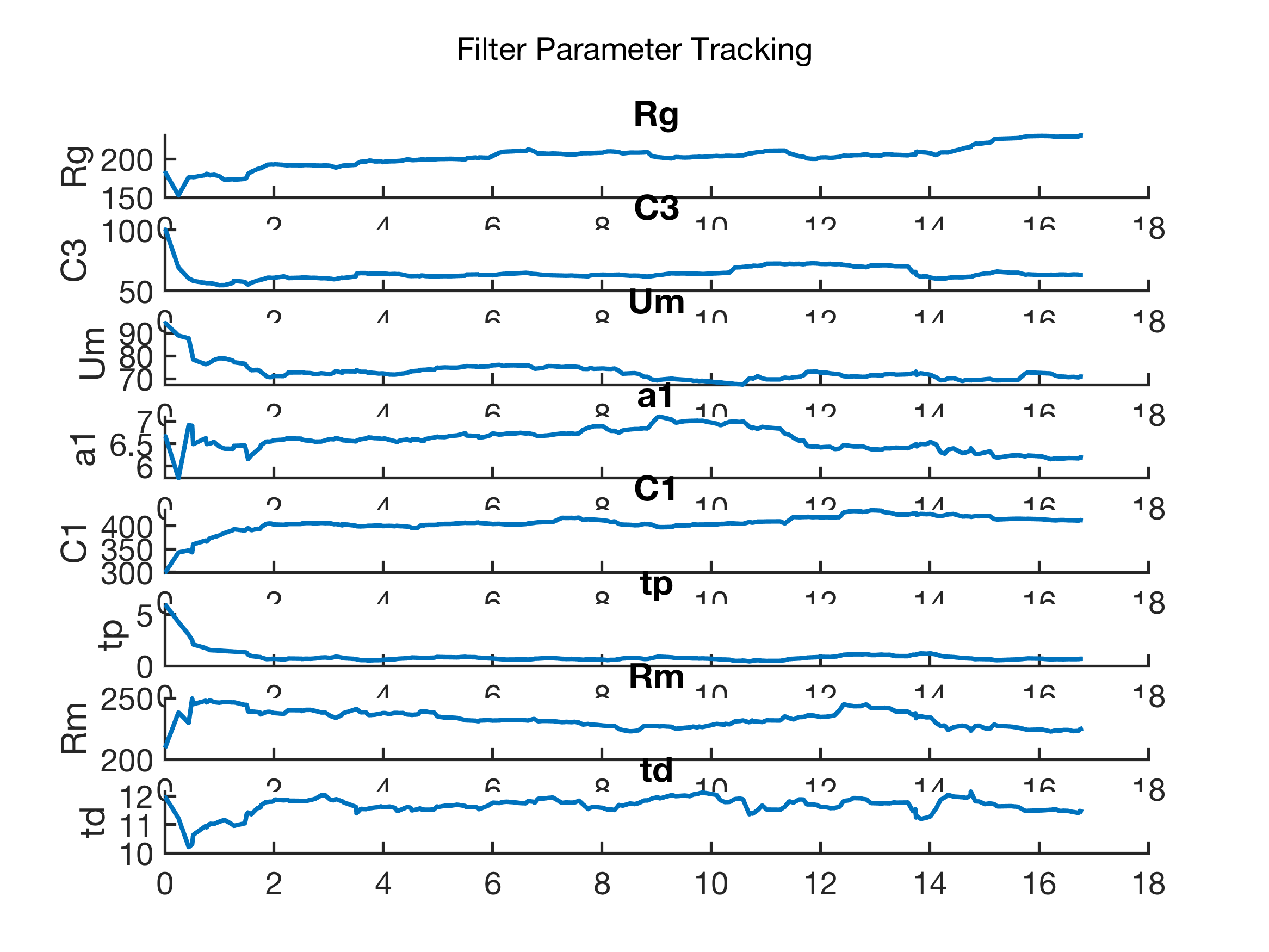}
		\caption{}
		\label{fig:motivation2c}
	\end{subfigure}

	\begin{subfigure}{0.5\textwidth}
		\includegraphics[width=0.9\linewidth]{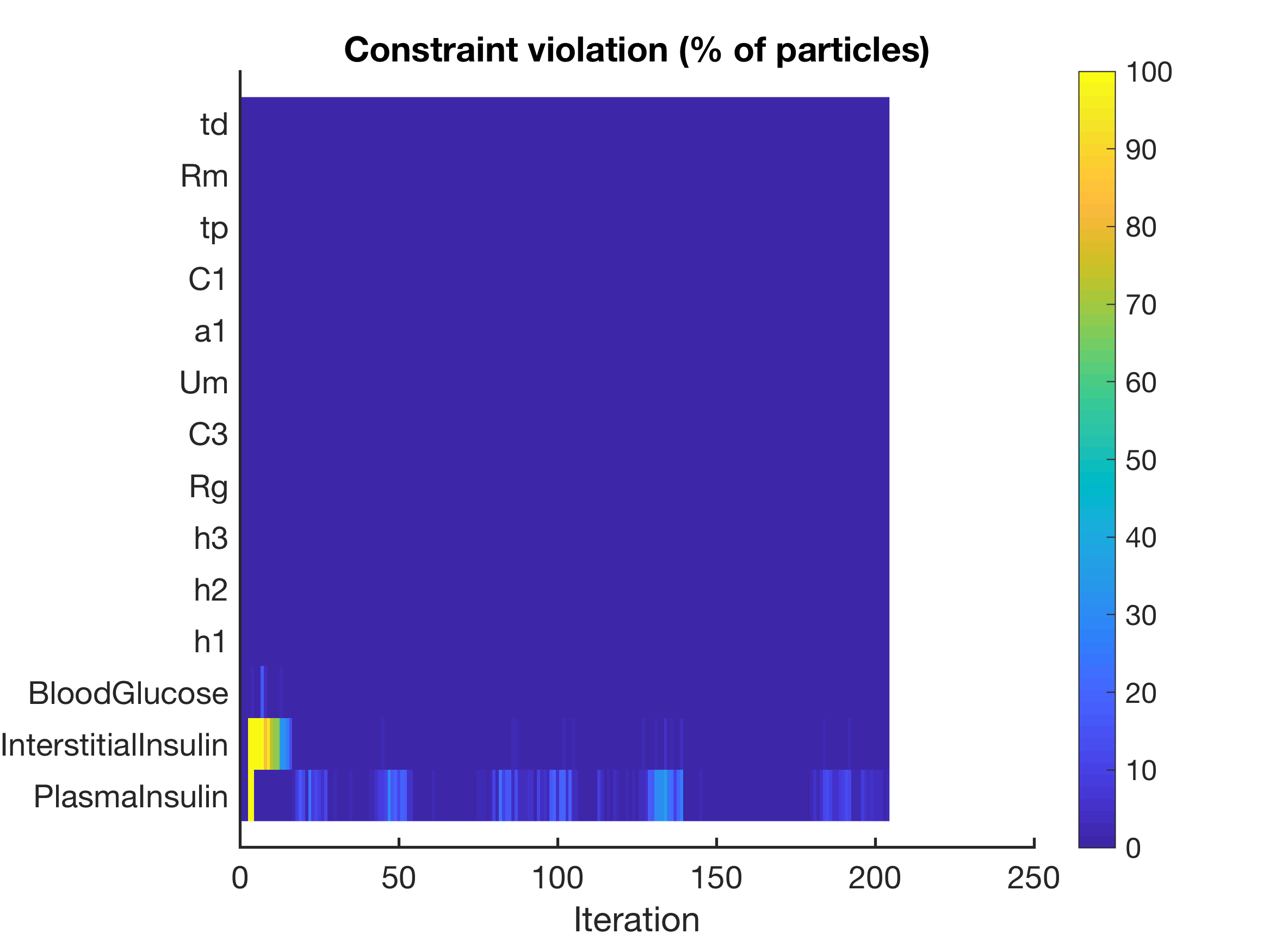}
		\caption{}
		\label{fig:motivation2d}
	\end{subfigure}%
	\begin{subfigure}{0.5\textwidth}
		\includegraphics[width=0.9\linewidth]{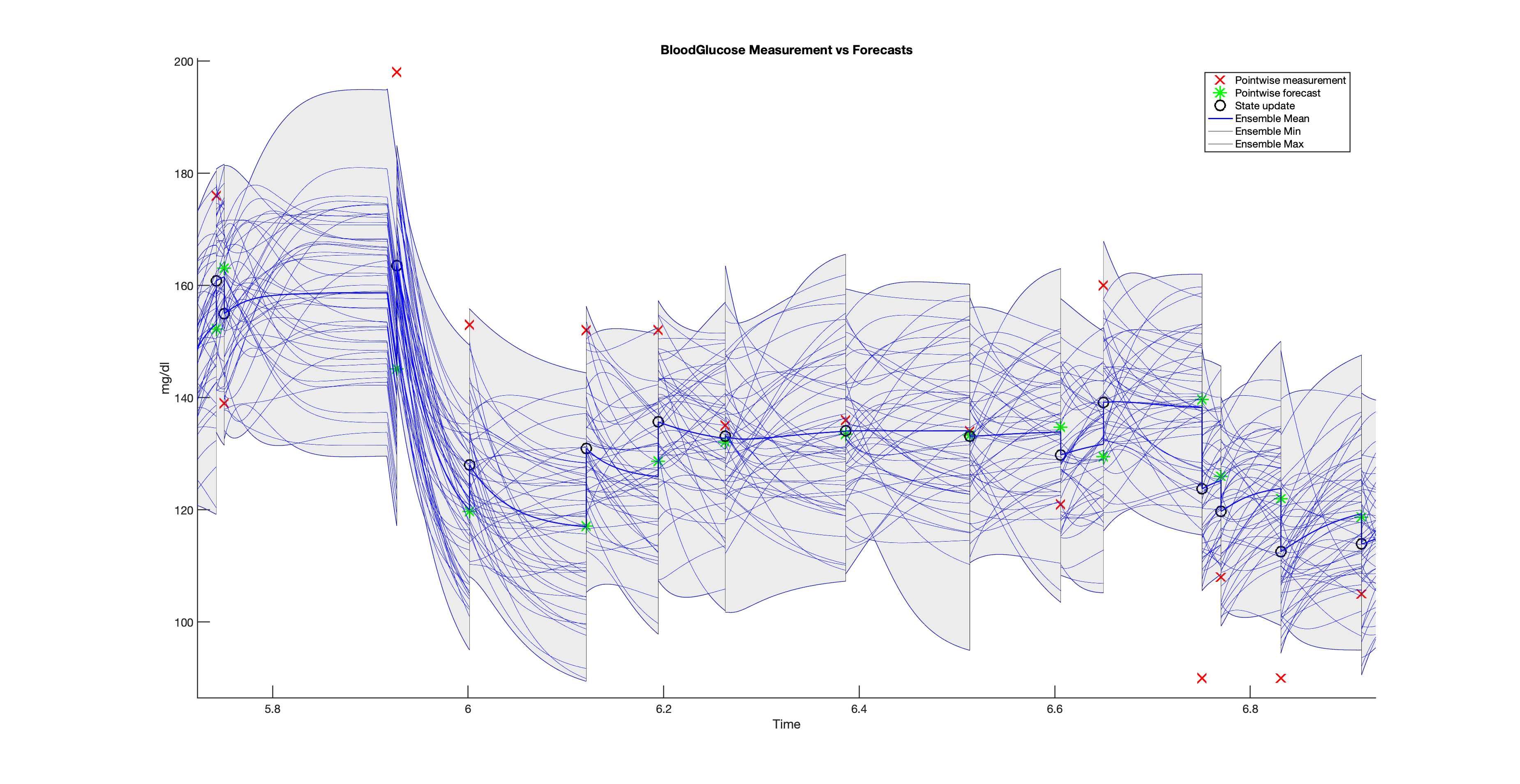}
		\caption{}
		\label{fig:motivation2e}
	\end{subfigure}

	\caption{Patient $593$'s glucose trajectory (a) and insulin states (b) estimated with the constrained EnKF, the constrained parameter estimate trajectory (c), the percentage of particles violating the constraints per \emph{data point} for estimated model states and parameters (d), and the individual ensemble particle trajectories for  $593$'s estimated glucose trajectory (e).}
	\label{fig:motivation2}
\end{figure}

One approach toward solving the first problem, identifiability failure due to unmeasured model states, is to constrain the unmeasured states to occupy plausible ranges by constraining the EnKF.  We developed a constrained version of the EnKF \cite{cenkf} allowing us to enforce constraints on the EnKF algorithm. This constrained EnKF (CEnKF) algorithm uses the standard EnKF framework, checks after the update step if any of ensemble particles lie outside the predefined constraint region, and if any particles lie outside the constraint region the algorithm invokes a quadratic program to solve for parameter values that bring the particles inside the constrained region.  The results of applying this method can be seen in Fig \ref{fig:motivation2}a. The model converges within about $1.5$ days compared to the $4.5$ days in the unconstrained case, a substantial two-thirds decrease in time to convergence \cite{ALBERS2023104477}. Comparing Figs. \ref{fig:motivation1} and Fig. \ref{fig:motivation2}b we can see the constraints' impact on the unmeasured insulin, forcing the insulin to take more reasonable values.  Additionally in Fig. \ref{fig:motivation2}c we can observe the parameter trajectory that shows a relatively steep parameter changes within the first $1.5$ days. And finally, in Fig. \ref{fig:motivation2}d we can see the proportion of particles that violate the constraints for each data point or iteration of the CEnKF.  Within about $20$ data points, the constraints are usually satisfied, and it was only the insulin values that were violating the constraints. However, imposing the constraints brings unintended consequences.  Figure \ref{fig:motivation2}e shows the particle trajectories within the EnKF ensemble, and the variation of the mean. It appears that the constraints have changed the model dynamics from oscillations to fixed points.  Meaning, the EnKF variance is mostly constructed from variations in the fixed point location according to the variation within the ensemble in parameter space.  \emph{The end result of this demonstrates that constraining the model states (a) does help solve the identifiability problem, (b) is not enough to fully solve the initialization problem and (c) may have unintended consequences related to the model dynamics.} Specifically, the optimal solution tends toward a fixed point of the dynamical system, a solution we know to be wrong, implying that the constrained model dynamics are not representative of the underlying system we know to exist.
 
\begin{figure}[t!]
	\centering
	\begin{subfigure}{0.5\textwidth}
		\includegraphics[width=\linewidth]{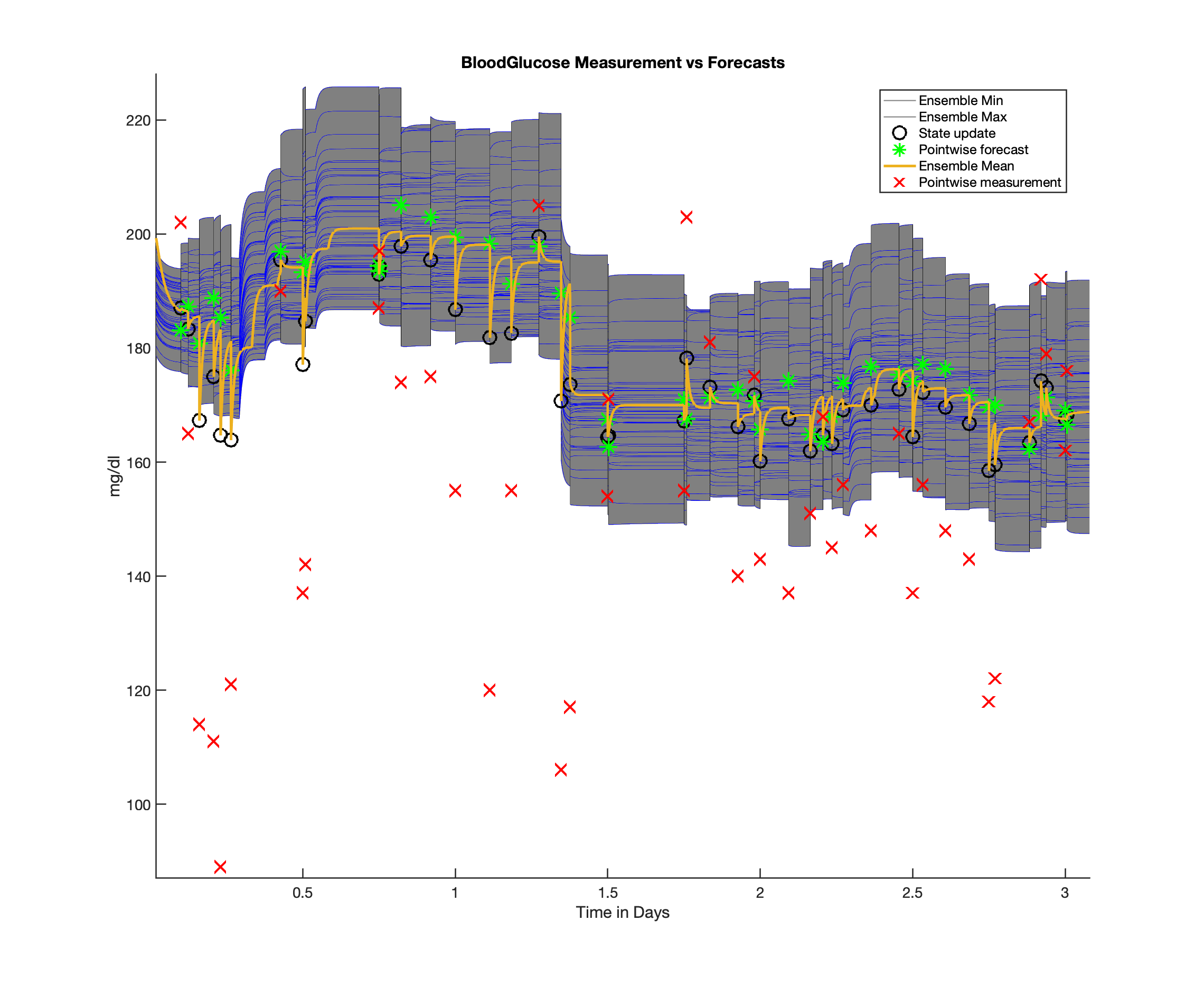}
		\caption{}
		\label{fig:motivation3a}
	\end{subfigure}%
	\begin{subfigure}{0.5\textwidth}
		\includegraphics[width=\linewidth]{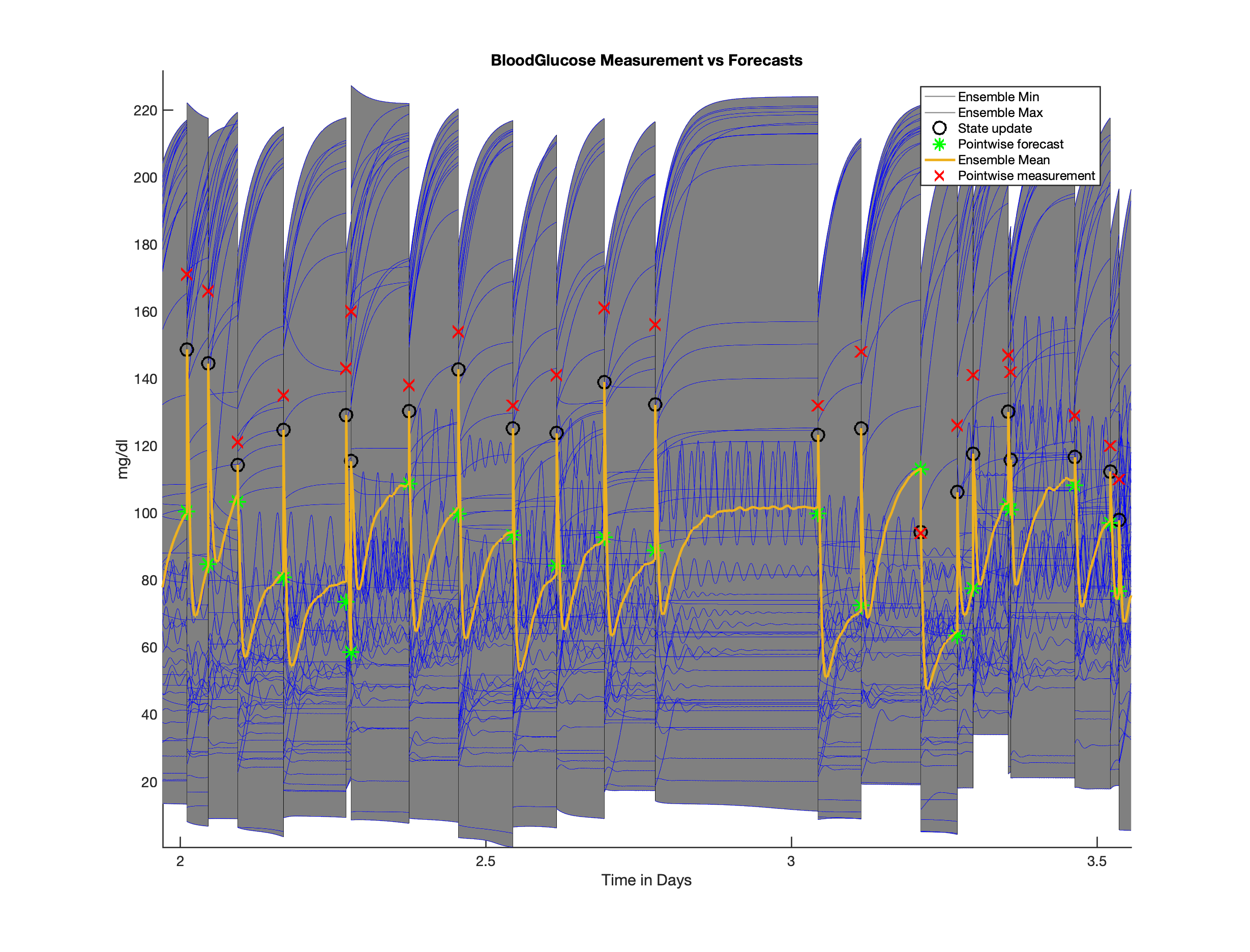}
		\caption{}
		\label{fig:motivation3b}
	\end{subfigure}

	\caption{Constrained EnKF forecasts using parameter initializations computed by optimizing over the patient's first $24$ hours of data. Plot (a) shows patient $593$'s individual ensemble particle trajectories of the estimated glucose trajectory,  (b) shows patient $426$'s individual ensemble particle trajectories of the estimated glucose trajectory. Note the particle ensembles are taken over estimated states and parameters. For patient $593$ the glycemic dynamics are wrong, they should be oscillatory instead of a strongly attracting fixed point. For patient $426$, the glycemic dynamics are representative of a bifurcation point demonstrating several topologically distinct orbits within the same CEnKF ensemble.}
	\label{fig:motivation3} 
\end{figure}

\emph{A naive idea we thought might work to solve the parameter initialization problem} was to run deterministic optimization, bootstrapping over a few thousand random initial conditions of the optimization using the first $24$ hours of data. Then, the outcome from this calculation would be used to initialize parameters and states for the EnKF. The results of this effort for both patients $426$ and $593$ are shown in Fig. \ref{fig:motivation3}. This initialization indeed reduces the MSE between the EnKF ensemble mean and these data. However, for patient $593$, shown in Fig. \ref{fig:motivation3}a, the oscillatory glycemic dynamics are completely gone and have been replaced with strongly attracting fixed points. We know these parameters estimates are wrong. Additionally, if we apply the same parameter initialization procedure to patient $426$'s data, we find something even more interesting.  The dynamics of the particles inside the ensemble for the model initialized with these optimized parameters, shown in Fig. \ref{fig:motivation3}b, have several simultaneously present, topologically different orbits ranging from fixed points with different attraction rates to periodic orbits. Meaning, \emph{the naive constrained DA did further minimize MSE, but it led to error minimizing parameters at a bifurcation point.} In hindsight, the algorithms are doing exactly what we would suspect given the sparse data and the limited constraints we have imposed.  A bifurcation point is a particularly useful location in parameter space from the functional approximation theory perspective because a bifurcation point provides a diverse set of functions to use to estimate the mean and variance of data. 


To address these problems we need to include more structure in the loss function beyond minimizing the MSE between an ensemble average and data subject to constraints. We want to minimize the MSE between an ensemble average and data while also \emph{(a)} ensuring that other global features such as measure-theoretic and topological properties are also accurately characterized and \emph{(b)}  forcing the estimation process to respect the model. This requires a more diverse loss function structure, and establishing how to formulate and estimate such a loss function is the point of this work.

\section{A new approach to DA smoothing and initialization}

We propose an objective function $L$ built as a weighted sum of components $L_k$, each serving a different purpose. Recall that $\{x^j\}\}$ represent the observable variables, $\{z^j\}\}$ represent the unobserved, latent variables, and $\{y^j\}\}$ represent the observations. Then,
\begin{enumerate}
\item[$L_1$] quantifies the point-wise agreement between the $\{x^j\}$ and the $\{y^j\}$. 

\item[$L_2$] quantifies distribution-wise agreement between the $\{x^j\}$ and the $\{y^j\}$, in terms of a slowly modulated invariant measure. $L_1$ and $L_2$ are the only components of $L$ involving $\{y^j\}$.

\item[$L_3$] and $L_4$ quantify the agreement of variables and parameters with the model provided, with two conditional distributions, one forecasting $\{x^j\}$ and the other the $\{z^j\}$ (In other words, in our model --as in any Gaussian model-- $x$ and $z$ are conditionally independent given the prior state and the value of the parameters.)

\item[$L_{k}$], with $k = 4+l$ and $l = 1,\dots n_p$, penalizes variations of the parameter $\alpha_l$ over time, where $n_p$ is the number of parameters estimated over their respective estimation ranges.

\end{enumerate}
The methodology is described in the following subsections.

\subsection{Defining data with measurement functions}

Data are accumulated by measuring the underlying system at times $\{t^j\}$ where $j$ indexes the ordering of these data. Measurements $\{y^j\}$ are taken according to the measurement function $h$ operating on the unknowable state of the system at time $t^j$, denoted by $\mathcal{Y}^{t^j}$, according to:
\begin{equation}
y^{j} = y\left(t^j\right) = h\left(\mathcal{Y}^{t^j}\right) + \eta^{t^j},
\end{equation}
where $\eta$ represents noise. More generally, a measurement $y$ is drawn from a probability distribution
$$ y^j \sim \rho^h\left(y|\mathcal{Y}^{t^j}, t^j\right). $$

These measurement functions, $h$ or $\rho^h$, can have dynamics of their own and be represented by dynamical systems or stochastic processes. In biomedicine broadly and biomedical informatics specifically,  the measurement process are part of the health care process \cite{jamia_phys_ehr}.  Human health can be non-stationary on multiple time scales, and patients are often only measured for a reason, but sometimes are measured according to clinical protocols \cite{51_protocols,george_eval}. Such processes can have very complex missingness properties \cite{missingness_stats_book,missingness_can_help_classification}. Additionally, measurement processes can impact effective time parameterizations \cite{time_parameterization_gh_ap_da},  and lead to signals being present due to inadvertently combining differently measured, statistically different processes that originate from the same source \cite{stat_dyn_diurnal_correlation_short}.  The net sum is that it can be important to conceptualize measurement processes as dynamical systems or stochastic processes with potentially complex properties and interdependencies.

\subsection{Point-wise agreement: $L_1$}

We would like our $x^j$ to be close to the $y^j$. Typically, we would translate this into the minimization of a loss function, such as
$$ \frac{1}{n} \sum_j \left(y^j - x^j\right)^2. $$
More generally, we can have a probability
$$ \rho^o(y|x) $$
of observing $y$ under current state $x$. Conceptually, this is the result of composing two probability distributions, one for the actual state $\mathcal{Y}$ given $x$, and the other for the observation $y$ given $\mathcal{Y}$:
$$ \rho^o(y|x) = \int \rho^h(y|\mathcal{Y}) \rho^m(\mathcal{Y}|x) \ d \mathcal{Y}. $$
Then we maximize the log-likelihood of these data:
$$ \max_x L_1 = \frac{1}{n}\sum_j \log\left[\rho^o\left(y^j | x^j \right)\right]. $$
Minimizing the sum of squared differences corresponds to a Gaussian assumption on $\rho^o$. Given measurement processes, this assumption may or may not be valid, but we will not directly address potential related issues here. Even so, this objective function is likely to be too restrictive given the realities of data, as it requires each $x^j$ to be close to the observed $y^j$: some $y^j$ may be outliers that we do not necessarily want our $x^j$ to adjust to and, more generally, the model may have a hard time adjusting so as to pass near all observations. To account for this, we mollify $\rho^o$, writing
$$ \max_x L_1 = \frac{1}{n}\sum_j \log\left[(1-\epsilon)\rho^o\left(y^j | x^j \right) + \epsilon \rho^0\left(y^j\right)\right], $$
where $\epsilon \sim O(\frac{1}{n})$ and
$$ \rho^0\left(y \right) = \frac{1}{n} \sum_l K^y\left(y, y^l\right)  $$
is a kernel-based estimation of the probability of $y$ lacking any accompanying $x$. For our experiments, we have adopted for $K^y$ a Gaussian kernel with bandwidth determined by the rule of thumb.

\subsection{Invariant measure agreement: $L_2$}

A more global constraint on the $\{x^j\}$ is that their invariant measure $\rho^x$ be close to the $\rho^y$ underlying the $\{y^j\}$. 

There is more than one sample-friendly way to enforce the equality of $\rho^x(z)$ and $\rho^y(z)$. A weak formulation for their equality is that
$$ \int F(z) \rho^x(z) \ dz =  \int F(z) \rho^y(z) \ dz $$
must hold for all measurable test functions $F$, with sample-based representation
$$ \sum_j F\left(x^j\right) = \sum_j F\left(y^j\right). $$

Yet we do not really need to try all possible test functions: the choice $F(z) = \rho^x(z) - \rho^y(z)$ yields
$$ \int F(z) \left(\rho^x(z) - \rho^y(z)\right) \ dz  = \int \left(\rho^x(z) - \rho^y(z)\right)^2 \ dz, $$
a non-negative quantity that vanishes if and only if $\rho^x$ and $\rho^y$ agree.
We implement this through a simple kernel density estimation:
$$ \rho^x(z) = \frac{1}{n} \sum_j K^y\left(z, x^j\right), \quad \rho^y(z) = \frac{1}{n} \sum_j K^y\left(z, y^j\right), $$
or, in order to allow the measures to evolve slowly over time, 
$$ \rho^x\left(z | t \right) = \frac{1}{\sum_l K^t\left(t, t^l\right)} \sum_l K^y\left(z, x^l\right) K^t\left(t, t^l\right),  $$
$$ \rho^y\left(z | t \right) = \frac{1}{\sum_l K^t\left(t, t^l\right)} \sum_l K^y\left(z, y^l\right) K^t\left(t, t^l\right),  $$
where the kernel $K^t$ has a bandwidth $T_l$ to be determined below.
Then we propose
$$ L_2 = -\frac{1}{n} \sum_j \left(F\left(x^j, t^j\right) - F\left(y^j, t^j\right)\right), \quad F(z, t) = \rho^x(z | t) - \rho^y(z | t). $$

\subsection{Agreement with the model and slow modulation: $L_3$, $L_4$, $L_{4+l}$}

The maximizations above were carried out over the $\{x^j\}$ alone, with no reference made so far to the model's parameters and its latent variables $\{z^j\}$. In order to account for these, we add the likelihood of the transition probabilities
$$ L_3 = \frac{1}{n}\sum_j \log\left[\rho^x\left(x^j | x^{j-1}, z^{j-1}, \Delta t^{j}, \alpha^{j} \right)\right], $$
$$ L_4 = \frac{1}{n} \sum_j \log\left[\rho^z\left(z^j | x^{j-1}, z^{j-1}, \Delta t^j, \alpha^{j} \right)\right]. $$
Here
$$ \Delta t^n = t^n - t^{n-1} $$
and we assume that the model's transitional probabilities $\rho^x$ and $\rho^z$, depending on parameters $\alpha$, are given in a closed form. For the examples below, we have adopted a general model for oscillatory behavior, described in the next subsection.

In order to allow the model's parameters $\alpha_l$ to evolve slowly over time, i.e., flex within an estimation time window, we introduce a transitional probability for them too:
$$ L_{4+l} = \frac{1}{n} \sum_j \log\left[\rho^l\left(\alpha^j_l | \alpha^{j-1}_l, 
\Delta t^j\right) \right], $$
for which we propose
$$ \alpha_l^{j+1} \sim N\left[d_s^j \alpha_l^j + \left(1-d_l^j\right) \tilde{\alpha}_l , \left(1-d_l^j\right) {\sigma_l}^2 \right]. $$ 
Here $\tilde{\alpha}_l$ is a prior value of $\alpha_l$ we would like the model to relax to when measurements take long to arrive --this parameter can be either externally provided, estimated or made part of the optimization--, $\sigma_l$ is the uncertainty in the parameter, and $d_l^n$ is a decaying weight starting at $1$ and ending up at $0$:
$$ d_l^j = e^{-\frac{\Delta t^{j}}{T_l}}, $$
%
%
with the same long time-scale $T_l$ adopted for the modulation of the invariant measure.


\subsection{Canonical oscillatory model}


In section \ref{sec:motivation} we posed a  motivating problem of estimating a system that oscillates but that is measured sparsely in time, leading to plausible model solutions that have qualitatively wrong dynamics. This motivational example included two notable assumptions. \emph{First}, we assumed data were sampled sparsely enough that dynamical properties such as the oscillatory frequencies would be difficult to exactly estimate. \emph{Second}, we assumed a highly specific parameterized family of models, an ODE model of glucose-insulin mechanics \cite{sturis_91}, to estimate these oscillatory data. The parameterized nature of mechanistic ordinary, partial, or stochastic differential equation-based models have an underlying rigidity that yields benefits such as interpretability and estimation with relatively little data. Such models  also have drawbacks such as solutions that minimize least squared errors that lie at bifurcation points or produce qualitatively wrong dynamics.


To address these problems, we proposed a new method for estimating parameterized models. This new filtering procedure assumes the availability of a parameterized stochastic model of the dynamics underlying $x$ and $z$. Additionally, the model must be a parameterized family whose transitional probabilities are given in closed form for which we can readily compute derivatives. For the experiments of this article, rather than adopting a highly-specific, complex, and biophysically motivated model as was used in the opening example, we developed a generic model for oscillatory behavior, in the spirit of the simple stochastic system of equations
\begin{eqnarray*}
  dx &=& \left(\omega\ z - \gamma\ x\right) \ dt + \sigma\ dW \\
  dz &=& -\omega\ x \ dt,
\end{eqnarray*}
where $\omega$ models the frequency of the oscillations and $\sigma$ their variability, which together with the damping parameter $\gamma$ determines their mean amplitude. This model can have both oscillatory and fixed point mean dynamics depending on model parameters, and importantly, oscillatory frequencies that are relatively rigid, limited, and parameter-dependent. This model also comes with the additional benefit of allowing us to test how much of the dynamics present in the data we can capture with a generic, simple, but still interpretable model.
Our methodology requires explicit expressions for the distributions of $x^{j+1}, z^{j+1}$ in terms of their values at the previous observation time $t^j$, which could be obtained by solving the system above. Instead, we propose an alternative system, which is given directly by a discrete transition probability model. Since oscillatory behavior is most naturally described in polar coordinates, we introduce the notation
$$ r^n = \sqrt{\left(x^n - b^n\right)^2 + \left(z^n\right)^2}, \quad \theta^n = \arg\tan\left(\frac{z^n}{x^n - b^n}\right), $$
where $b$ represents the local mean around which $x$ oscillates. Then we propose the following transition probabilities:
$$ x^{n+1} \sim N\left[b^{n+1} + r_+^{n+1} \cos\left(\theta^{n+1} \right), \sigma^2 \right], $$
$$ z^{n+1} \sim N\left[r_+^{n+1} \sin\left(\theta^{n+1} \right), \sigma^2 \right], $$
where
$$ r_+^{n+1} =  \left(1 - d_s^n\right) a^{n+1} + d_s^{n}\ r^n, \quad d_s^n = e^{-\frac{\Delta t^n}{T_s}}$$
%
models relaxation of $r$ toward the local oscillation amplitude $a$, with a time-scale $T_s$ typically shorter than the modulation scale $T_l$, and
$$ \theta_+^{n+1} = \theta^n + \omega^n \Delta t^{n}. $$
models a linear evolution of the phase $\theta$ at the local frequency $\omega$.

Thus we pose two time-scales: $T_s$ and $T_l$, for ``short'' and ``long.'' The first is a component of the model, quantifying how fast the current state is forgotten. The second is the time-scale over which the model itself changes, through the slow modulation of its parameters. 

\subsection{Extension to external drivers or kicks}

Next we extend the methodology to systems that have a driver, such as nutrition consumption driving the evolution of blood glucose, that are not explicitly incorporated into the dynamical model. We consider drivers that can be modeled as acting instantaneously, through``kicks'' at known times $k_j$ with specified intensities $I_j$. Rather than modeling the effect of the kick on the system directly, we will think of it as a break that partially decouples the system's state and its ruling parameters before and after the kick, with the strength of the decoupling proportional to the kick's intensity.

Since most of the coupling between successive times in our model depends exponentially of the interval $\Delta t_j$,  increasing this interval to
$$ \Delta t_i \rightarrow \Delta t_i + \alpha I_j $$
whenever the time $k_j$ falls within the interval $\Delta t_i$ would do the job. (We need to take care to keep the original $\Delta t_i$ when they multiply the frequency $\omega$, as this corresponds to the phase of the oscillation evolving over time.) Our choice for the parameter $\alpha$ is
$$ \alpha = \frac{T_s}{\tilde{I}}, $$
where $\tilde{I}$ is a typical intensity, and $T_s$ is the time scale over which things decouple naturally.

This decoupling also affects the kernel $K^t$, where the argument $t_i - t_j$ must be replaced by
$$ \left|t_i - t_j \right| + \alpha \sum_l I_l, $$
where the sum is over all kicks with $k_l$ falling between $t_i$ and $t_j$.

\medskip

A posteriori, we can use the model's resulting jumps in the various variables at the kicking times, to tune a model that predicts the distribution of these jumps as a function of the parameters of the kicks. More generally, this suggests a potential use of our model: to translate the raw data into a coherent evolution of a set of meaningful parameters --typically slow and punctuated by jumps-- that can facilitate further parameterizations of the process under study in terms of more detailed external factors not available, or at least not used by the model.

\subsection{Computational implementation procedure}
\label{sec:comp_proceedure}

Putting together all the ingredients above, we propose the following procedure.
We are given a set of observations 
$$ y^j = y\left(t^j\right), \quad j \in \{0,1, \ldots, n\} $$
and a model for the transitional probability between states
$$ \rho^x\left(x^{j+1} \big| x^{j}, z^{j}, \Delta t^{j+1}, \alpha^{j+1} \right), \quad 
\rho^z\left(z^{j+1} \big| x^{j}, z^{j}, \Delta t^{j+1}, \alpha^{j+1} \right), \quad 
\rho^l\left(\alpha^j_l | \alpha^{j-1}_l, \Delta t^j\right). $$
Here the $\{x^j\}$ are the model's surrogate for the observations $\{y^j\}$, the $\{z^j\}$ a set of unobserved, latent variables, and $\{\alpha^j_l\}$ the model parameters, $l = 1, \dots, n_p$.

We introduce two kernel functions: $K^y(\cdot,\cdot)$, with bandwidth tuned to the observed $y$'s, and $K^t(\cdot,\cdot)$, with bandwidth $T_l$ associated to a slow modulation of the model's parameters and of the invariant measure. 

We initialize $ x^j= y^j$, $z^j = 0$ and the $\{\alpha_l^j\}$ to default model parameters or initial estimations thereof, through a procedure described below. Then we perform gradient ascent over the $\{x^j\}$,  $\{z^j\}$ and $\{\alpha_l^j\}$ of an objective function $L$:
$$ x_{k+1}^j = x_k^j + \eta \frac{\partial L}{\partial x^j}, $$
$$ z_{k+1}^j = z_k^j + \eta \frac{\partial L}{\partial z^j}, $$
$$ {\alpha_l^j}_{k+1} = {\alpha_l^j}_k + \eta \frac{\partial L}{\partial \alpha_l^j}, $$
with a learning rate $\eta$ that can be evolved adaptively or determined by line search at each step.

Here
$$ L = \sum_{k=1}^{4+n_p} \lambda_k L_k, $$
where the weights $\lambda_k$ may vary over algorithmic time. In particular, we set $\lambda_3$, $\lambda_4$ and the $\{\lambda_{4+l}\}$ intially to zero, so that initially the values of the $\{x^j\}$ are anchored at the corresponding $\{y^j\}$ while the $\{z^j\}$ evolve, and only later the $\{x^j\}$  and $\{\alpha_l^j\}$ start to drift toward their optimal values. We have
$$ L_1 = \frac{1}{n} \sum_j \log\left[(1-\epsilon) K^y\left(y^j , x^j \right) + \epsilon \frac{1}{n} \sum_{i=1}^n K^y\left(y^j, y^i\right)\right], $$ 
enforcing point-wise agreement,
$$
L_2 =  -\frac{1}{n}\sum_{i,j} \left[\left(K^y\left(x^i , x^j \right) - K^y\left(y^i , x^j \right)\right) - \left(K^y\left(x^i , y^j \right) - K^y\left(y^i , y^j \right)\right)\right] \frac{K^t\left(t^i, t^j\right)}{\sum_l K^t\left(t^i, t^l\right)} 
$$
enforcing agreement in measure, and
$$ L_3 = \frac{1}{n} \sum_j \log\left[\rho^x\left(x^j | x^{j-1}, z^{j-1}, \Delta t^j, \alpha^j \right)\right],$$ 
$$ L_4 = \frac{1}{n} \sum_j \log\left[\rho^z\left(z^j | x^{j-1}, z^{j-1}, \Delta t^j, \alpha^j \right)\right]$$ 
and
$$ L_{4+l} = \frac{1}{n} \sum_j \log\left[\rho^l\left(\alpha^j_l | \alpha^{j-1}_l, \Delta t^j\right) \right], $$
enforcing agreement with the model while allowing the parameter values to vary a limited degree over the time interval defined by the times $t^j \in [t^1, t^n]$.

Let us further detail this procedure for our motivating example:

\begin{enumerate}

\item {\bf Choice of kernels.} We adopt the simplest default option: Gaussian kernels, with bandwidth $T_l$ for $K^t$, and determined for $K^y$ by the rule of thumb applied to the $\{y^j\}$:
$$ K^y\left(x_1, x_2\right) = \frac{1}{\sqrt{2\pi} h} e^{-\frac{\left(x_2 - x_1\right)^2}{2 h^2}},
\quad h = \frac{\sigma}{n^{\frac{1}{5}}}, \quad \sigma = \sqrt{\frac{1}{n} \sum_j \left(y^j - \bar{y}\right)^2}, $$ 
$$ K^t\left(t_1, t_2\right) = \frac{1}{\sqrt{2\pi} T_l} e^{-\frac{\left(t_2 - t_1\right)^2}{2 T_l^2}}. $$

\item {\bf Preliminary calculations and initialization} 

We can pre-compute the pairwise kernels
$$ K^y_{ij} = K^y\left(y^i, y^j\right), \quad K^t_{ij} = K^t\left(t^i, t^j\right), $$
the reference densities
$$ \rho_0^i = \frac{1}{n} \sum_j K^y_{ij}, \quad \mu^i = \frac{1}{n} \sum_j K^t_{ij}, $$
and the factors for exponential decay,
$$ d_l^n = e^{-\frac{\Delta t^{n}}{T_l}}, \quad d_s^n = e^{-\frac{\Delta t^{n}}{T_s}}, $$
with time-scales to be determined below. We initialize the states and parameters as follows:
\begin{enumerate}

  \item Set $x^i = y^i$ and $z_i = 0$.
  
  \item Set $ \tilde{b} = \frac{1}{n} \sum_i y^i$, and initialize $b^i$ through kernel regression:
  $$ b^i = \frac{\sum_j y^j K^t\left(t^i, t^j\right)}{\sum_j K^t\left(t^i, t^j\right)}, $$
  settting $\sigma_b$ to the standard deviation of the $\{y^i\}$.
  
  \item Correspondingly, set
  $$ a^i = \frac{\sum_j \hat{a}^j K^t\left(t^i, t^j\right)}{\sum_j K^t\left(t^i, t^j\right)}, \quad  \hat{a}^i = \max_{|t^j-t^i| < T} \left|y^j - b^j\right|, $$
  and $\sigma_a = \sigma_b$.
  Set either $\tilde{a} = 0$, so that oscillations die out in the absence of observations, or to the mean of the $\hat{a}^i$.
  
  \item Regarding $\omega^i$, we first compute $\tilde{\omega}^j = \frac{\pi}{P^j}$, where $P^j$ is the time interval between two sign changes of $y^i - b^i$ that includes $t^j$, and then propose
$$ \omega^i = \frac{\sum_j \tilde{\omega}^j K^t\left(t^i, t^j\right)}{\sum_j K^t\left(t^i, t^j\right)}. $$

Set $\tilde{\omega}$ to the average of the $\omega^i$, and $\sigma_{\omega} = {\tilde{\omega}}$.

\item The natural unit for the time-scales $T_s$ and $T_l$ is the average period $P = \frac{2\pi}{\omega}$ of the oscillations. We will adopt
$$ T_s = \frac{2\pi}{\tilde{\omega}}, \quad T_l = 4 \frac{2\pi}{\tilde{\omega}}. $$

\item The parameter $\epsilon \ll 1$ to account for outliers is somewhat arbitrary. We will pick 
$$ \epsilon = 0.1. $$

\item For the standard deviations $\sigma$ and $\sigma_0$, we adopt
$$ \sigma = \bar{a}, \quad \sigma_0 = 0.1\ \sigma. $$
%

\end{enumerate}

\item {\bf Order in which to carry the various descent processes.}

After the initialization above, most of the model parameters should be roughly in the neighborhood of their optimal values, since the $\{x_i\}$ should be close to the observations $\{y_i\}$, and the other parameters were roughly tuned to the data.  The exception however are the $\{z_i\}$, set initially to non-informative values. Starting the full descent process at once from this initialization has a number of problems:
\begin{enumerate}

\item The model's parameters will rapidly deviate from their initial values, since the current $\{z_i\}$ are not consistent with them. So will the $\{x_i\}$, unless the $L_1$ and $L_2$ hold enough strength to keep them in place.

\item The variance $\sigma$ of the models need to be set quite high, for otherwise the current $\{z_i\}$ and $\{x_i\}$ will tend to fall far in their tails, unbalancing the descent process.

\end{enumerate}

Thus it makes sense to start descending only over the $\{z_i\}$, using large values for the variances, and only afterwards set the variances to better tuned values and switch in the other descent process. In this spirit we propose:

\begin{enumerate}

\item First descend only $L_3$, then also $L_4$, both only over the $\{z_i\}$, with standard deviation
$$ \sigma = 2 \bar{a}. $$

\item Then reset the standard deviation to
$$ \sigma = \bar{a} $$
and descend the full objective function over all parameters.

\end{enumerate}

\end{enumerate}

\section{Numerical results}

\subsection{Measurement functions}
\label{sec:meas_func}

Data or measurements are controlled with measurement processes or functions.  A natural conceptualization of a measurement function is as a stochastic process that controls when and how measurements are taken. Here the measurement process is one aspect of the \emph{health care process}\cite{jamia_phys_ehr}. The times at which data are collected can be pre-established, random or correlated to the data or some other process\cite{time_parameterization_gh_ap_da,missingness_stats_book,missingness_can_help_classification}. Because we are concerned with the impacts of sparsity and measurement noise, we will define three distinct measurement functions that specify the measurement times for use in our experiments as:
\begin{enumerate}
\item $h_1$ where measurements are taken as directed by a human according to need or a protocol \cite{george_eval,melike2019simple,melike_g_control,51_protocols}, corresponding to the subject's real, finger-stick or IV measurement times;
\item $h_2$ where measurements are taken at random times where the difference between two consecutive measurements is uniformly distributed over $[ 60,90]$ minutes;
\item $h_3$ where measurements are taken every five minutes, simulating continuous glucose monitor (CGM) measurements. 
\end{enumerate}
When we have CGM data, $h_3$ will be the measurement times defined by the CGM and $h_2$ will be a subsampling of these measurement times. In the case where we do not have CGM data, $h_2$ and $h_3$ will be extracted from the continuously simulated data.

\subsection{Real and simulated data}

To evaluate our method's performance with different generating dynamics and sampling patterns, we will use two different sources of real and simulated glucose-insulin data.  While we chose blood glucose dynamics because estimating them was the initial motivation for this work, we also think that glycemic dynamics are varied and represent many other situations. To add to the generalizability of the context, we will focus on two cases of glycemic dynamics, \emph{(i)} tube-fed patients in the ICU and \emph{(ii)} normal patients in the ``wild''.  The blood glucose levels of intensive care unit (ICU) patients are often oscillatory when the patient is given constant nutrition \cite{sturis_91} and are highly non-stationary because of the effects of critical illness and related interventions. The blood glucose levels of patients in the wild are more akin to a damped-driven oscillatory with noise because, e.g.,  nutrition consumption causes a rapid increase in blood glucose followed by an oscillatory and noisy return to glycemic homeostasis.


\subsubsection{Real world data}
\label{sec:real_data}


The first data set includes ICU data extracted from the Columbia University Medical Center Clinical Data Warehouse for a previous study \cite{melike2019simple,ALBERS2023104477,WANG2023104547,cenkf,melike_g_control}. All the patients in this dataset were fed through an enteral tube and were delivered intravenous (IV) insulin for glycemic management. We selected one patient, patient $593$, as a representative for highly non-stationary glycemic behavior and challenging model estimation. We used this patient's data in Section \ref{sec:motivation} as one of our motivating examples. These data consist of the point-of-care blood glucose measurements, carbohydrate records delivered through the enteral tube, and exogenous IV insulin records. Additional details about the data are presented in Table \ref{tab:1}.

The second data set was collected from a healthy subject representing ``normal'' glycemic dynamics observed in the ``wild'', used in several previous studies \cite{albers_plos_comp_bio_DA_I,jamia_da,levine_albers_plos_comp_bio_DA_II,me_lena_1,lena_reasoning_book_2017,BURGERMASTER2023104419}. The words ``in the wild'' refer to the fact that this individual measured their nutrition and glycemic dynamics while living otherwise normally. We selected this subject's data as a representative of data collected in usual living conditions and because the dynamics are distinctly different. In particular, glycemic dynamics of individuals not being continuously fed resemble a noisy damped, driven oscillator with a noisy glycemic equilibrium.  These data consist of sparse finger stick BG measurements, continuous glucose monitor (CGM) data, and carbohydrate intake of meals and snacks reviewed and verified by nutritionists. We provide additional details in Table \ref{tab:1}. The data sets vary according to different measurement functions.  For measurement function $h_1$ we restrict the model-estimated data to the finger-stick measurements taken by the patient while data corresponding to $h_2$ are randomly down-sampled from the CGM data, and data corresponding to $h_3$ are the raw CGM data.

\subsubsection{Simulated data}
\label{sec:sim_data}

In addition to the real data, we used simulated data to evaluate this methodology's performance for the ICU patient because this patient did not have CGM data and we evaluate the estimation algorithm assuming different data sampling patterns via the different measurement functions.  Additionally, simulated data allow us to have a completely knowable ground truth for at least one of our evaluative cases, even if the model's representation of these data is not ideal.  






\paragraph{High level computational workflow for generating simulated data}  We estimate glucose by fitting plasma glucose and three parameters ($t_p$, $R_g$, $a_1$) of the Ultradian model \cite{sturis_91} while setting the remaining parameters at their default values. We selected these parameters because they control the large subsystems of the model. Specifically $t_p$ controls insulin clearance, $R_g$ controls insulin resistance, and $a_1$ controls insulin secretion. A detailed description of this model is in the Appendix. For the relatively stationary data from the patient in the wild we used standard Metropolis Markov Chain Monte Carlo (MCMC) with three chains and 10,000 iterations for each chain.  Not every chain converges, so we only use the empirical density functions of the converged chains to sample parameters to create the simulated data. We consider chains converged using the Geweke statistic \cite{geweke,bda3} and minimizing mean squared error, following the methodology presented in \cite{WANG2023104547}.  For the severely nonstationary ICU patient, we used an ensemble Kalman filtering (EnKF) with an ensemble size of $100$ similar to what was used in \cite{ALBERS2023104477}. We used the EnKF because the parameter values move too much over the week-long time window we used to estimate the model. In both cases we then apply the measurement functions $h_1$, $h_2$, and $h_3$ to generate the data sets we estimate with our new method.  
 
\paragraph{ICU data generation} We used one week of  data collected from ICU patient 593 from our opening example. These data are documented in Table \ref{tab:1}. The average parameter values estimated over this period were $t_p = 5.5$, $a_1 = 7.5$, and $R_g = 225$. We then simulated blood glucose data every minute with these parameters and subsampled the data according to $h_1$, $h_2$ and $h_3$.

\paragraph{In-the-wild patient data generation} We did not use simulated data for the ``in the wild'' patient because we have CGM data to compare against the model estimates for different measurement functions. 


\begin{table}[ht]
	\centering
	\begin{tabular}{|p{0.2\linewidth}|c|p{0.2\linewidth}|c|c|} \hline
		\multicolumn{2}{|c|}{\multirow{2}{*}{Healthy subject}} & \multicolumn{3}{|c|}{ICU patients} \\ \cline{3-5}
		\multicolumn{2}{|c|}{} & & patient 426 & patient 593 \\ \hline
		total data recording interval (days) &  \multirow{2}{*}{31.6} & total data recording interval (days) &  \multirow{2}{*}{14.0} &  \multirow{2}{*}{21.2} \\ \hline
		total \# of BG & finger stick: 120 & total \# of BG & \multirow{2}{*}{177} & \multirow{2}{*}{249} \\
		measurements & CGM: 7722 & measurements & & \\ \hline
		measured BG val- & finger stick: 102.5$\pm$17.6 & measured BG val- & \multirow{2}{*}{141.2$\pm$18.0} & \multirow{2}{*}{149.7$\pm$32.0} \\
		ues (mean$\pm$stdev) & CGM: 95.1$\pm$13.8 & ues (mean$\pm$stdev) & & \\ \hline
		total \# of meal/ snack recordings &  \multirow{2}{*}{73} & total amount of recorded tube-fed carbohydrates (g) &  \multirow{3}{*}{1938.9} & \multirow{3}{*}{2139.1} \\ \hline
		\multicolumn{2}{|c|}{} & total amount of delivered IV insulin (unit) &  \multirow{3}{*}{0} &  \multirow{3}{*}{1498.1} \\ \hline
	\end{tabular}
	\caption{Patient data descriptions for data used to test and validate our optimization methodology.}
	\label{tab:1}
\end{table}

\subsection{Primary computational questions and the design of data sets for numerical experiments}

Our method is motivated by two goals in the context of sparse data: \emph{(a)} model estimation that balances point-wise accuracy with the preservation of global properties such as agreement of invariant measures, and \emph{(b)} import model flexibility to manage measurement time errors, non-stationarity and sparsity, by allowing temporal parameter flex over a fixed estimation window.  We claim that our method is able achieve the goals of preserving estimation of the invariant measure while also supporting very accurate point-wise estimation and accounting for external noise, shocks, and errors in measurement times  by \emph{(i)} giving up a little accuracy estimating the model both point-wise and in-distribution and \emph{(ii)} allowing parameters to flex over the estimation interval near data points. To verify that our method is able to achieve these goals we designed a series of computational experiments to address the following questions:
\begin{enumerate}

\item[Q1.] \textbf{can the new method include knowledge of global properties of the system to improve model estimation by allowing estimation of global properties that:}
\begin{enumerate}
\item[a.]  we either know or can be derived from the entirety of given data?
\item[b.] are compatible with reasonable point-wise data estimates?
\item[c.] lead to more reasonable model estimates from a face-validity standpoint?
\end{enumerate}

\item[Q2.] \textbf{does accounting for measurement function dependence and measurement timing errors---specifically measurement density and targeted vs random measurement times---impact:}
\begin{enumerate}
\item[a.] the parameter flex over the estimation window; 
\item[b.] the estimation accuracy of the invariant measures;
\item[c.] the estimation accuracy of the trajectory.
\end{enumerate}

\item[Q3.] \textbf{does allowing enough flex in the model parameter estimation on a time window:}
\begin{enumerate}
\item[a.] allow the model estimation to accommodate nonstationarity within the estimation window?
\item[b.] vary depending on the measurement function?
\item[c.] impact the estimation accuracy of the trajectory?
\end{enumerate}
\end{enumerate}

To address these questions we will estimate two types of dynamics, \emph{(i)} oscillatory dynamics that relax to a periodic orbit when driven and to a fixed point when not driven, and \emph{(ii)} damped driven oscillatory dynamics that relax to a fixed point.  Data will be collected via the measurement functions defined in section \ref{sec:meas_func}.

%
%
%

\subsection{Estimation evaluation}

This paper addresses situations where the parameters that minimize an error such as least squared may not be the model parameters we would most like to use. Additionally, because the loss function minimizes a weighted sum or errors, it will be rare that the optimal parameters our method selects correspond to the parameter set that also minimizes one of those components, such as the least squared error. Because of this, standard methods for evaluation will not suffice.  Instead, we will evaluate the method's usefulness qualitatively.  Specifically, we will consider whether the method is able to achieve our goals such as estimate both the invariant measure of the data while also accurately estimating these data point-wise, qualitatively reproduce the type of dynamics present in the generating system in situations where solutions exist that could have lower least squared error, and can incorporate external knowledge of the underlying dynamics.

\subsection{Incorporating global information (\emph{Q1}) and impacts of different measurement functions on model estimation (\emph{Q2}) }

\begin{figure}
	\centering
	\begin{subfigure}{0.33\textwidth}
		\includegraphics[width=\linewidth]{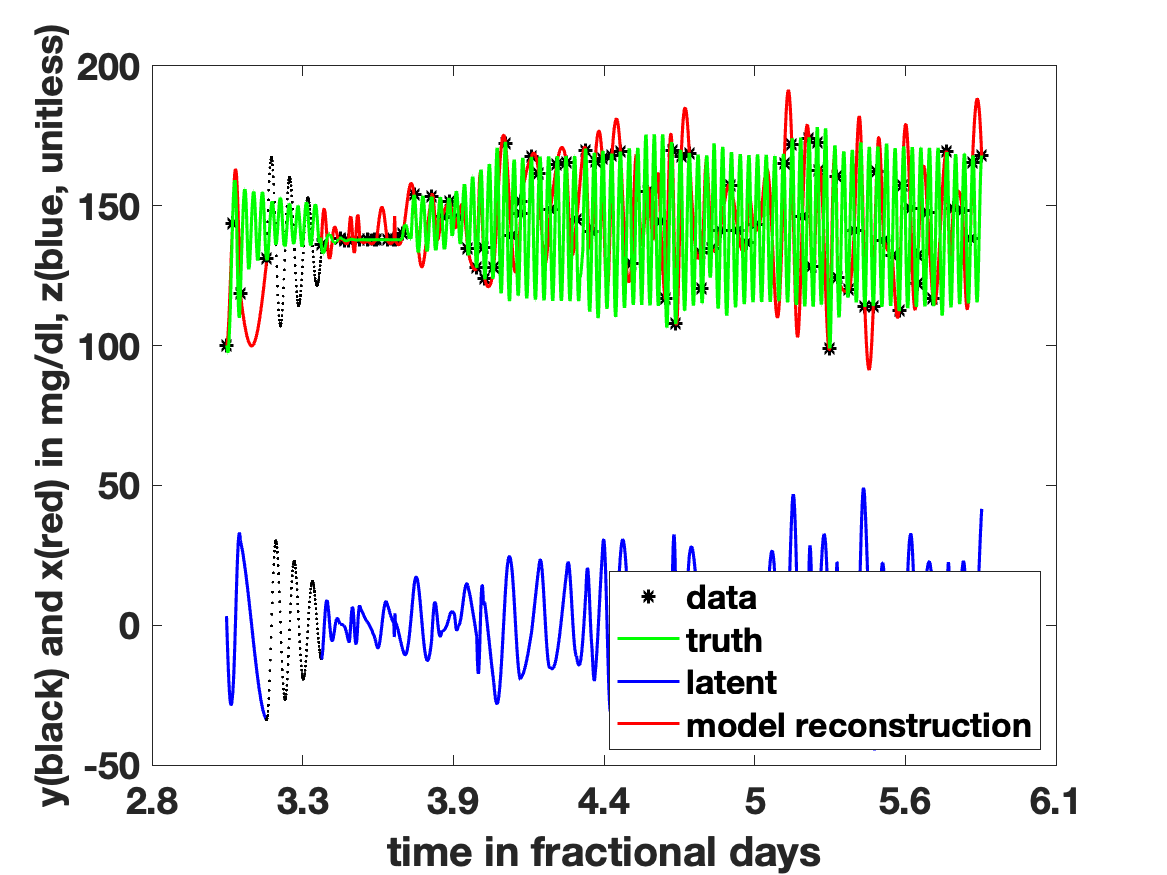}
	\end{subfigure}%
	\begin{subfigure}{0.33\textwidth}
		\includegraphics[width=\linewidth]{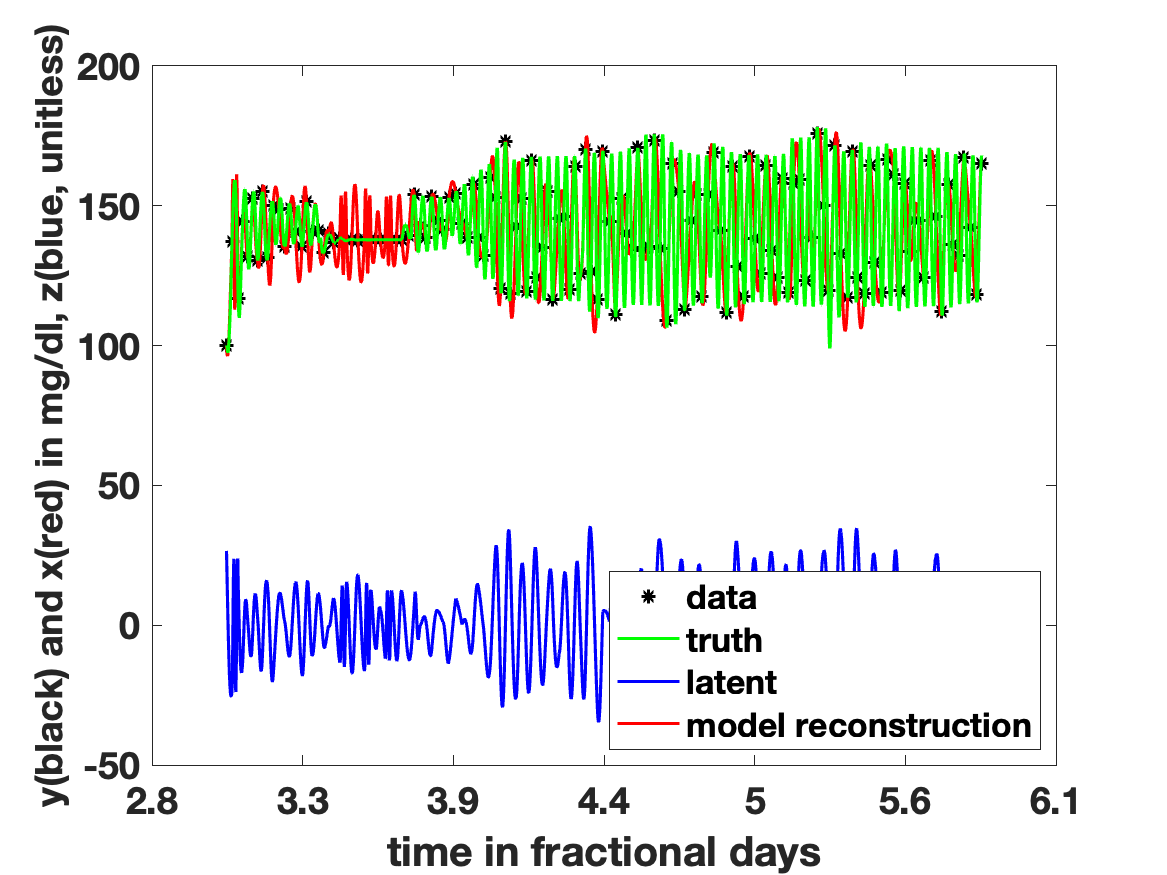}
	\end{subfigure}%
	\begin{subfigure}{0.33\textwidth}
		\includegraphics[width=\linewidth]{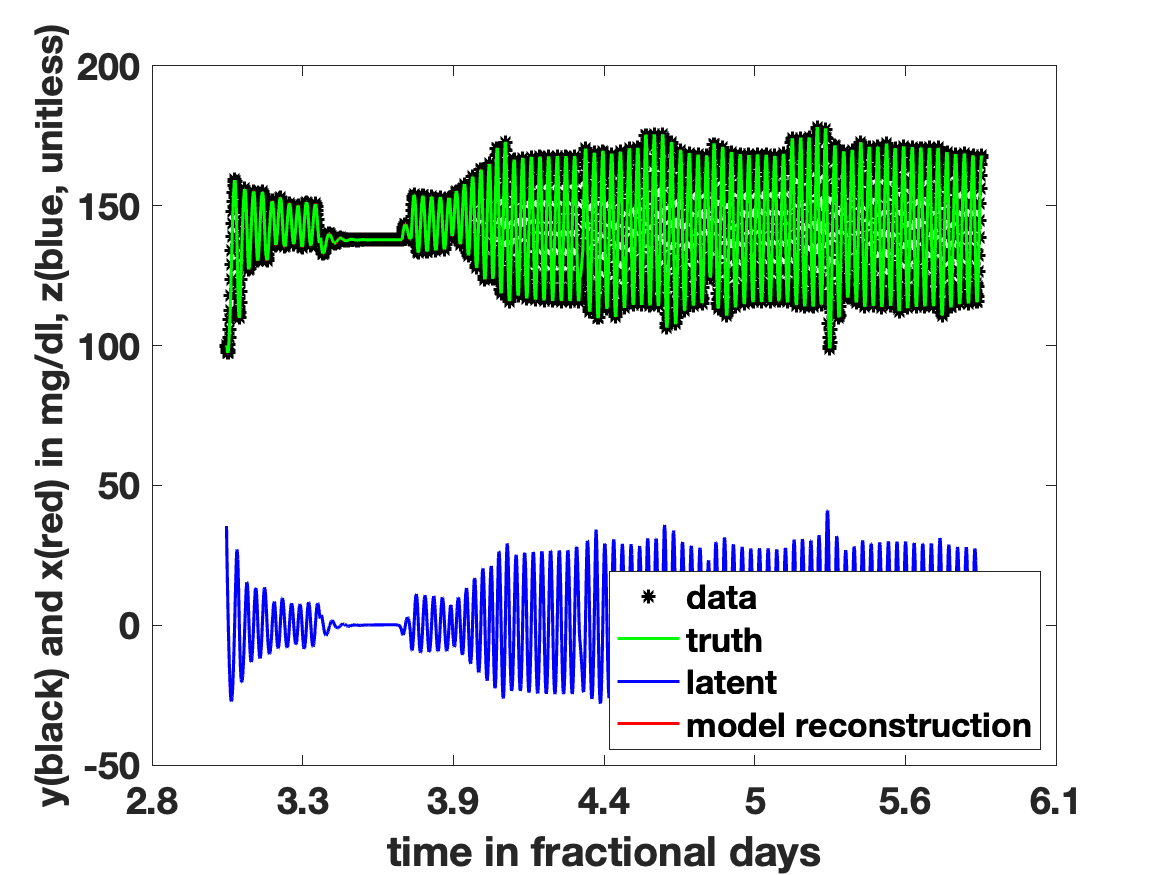}
	\end{subfigure}
	
	\begin{subfigure}{0.33\textwidth}
		\includegraphics[width=\linewidth]{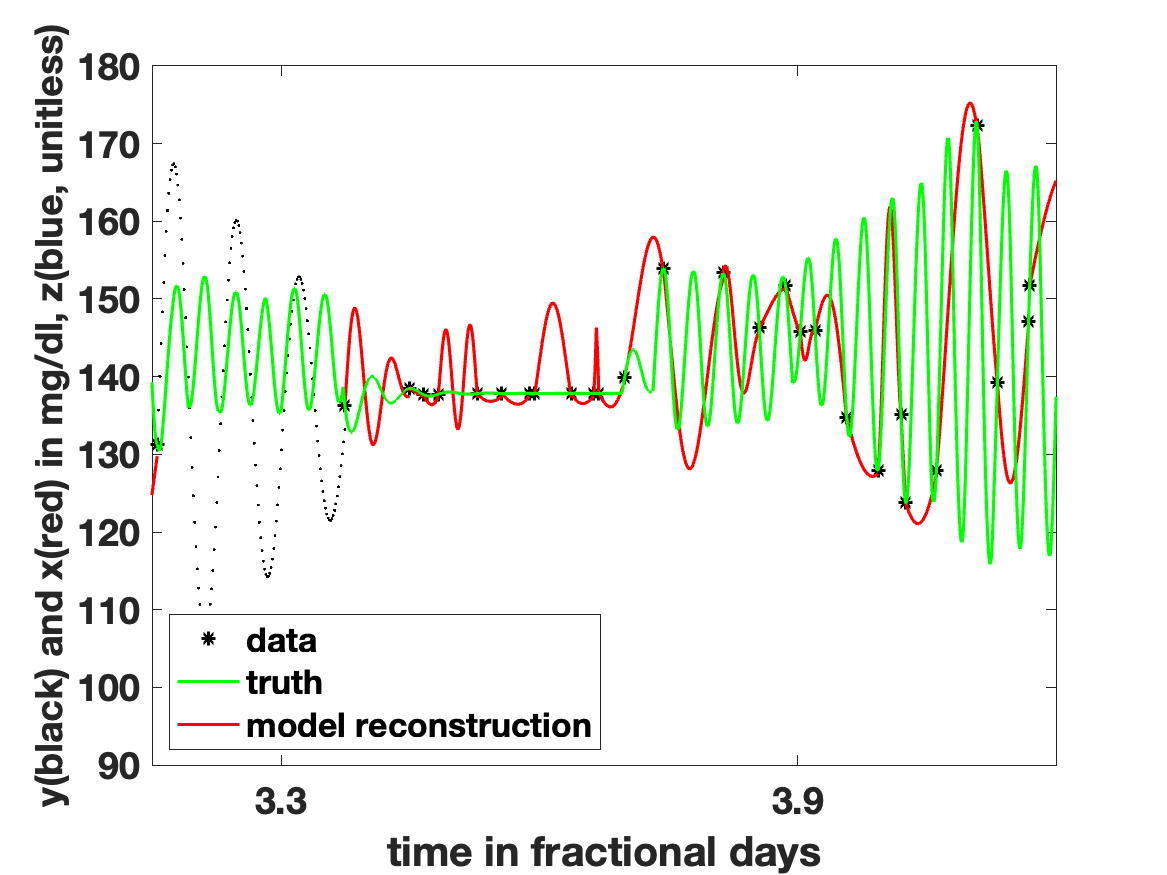}
	\end{subfigure}%
	\begin{subfigure}{0.33\textwidth}
		\includegraphics[width=\linewidth]{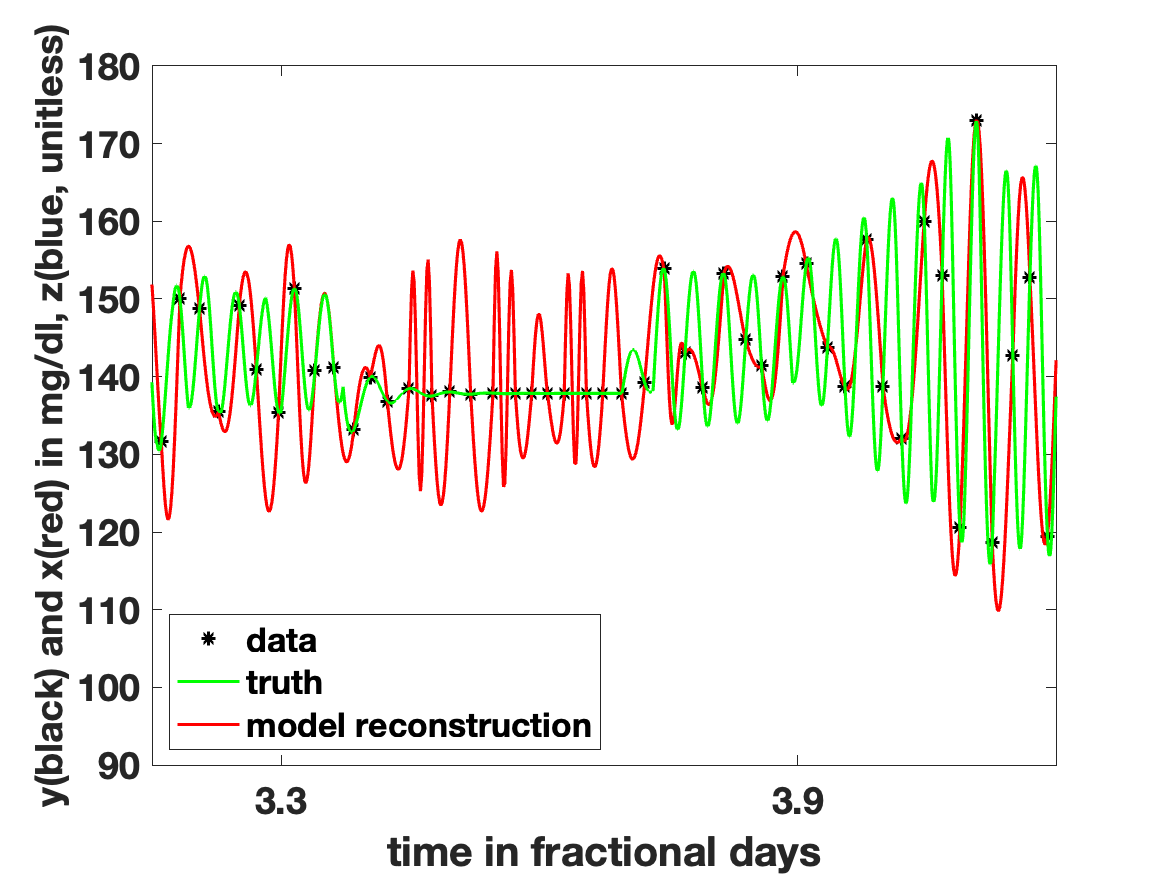}
	\end{subfigure}%
	\begin{subfigure}{0.33\textwidth}
		\includegraphics[width=\linewidth]{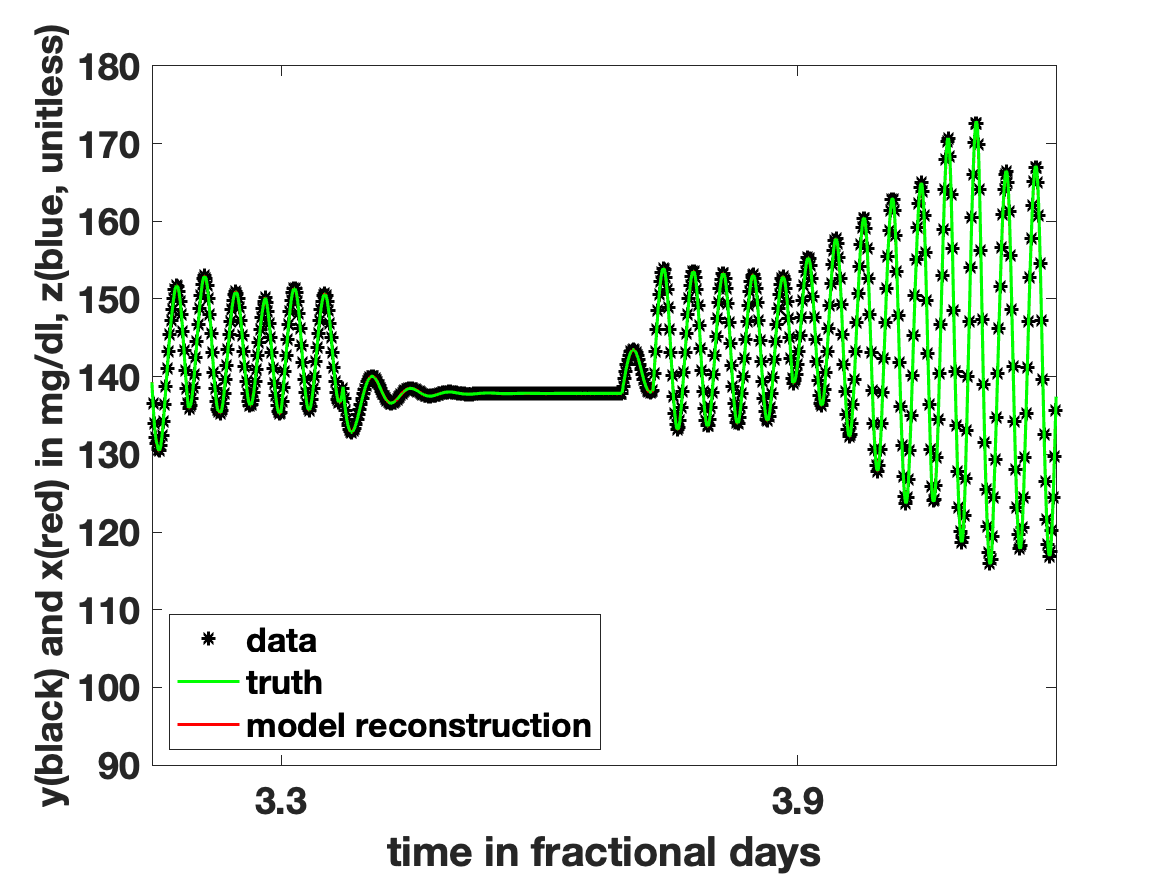}
	\end{subfigure}

	\caption{Given the oscillatory ICU glycemic dynamics and the model basic state assuming oscillatory dynamics, we see the model estimating the simulated glucose measured according to $h_1$ (left), $h_2$ (center) and $h_3$ (right).  We can see that even for the sparse data cases ($h_1$, $h_2$), the model produces oscillatory dynamics with reasonable mean and amplitude while for the densely measured case ($h_3$) the model tracked these data precisely. The point-wise estimates remain accurate in all cases. The dotted lines signal reconstruction for periods without data longer than a prescribed threshold.}
	\label{fig:icu_traj}
\end{figure}

\begin{figure}
	\centering
	\begin{subfigure}{0.33\textwidth}
		\includegraphics[width=\linewidth]{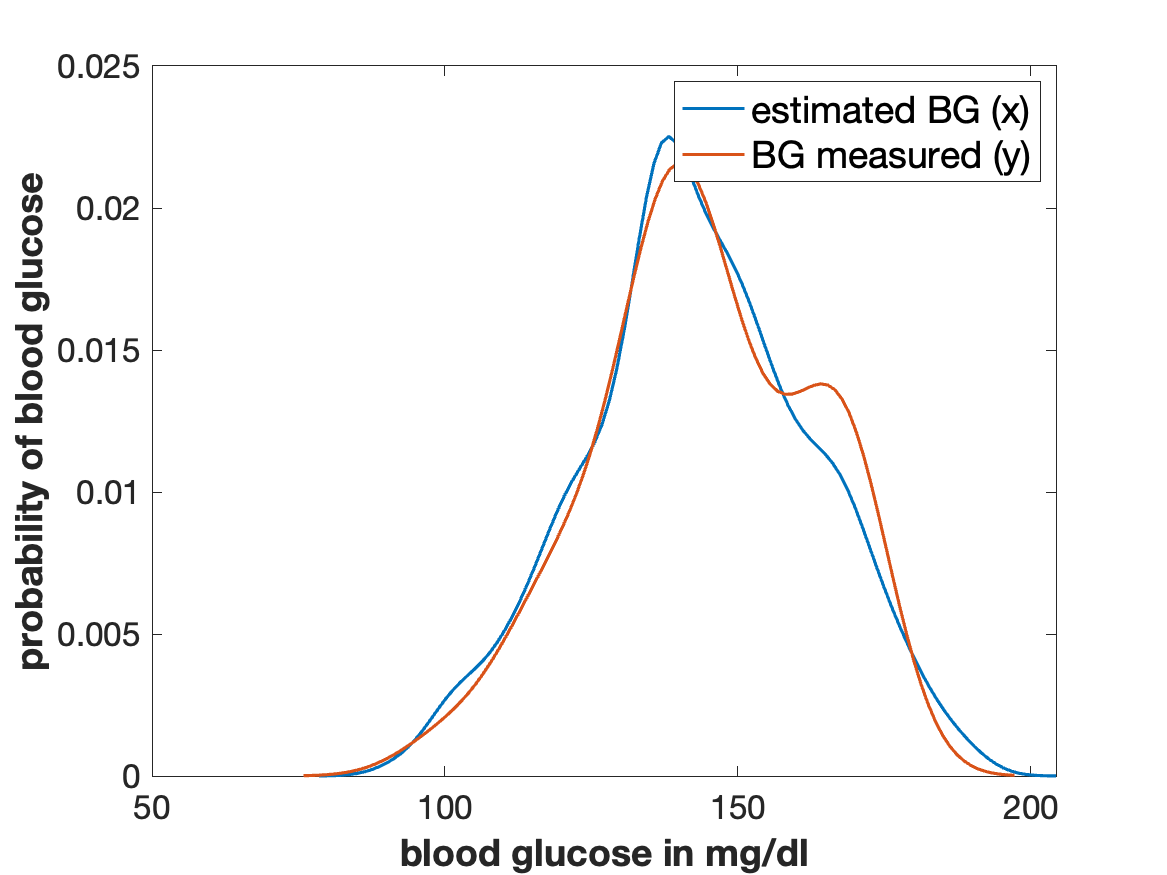}
	\end{subfigure}%
	\begin{subfigure}{0.33\textwidth}
		\includegraphics[width=\linewidth]{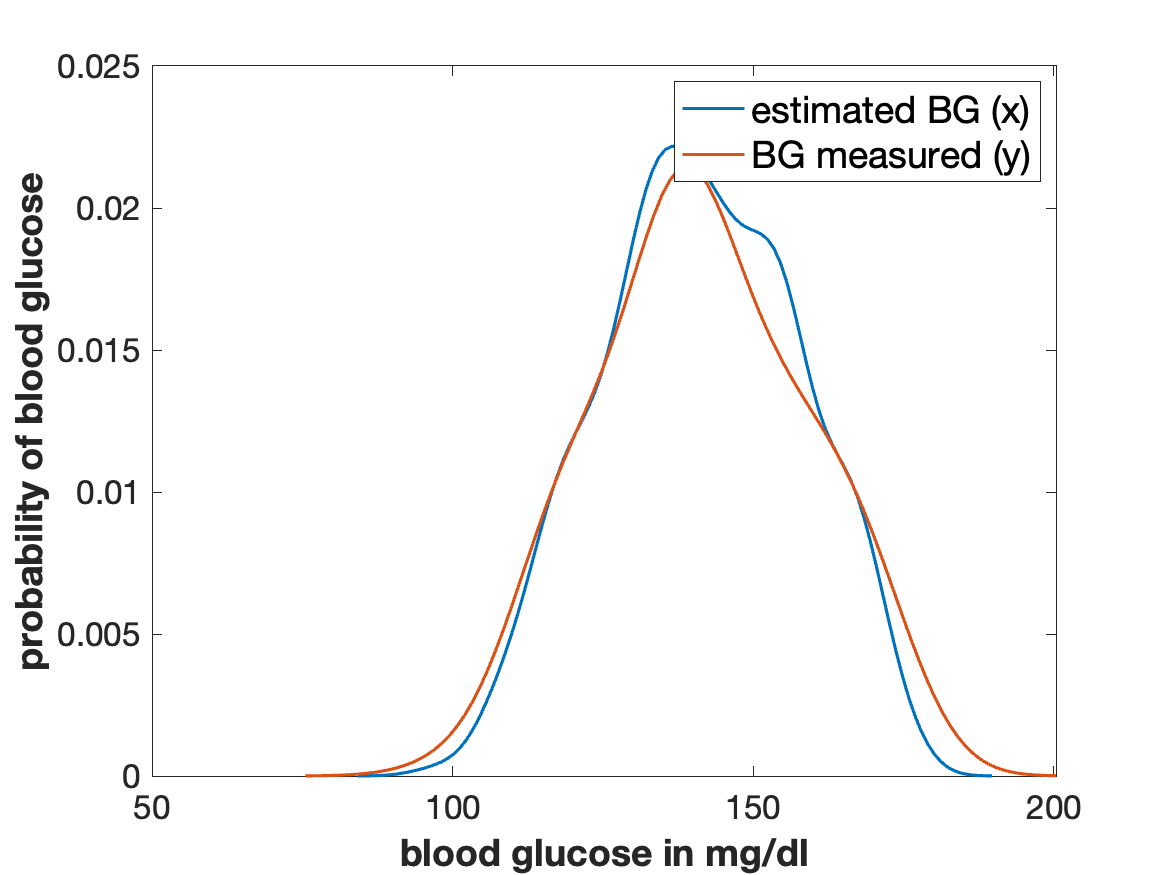}
	\end{subfigure}%
	\begin{subfigure}{0.33\textwidth}
		\includegraphics[width=\linewidth]{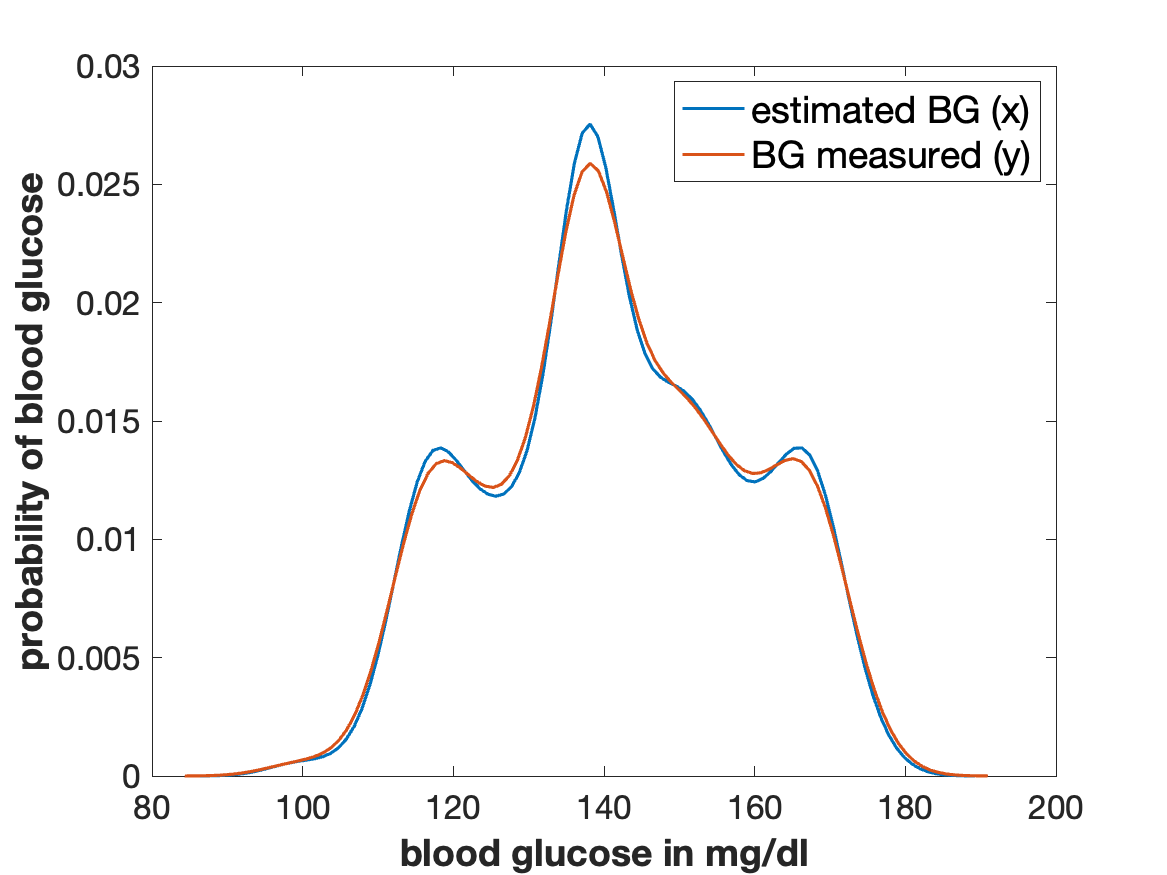}
	\end{subfigure}

	\caption{Given the oscillatory ICU glycemic dynamics and the model basic state assuming oscillatory dynamics, we see the model estimating the simulated glucose measured according to $h_1$ (left), $h_2$ (center) and $h_3$ (right). Note \emph{BG measured} denotes data available to the model when it is estimated and \emph{estimated BG} denotes the model-estimated invariant measure of the data.  The densely measured case ($h_3$) is likely the closest representative of a gold standard baseline, again for data measured frequently in time. We can see that even for the sparse data cases ($h_1$, $h_2$), the model produced an accurate representation of the invariant measure that was not particularly dependent on the sparse measurement function while for the densely measured case ($h_3$) the model estimated all the details of the invariant measure well.}
	\label{fig:icu_pdf}
\end{figure}

The global information we incorporate into the model here includes the invariant measure of the data and the qualitative orbit type (oscillatory versus fixed point, etc.).  We will use the ICU data to show how the new method is able to include global distributional information. We will use the wild case to additionally show how the new method can also be seeded with external knowledge of the underlying dynamics and how this impacts model estimation with sparse data.

Starting with the ICU case, where the generating process produces oscillatory dynamics except for a brief time when the oscillations disappear because the nutrition is turned off, we estimated the model that assumed underlying oscillatory dynamics by setting the relaxed state of the model to have a non-zero amplitude. The estimated model dynamics, shown in Fig. \ref{fig:icu_traj}, reproduced the mean and amplitude of oscillations well while the frequency of oscillations was only correctly represented when the data were measured frequently ($h_3$).   The brief time when the generating processes was not oscillating was only represented with the correct dynamics when the data were measured frequently ($h_3$). Interestingly, the sparse clinician measurement function, $h_1$, reproduced the non-oscillatory dynamics substantially more accurately than the random measurement function, $h_2$, potentially implying information in the clinician-driven measurement times \cite{time_parameterization_gh_ap_da,pop_phys,levine_lagged_regression_AMIA,george_lagged_correlation_jamia,patient_aggregation_paper}. It is likely that a specifically glucose-insulin model whose oscillations are directly controlled by nutrition would have faired differently. See Figure \ref{fig:icu_pdf} where the model-estimated glucose distribution is compared with the data used to estimate the model. The invariant measures of these highly non-Gaussian  data were reproduced well regardless of the measurement function. Together, Figs. \ref{fig:icu_traj} and \ref{fig:icu_pdf} demonstrate that both point-wise and in-distribution estimation were quite accurate, producing affirmative answers to \emph{Q1a}, \emph{Qb} and \emph{Qc}.

\begin{figure}
	\centering
	\begin{subfigure}{0.33\textwidth}
		\includegraphics[width=\linewidth]{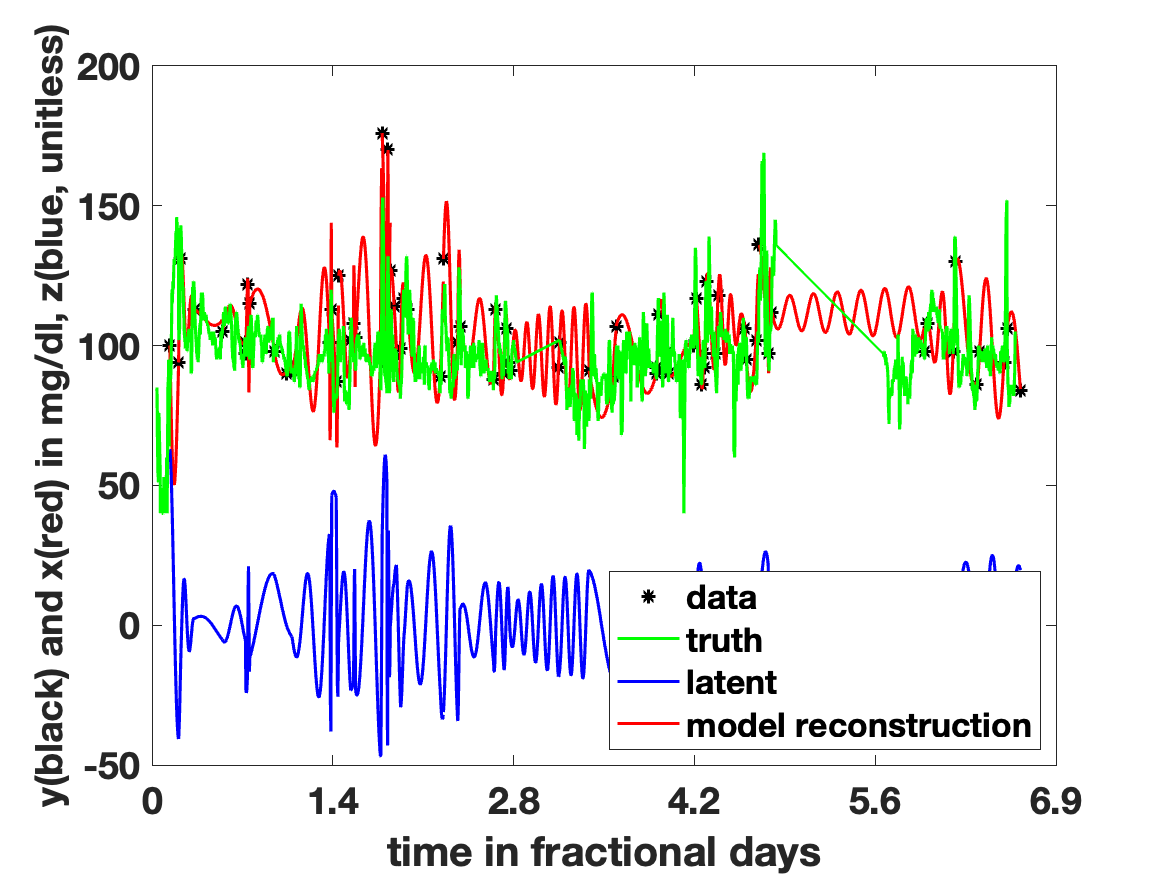}
	\end{subfigure}%
	\begin{subfigure}{0.33\textwidth}
		\includegraphics[width=\linewidth]{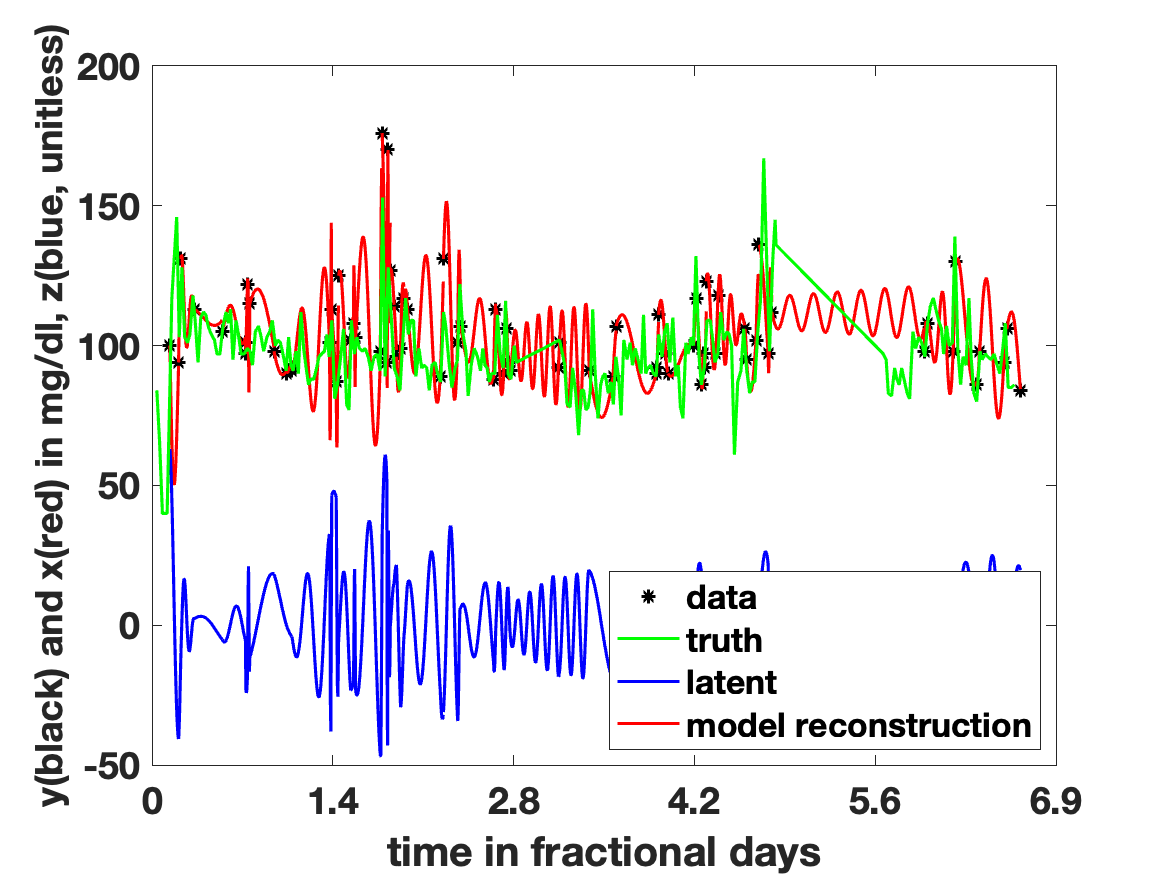}
	\end{subfigure}%
	\begin{subfigure}{0.33\textwidth}
		\includegraphics[width=\linewidth]{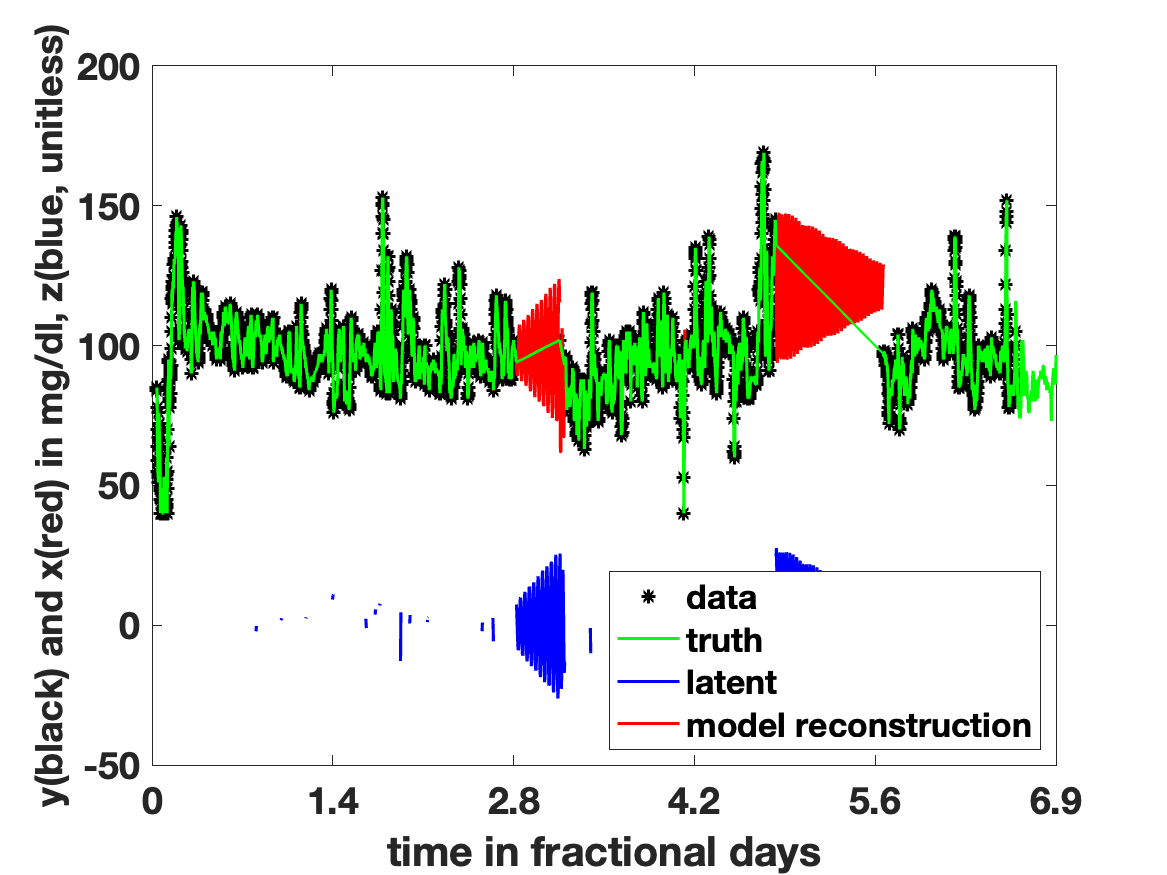}
	\end{subfigure}
	
	\begin{subfigure}{0.33\textwidth}
		\includegraphics[width=\linewidth]{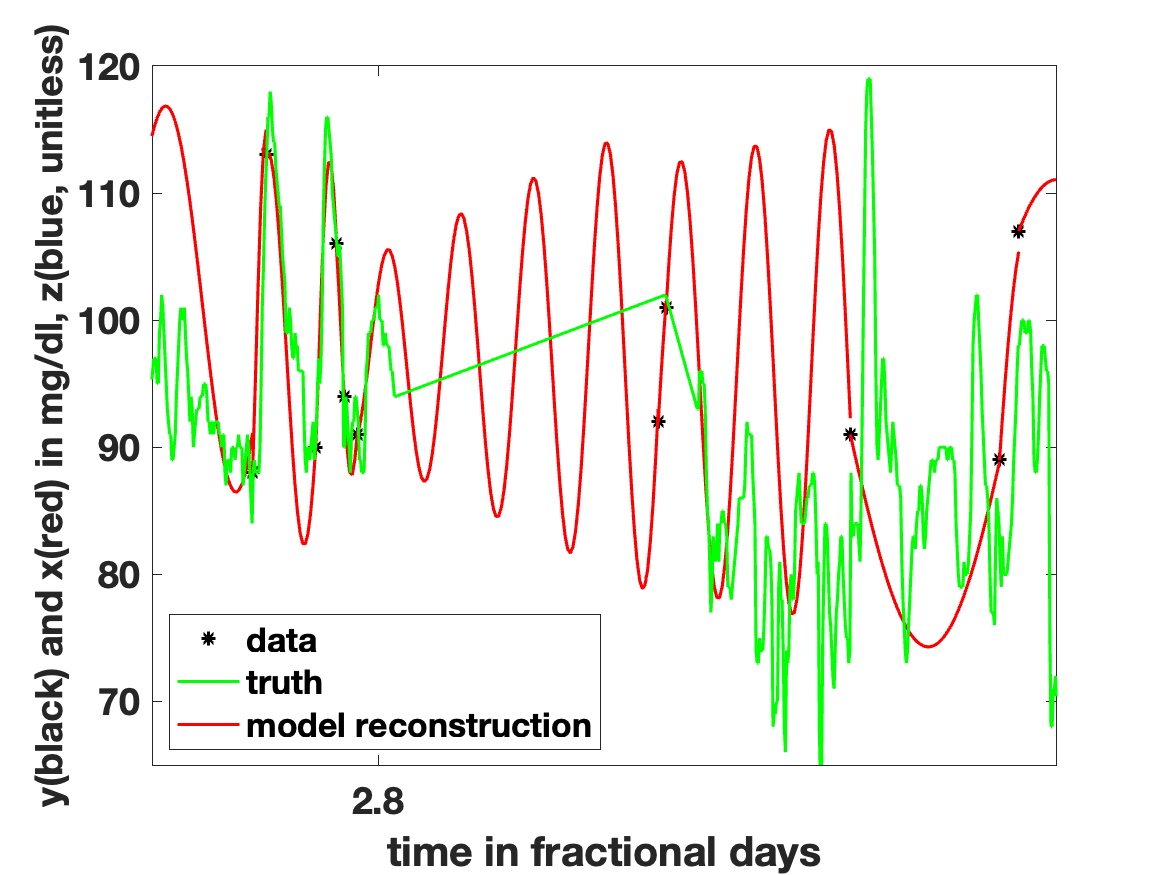}
	\end{subfigure}%
	\begin{subfigure}{0.33\textwidth}
		\includegraphics[width=\linewidth]{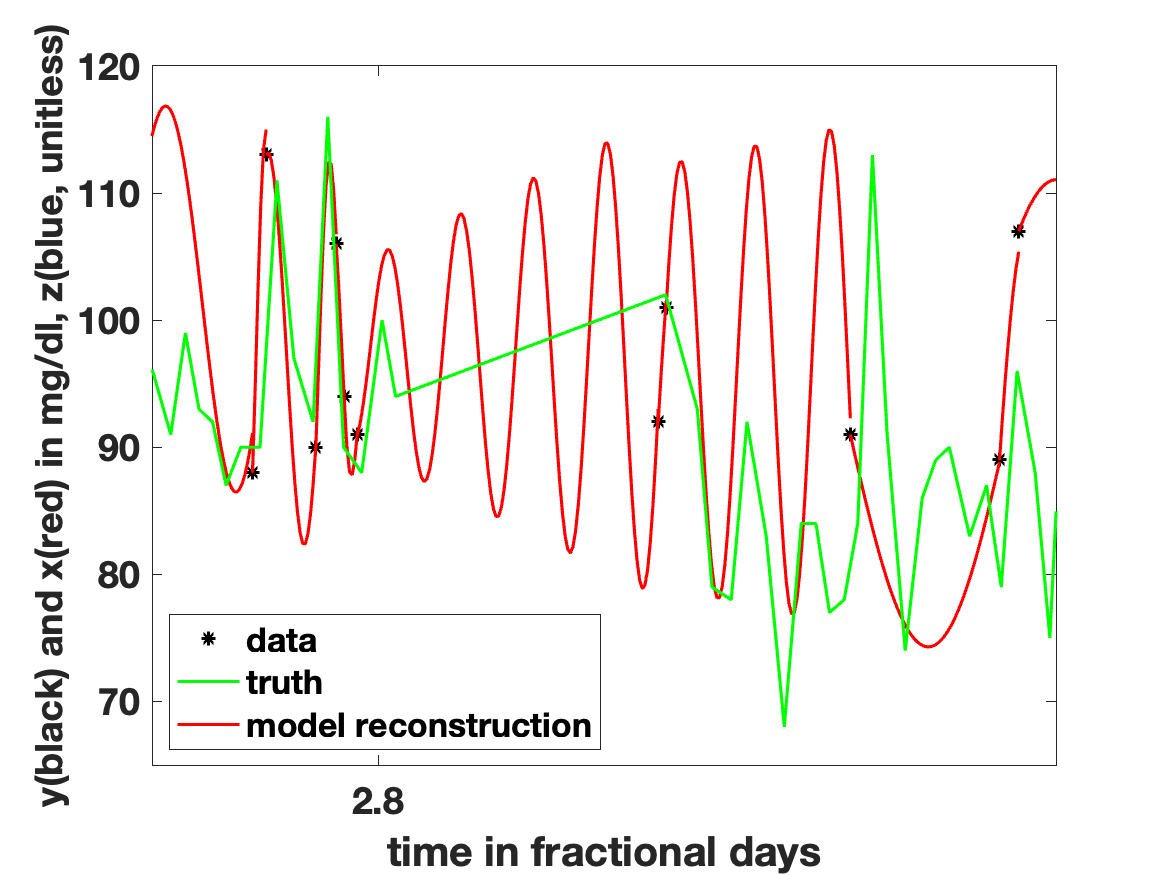}
	\end{subfigure}%
	\begin{subfigure}{0.33\textwidth}
		\includegraphics[width=\linewidth]{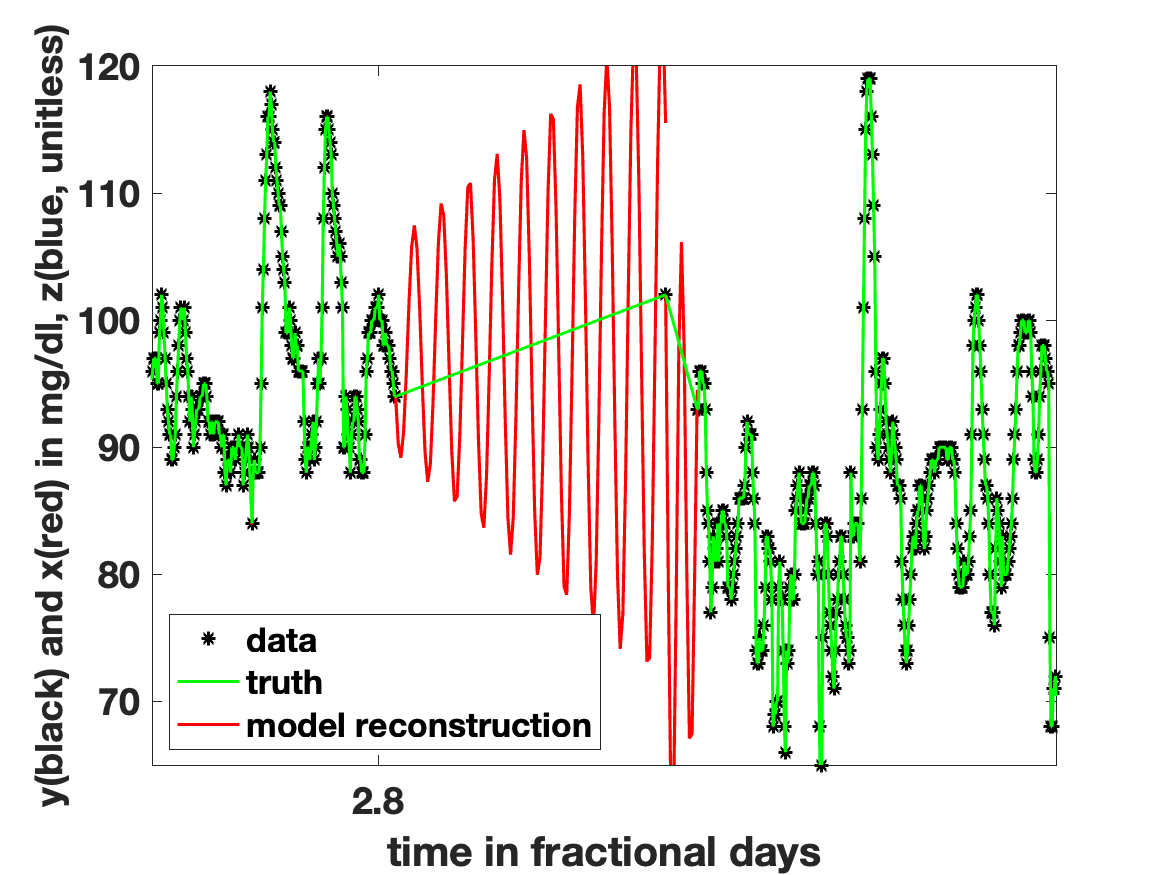}
	\end{subfigure}

	\caption{Given the damped, driven wild patient data, the model was estimated assuming an underlying oscillatory solution (set according to section \ref{sec:comp_proceedure}) for the three measurement functions, sparse clinician measured ($h_1$, top/bottom left), random ($h_2$, top/bottom center) and dense ($h_3$, top/bottom right). We can see the point-wise estimates are good in the lower plots, the model captures the large peaks and troughs but only captures the correct frequency for densely measured data. When data are missing the model relaxes to oscillatory dynamics. Note that the straight lines in green are not the truth but a linear interpolation in the absence of data.}
	\label{fig:wild_traj_osc}
\end{figure}

\begin{figure}
	\centering
	\begin{subfigure}{0.33\textwidth}
		\includegraphics[width=\linewidth]{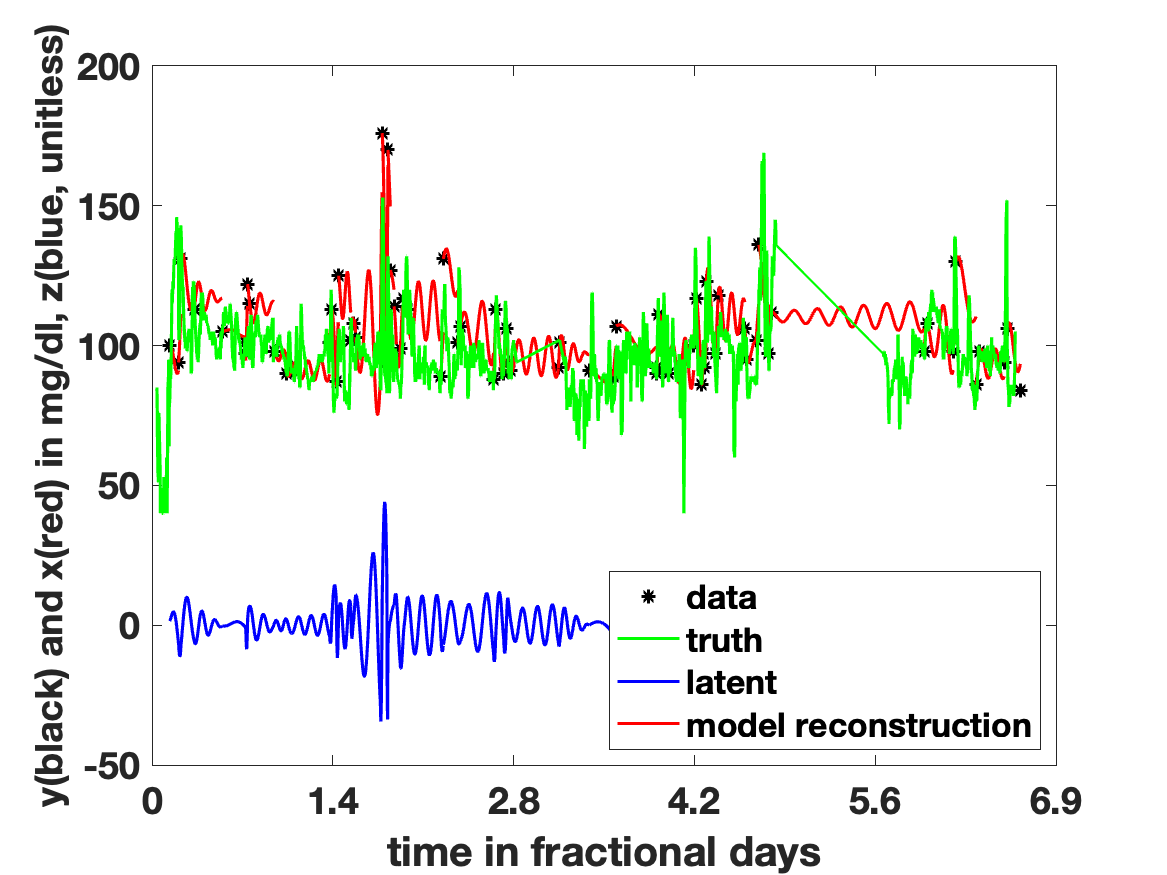}
	\end{subfigure}%
	\begin{subfigure}{0.33\textwidth}
		\includegraphics[width=\linewidth]{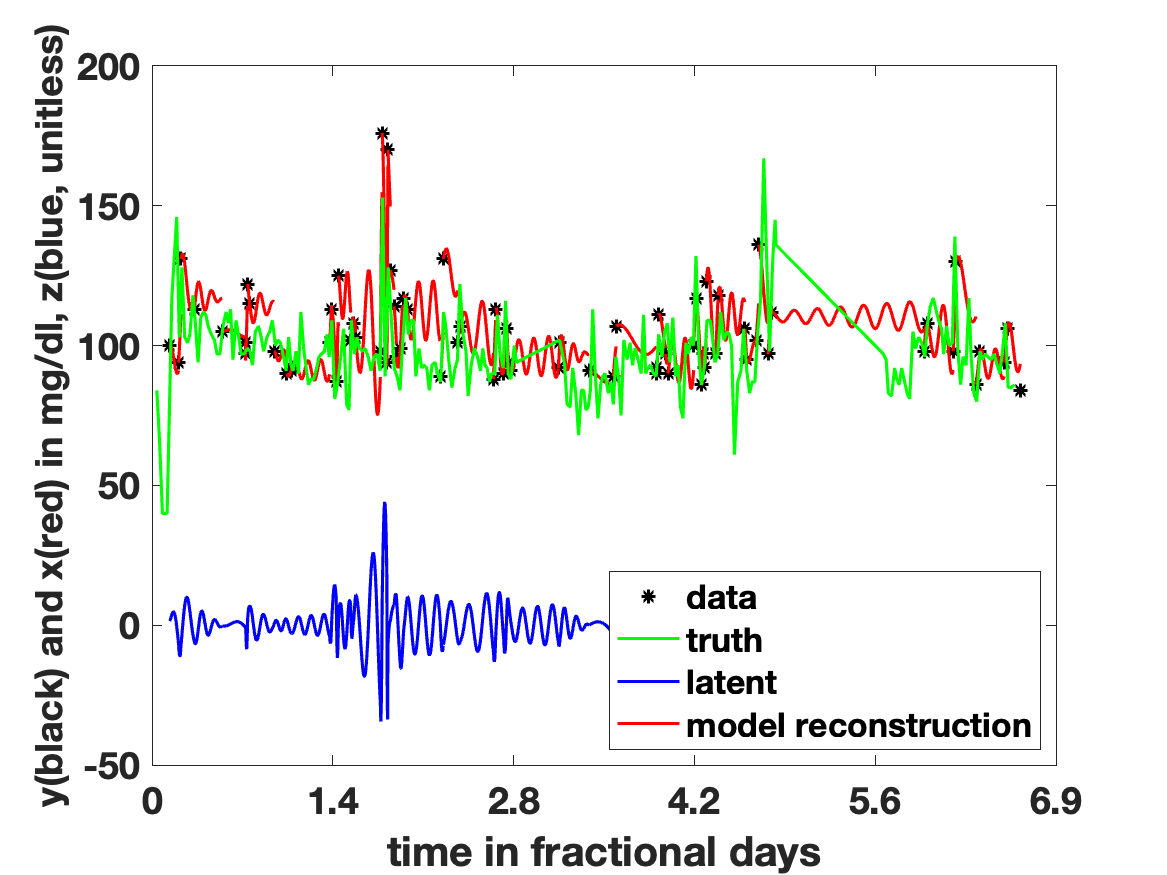}
	\end{subfigure}%
	\begin{subfigure}{0.33\textwidth}
		\includegraphics[width=\linewidth]{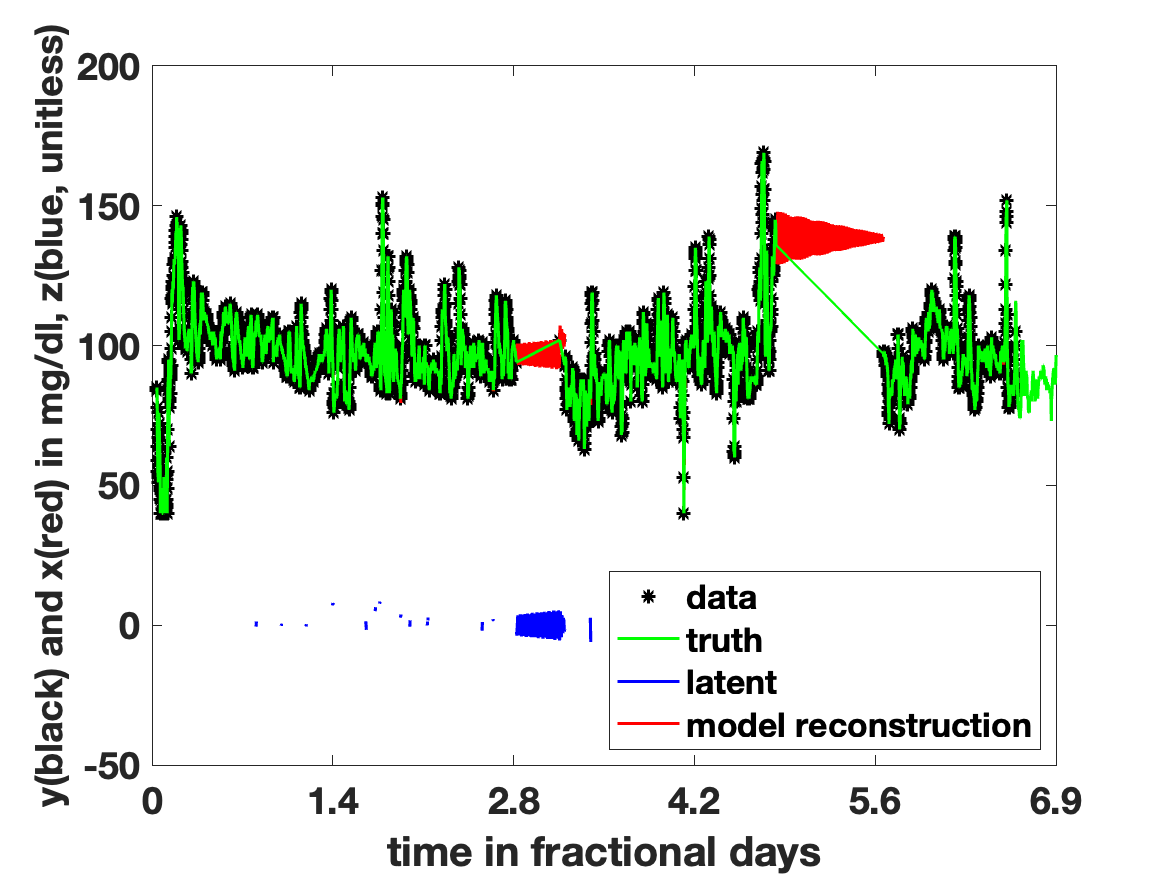}
	\end{subfigure}
	\begin{subfigure}{0.33\textwidth}
		\includegraphics[width=\linewidth]{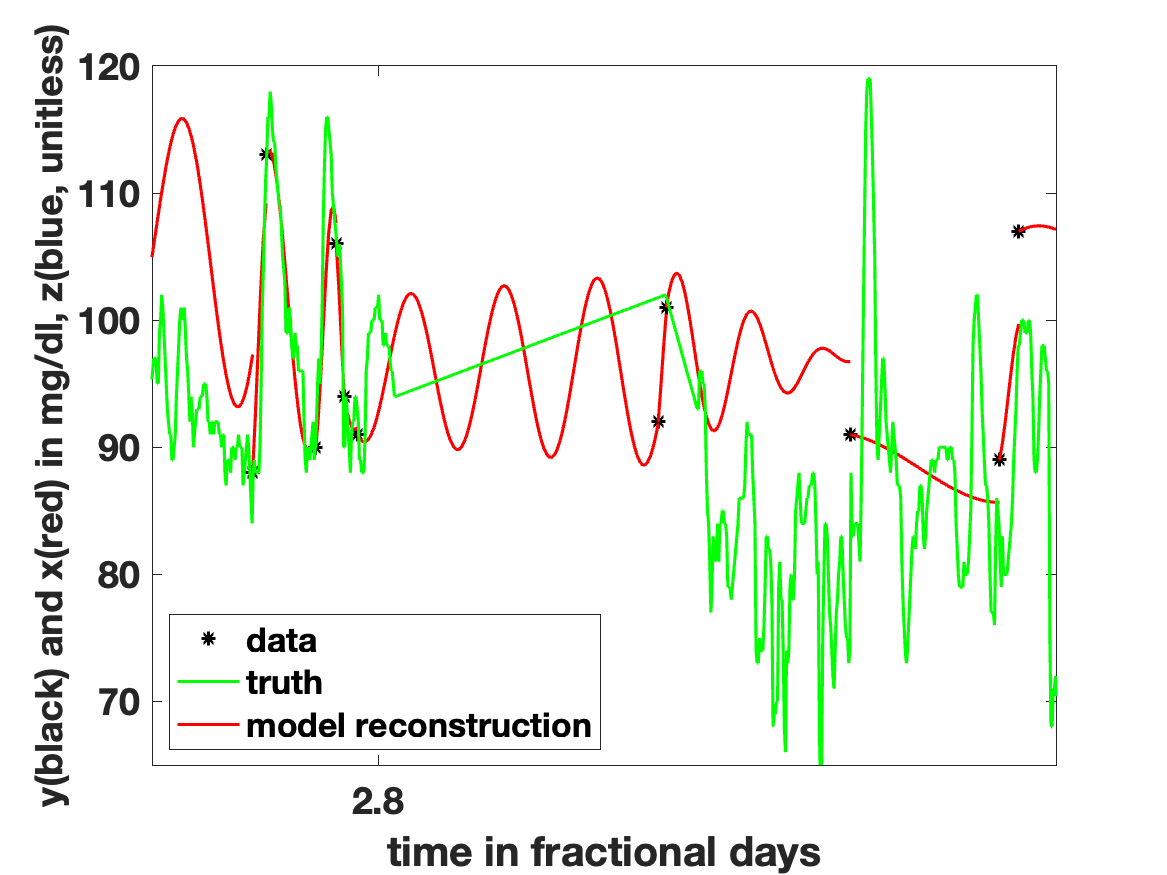}
	\end{subfigure}%
	\begin{subfigure}{0.33\textwidth}
		\includegraphics[width=\linewidth]{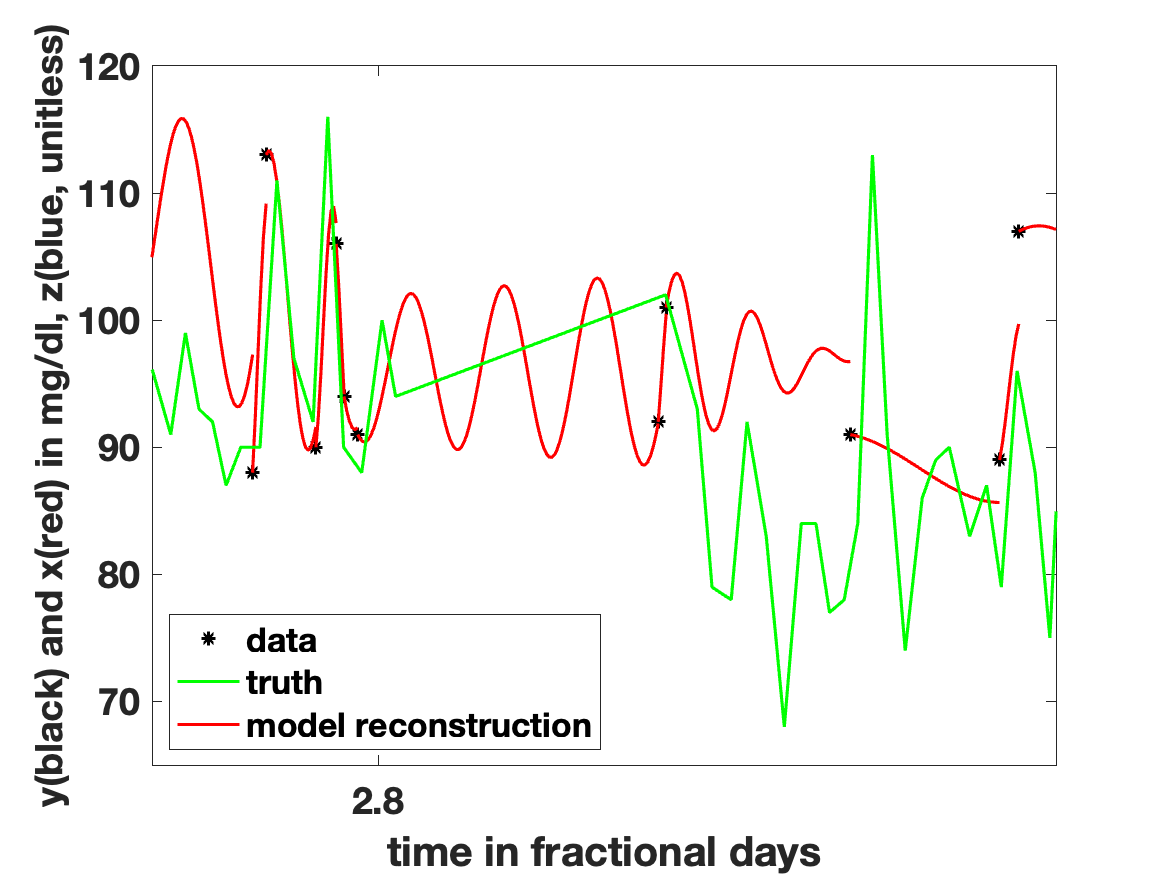}
	\end{subfigure}%
	\begin{subfigure}{0.33\textwidth}
		\includegraphics[width=\linewidth]{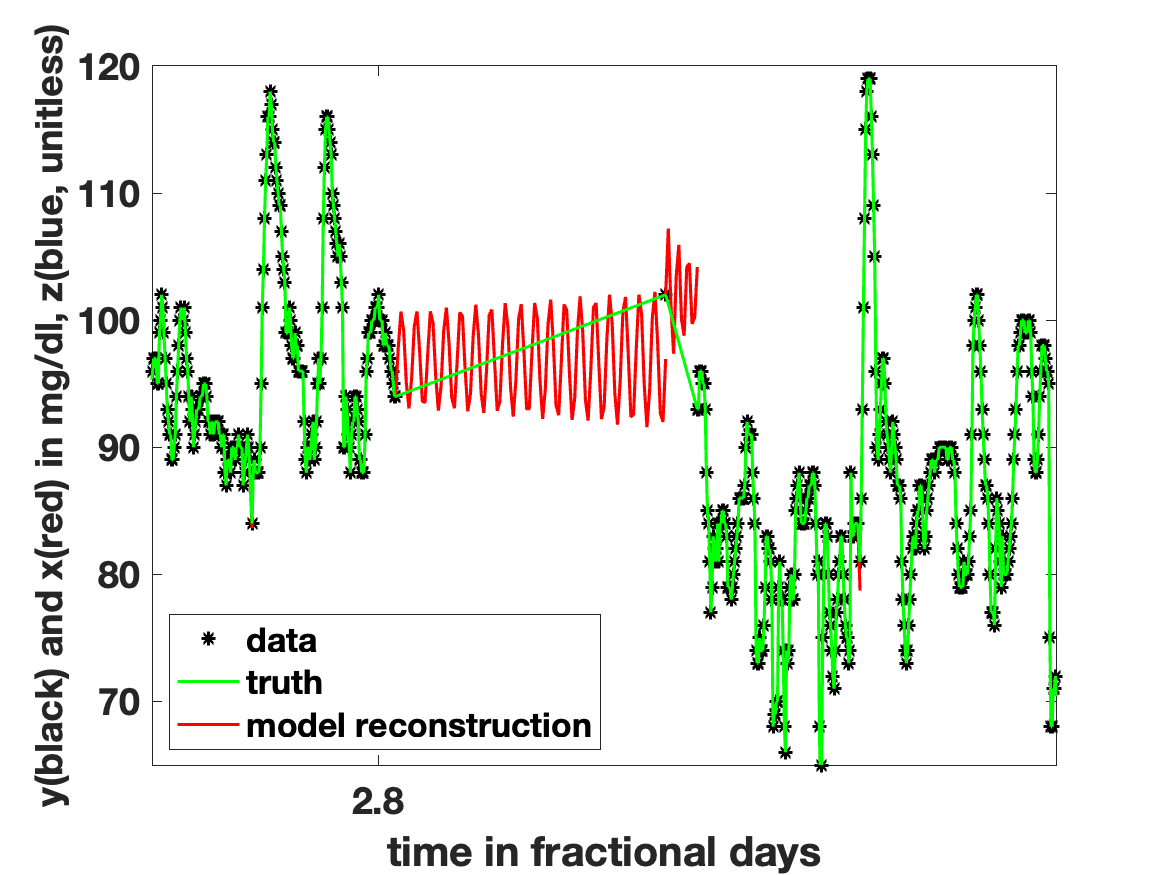}
	\end{subfigure}

	\caption{Given the damped, driven wild patient data, the model was estimated assuming no oscillatory solution and relaxed to its non-oscillatory solution controlled by the $\alpha$ parameter and by setting the baseline amplitude $a$ to zero, for the three measurement functions, sparse clinician measured ($h_1$, top/bottom left), random ($h_2$, to/bottom, center) and dense ($h_3$, top/bottom right). We can see the point-wise estimates are good in the lower plots, the model captures the large peaks and troughs but only captures the correct frequency for densely measured data. When data are missing the model relaxes to non-oscillatory dynamics.  Again note that the straight lines in green are not the truth but a linear interpolation in the absence of data.}
	\label{fig:wild_traj_no_osc}
\end{figure}

\begin{figure}
	\centering
	\begin{subfigure}{0.33\textwidth}
		\includegraphics[width=\linewidth]{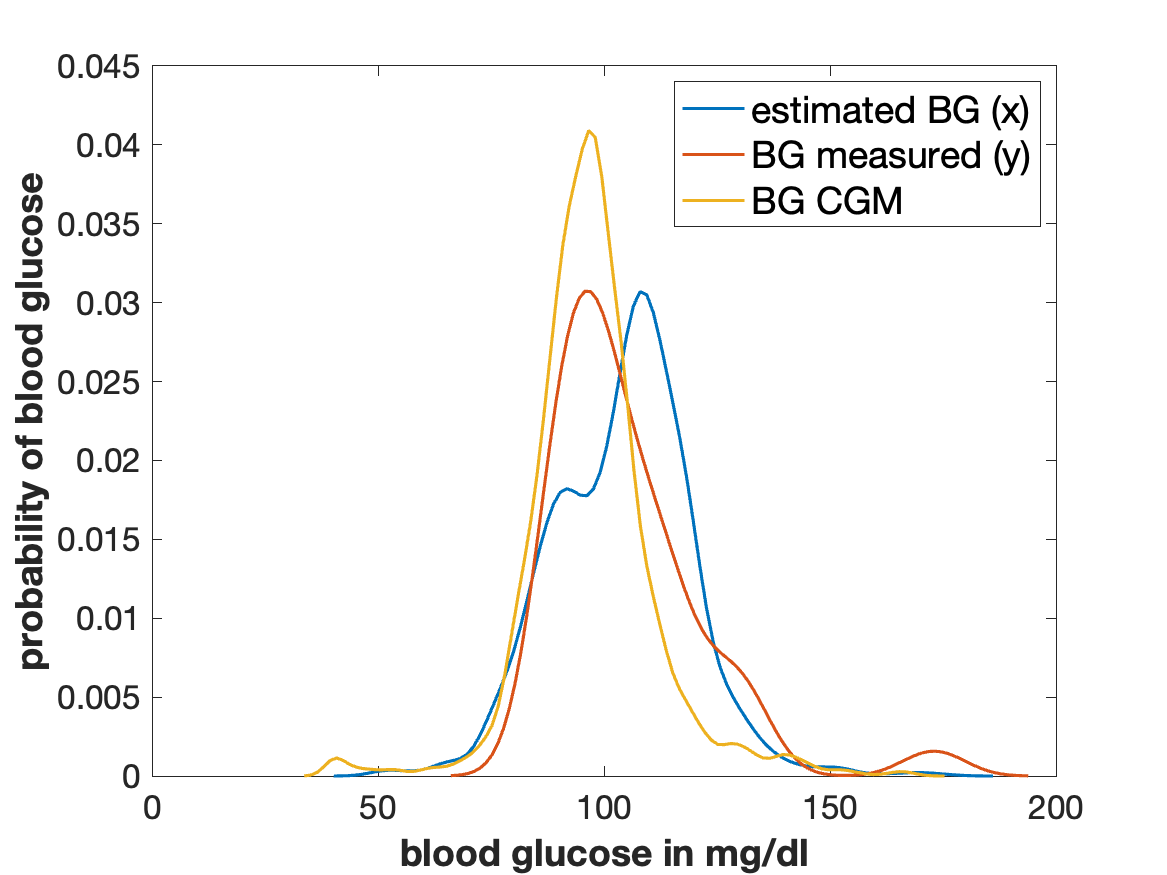}
	\end{subfigure}%
	\begin{subfigure}{0.33\textwidth}
		\includegraphics[width=\linewidth]{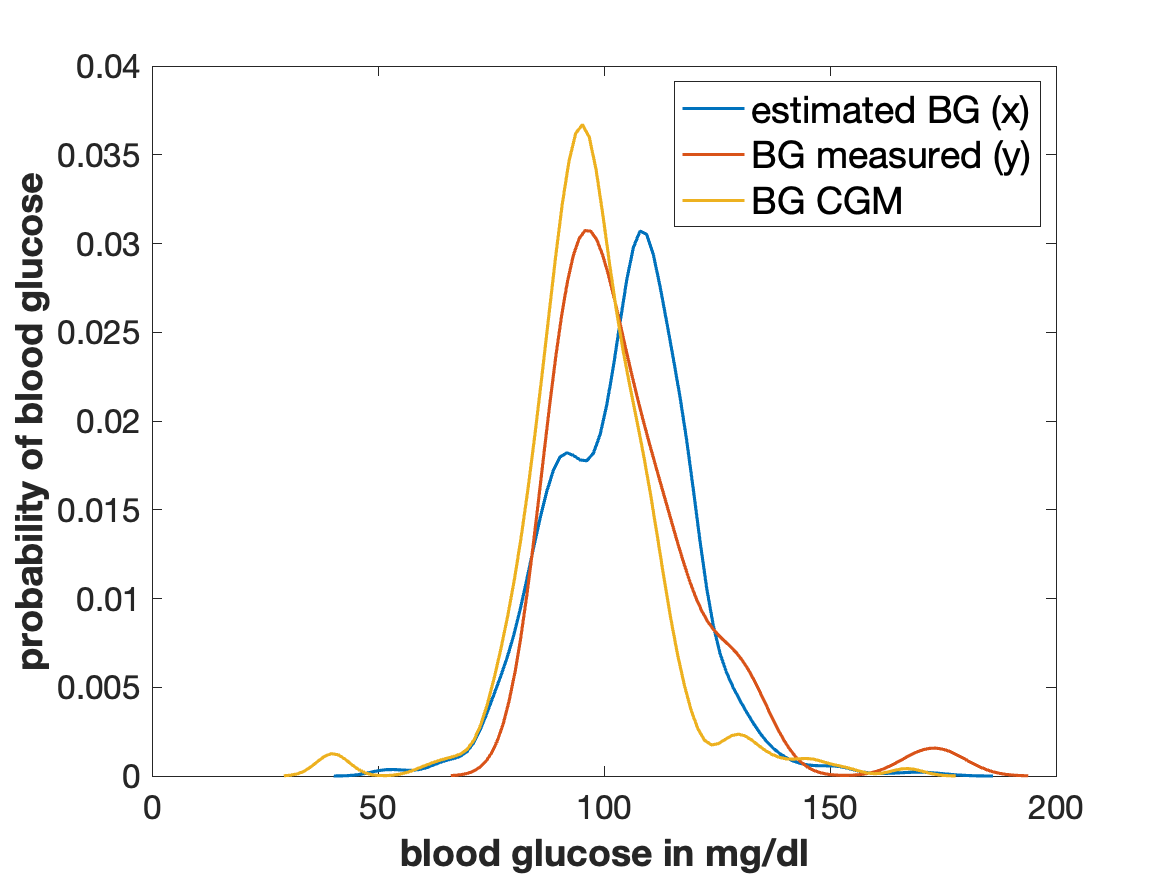}
	\end{subfigure}%
	\begin{subfigure}{0.33\textwidth}
		\includegraphics[width=\linewidth]{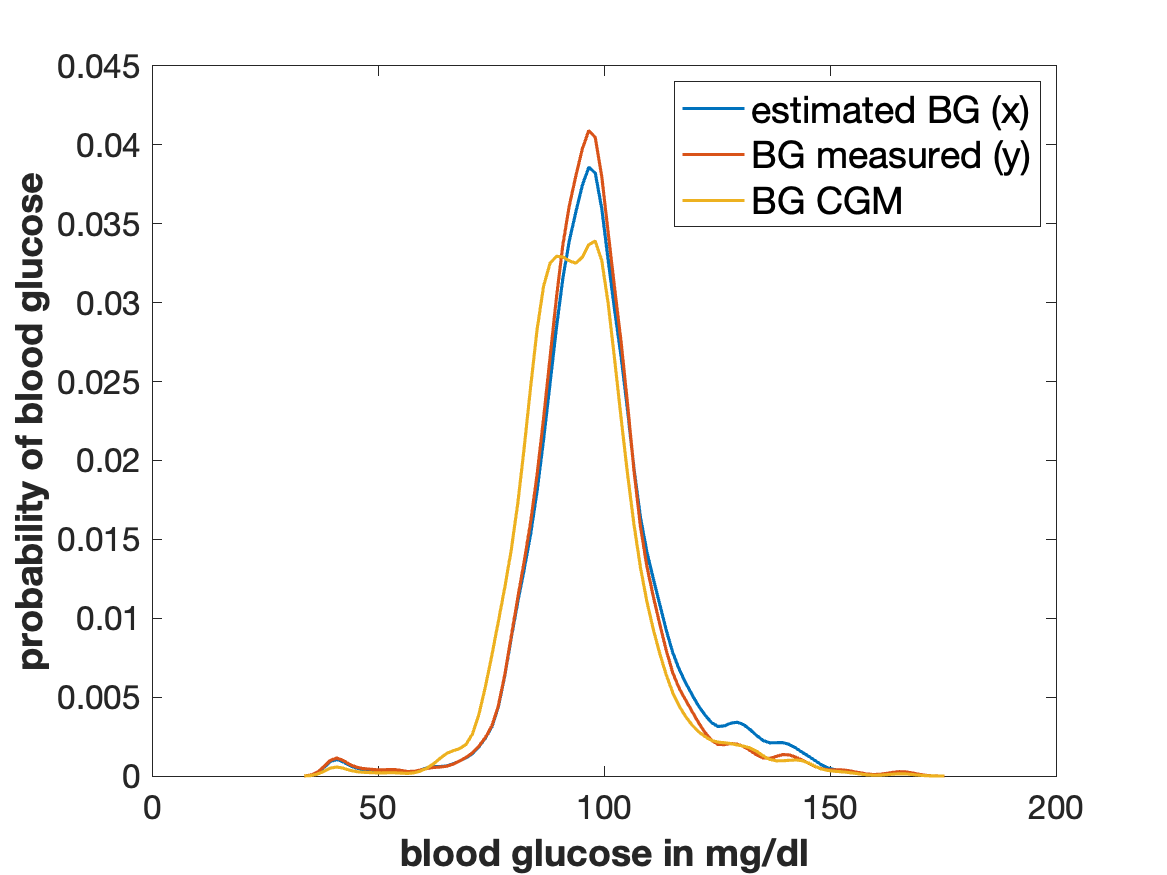}
	\end{subfigure}

	\caption{Given the damped, driven wild patient data, the model was estimated assuming an oscillatory solution and relaxed to its oscillatory solution (set according to section \ref{sec:comp_proceedure}) for the three measurement functions, sparse clinician measured ($h_1$, left), random ($h_2$, center) and dense ($h_3$, right). Note \emph{BG measured} denotes data available to the model when it is estimated and \emph{estimated BG} denotes the model-estimated invariant measure of the data. The densely measured case ($h_3$) is likely the closest representative of a gold standard baseline, again for data measured frequently in time.  We can see the model captures the distributions well. While the clinician-driven measurement times, $h_1$, may outperform random measurement times, $h_2$, neither are perfect as they are balanced against the point-wise and model-coherence loss functions. Additionally, the case here where the solution is assumed to be oscillatory did not perform as well as the case where the solution was assumed to be non-oscillatory for the sparse measurement functions $h_1$ and $h_2$.   When the model is estimated with data from the dense measurement function, the invariant measure of these data is estimated very accurately and was not appreciably different from the case where the model assumed non-oscillatory baseline dynamics shown in Fig. \ref{fig:wild_pdf_no_osc}.}
	\label{fig:wild_pdf_osc}
\end{figure}

\begin{figure}
	\centering
	\begin{subfigure}{0.33\textwidth}
		\includegraphics[width=\linewidth]{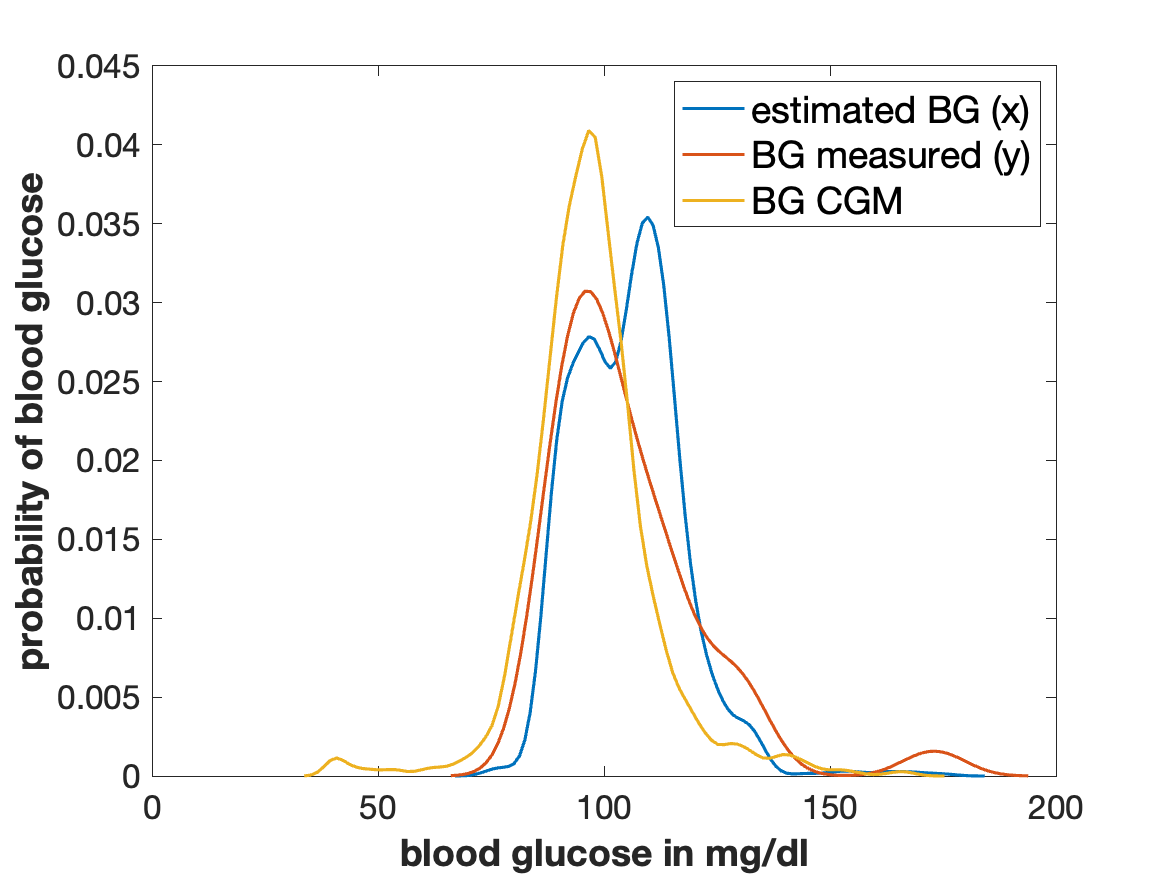}
	\end{subfigure}%
	\begin{subfigure}{0.33\textwidth}
		\includegraphics[width=\linewidth]{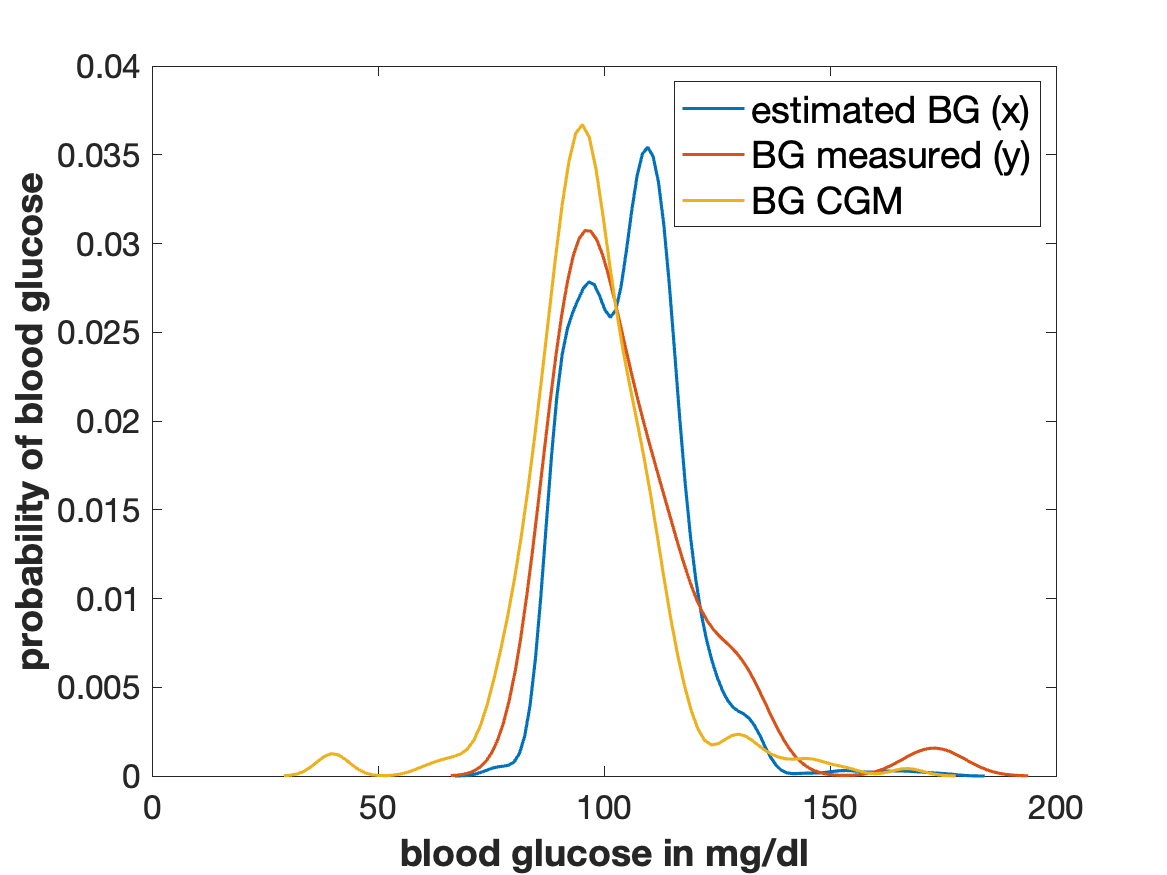}
	\end{subfigure}%
	\begin{subfigure}{0.33\textwidth}
		\includegraphics[width=\linewidth]{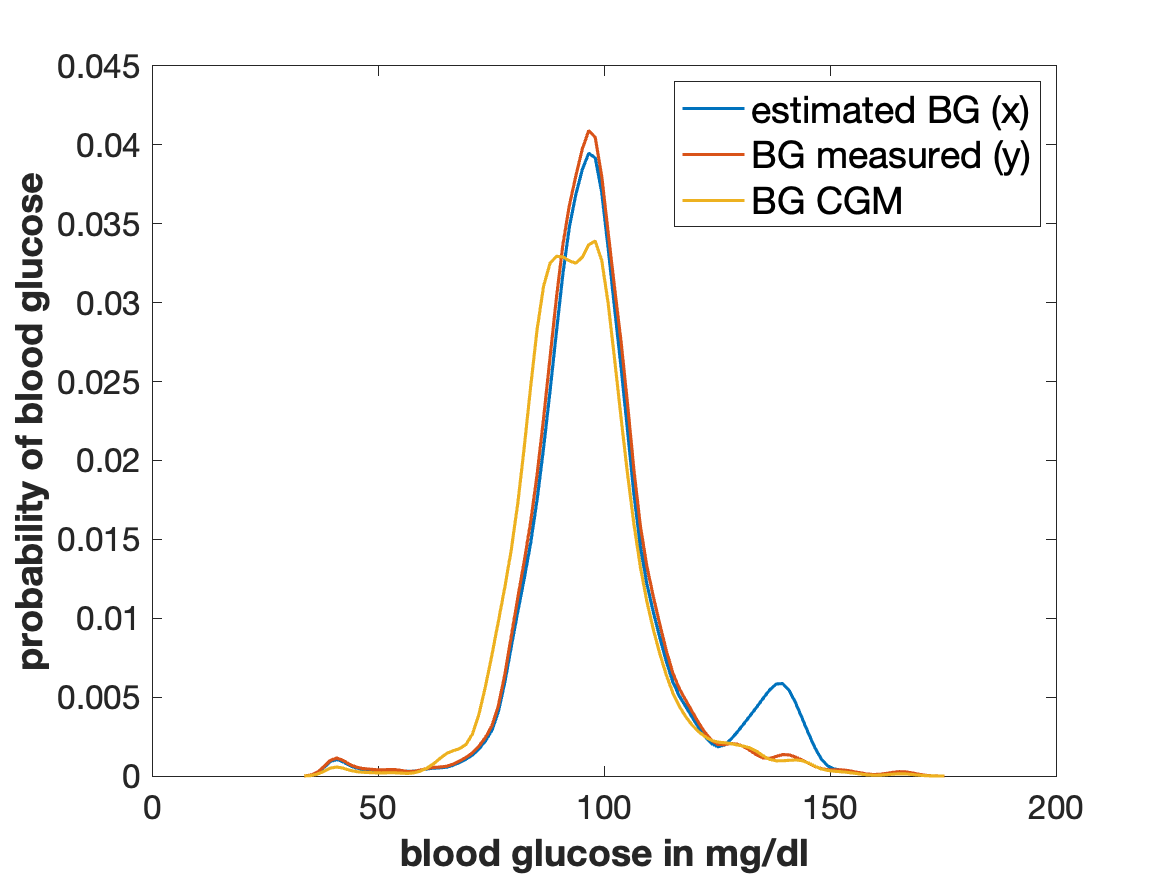}
	\end{subfigure}

	\caption{Given the damped, driven wild patient data, the model was estimated assuming no oscillatory solution and relaxed to its non-oscillatory solution (baseline $a=0$) for the three measurement functions, sparse clinician measured ($h_1$, left), random ($h_2$, center) and dense ($h_3$, right). Note \emph{BG measured} denotes data available to the model when it is estimated and \emph{estimated BG} denotes the model-estimated invariant measure of the data. The densely measured case ($h_3$) is likely the closest representative of a gold standard baseline, again for data measured frequently in time.  We can see the model captures the distributions well. While the clinician-driven measurement times, $h_1$, may outperform random measurement times, $h_2$, neither are perfect as they are balanced against the point-wise and model-coherence loss functions. When the model is estimated with data from the dense measurement function, the invariant measure of these data is estimated very accurately.}
	\label{fig:wild_pdf_no_osc}
\end{figure}

Moving to the in-the-wild case where the generating dynamics are damped, driven---by punctuated nutrition consumption---oscillations, we estimated the model assuming both non-oscillatory and oscillatory (resp. Figs. \ref{fig:wild_traj_no_osc} and \ref{fig:wild_traj_osc}) dynamics, by controlling the model's amplitude and relaxations parameter $\alpha$. In the cases where measurements were sparse, $h_1$ and $h_2$, assuming the underlying dynamics were non-oscillatory lead to better performance compared with assuming oscillatory dynamics both point-wise, Figs. \ref{fig:wild_traj_osc}-\ref{fig:wild_traj_no_osc}, and distributionally, Figs.   \ref{fig:wild_pdf_osc}-\ref{fig:wild_pdf_no_osc}. When the data were measured densely with $h_3$ there was not an appreciable difference between the assumed oscillatory and non-oscillatory cases with respect to point-wise and distributional fits. The dependence of the results on the measurement functions highlights the important role measurement functions can play on model estimation. Together these results imply that we can include global properties and that their inclusion does improve estimation, supporting an affirmative answer to \emph{Q1a}, \emph{Q1b}, and \emph{Q1c}. Additionally, related to  \emph{Q2b} and \emph{Q2c}, the measurement functions clearly impact model estimation performance.  Specifically, accuracy of estimation of the off-data trajectory is highly dependent on the measurement function, but including information about the correct underlying dynamics and global properties such as the invariant measure of the data does qualitatively improve estimation (\emph{Q2c}). Additionally, regardless of measurement function, we are able to accurately estimate the invariant measure of the data (\emph{Q2b}), yet setting the underlying dynamics correctly does improve estimation of the invariant measure of these data when these data are particularly sparsely measured.

\subsection{Impact of parameter flex over the estimate window to account for measurement time errors and nonstationarity (\emph{Q3}) and impacts of different measurement functions on model estimation (\emph{Q2}) }


We designed our new algorithm to allow parameters to vary within the estimation window according to the $L_{4+l}$ loss function that models the transition probabilities of the $l^{th}$ parameter value between measurements with a decay rate $d^n_l$ and uncertainty $\sigma_l$. Modeling and estimating the transition probabilities between measurements was designed to solve two problems, both related to the goal of  giving the model enough flexibility to estimate point-wise and global qualitative dynamics with sparse data. \emph{First,} data measurement recording times have errors while mechanistic models with oscillatory dynamics have relatively rigid orbits with rigid frequencies that can be tuned, to some degree, according to model parameters. Allowing parameters to vary allows the models to have more flexible trajectories that can accommodate measurement error times without resorting to zeroing out the oscillations. \emph{Second,} when data are sparse, to maximize the number of data points included we are forced to estimate the model over longer time periods, increasing the potential impact to nonstationarity. Allowing the parameters to flex over the estimation window can potentially lead to more productive parameter estimation in the presence of this nonstationarity. Here we are considering specifically \emph{Q3}, the impact of nonstationarity and measurement functions on our flexible parameter estimation methods. Additionally and again, the ICU patient has nonstationary physiological mechanics and oscillatory  dynamics except when tube-feeding is interrupted while the in-the-wild patient is likely quite stationary with dynamics that are roughly a noisy damped, driven oscillator.

\paragraph{\textbf{ICU}}
All measurement functions impact parameter estimation. The random $h_2$ and dense, $h_3$ measurement functions estimate the \emph{amplitudes}, $a$ and $b$, similarly and both differ from the clinician-driven measurement function, $h_1$.  All three measurement functions differ in their estimation of  frequency of oscillations, $\omega$, which is both interesting because the sparse measurement functions differed and expected because only the dense measurement function could truly resolve the frequency. Nevertheless, all measurement functions placed the amplitude parameters in roughly the same range, which is particularly important for two reasons. First, it was the amplitude estimation that was failing in our motivating example, and here all the parameter estimates included oscillatory solutions. Second, the most important feature of the blood glucose dynamics for clinical decision-making is the glycemic range defined by the amplitude of oscillations \cite{george_eval}, and regarding this, all measurement functions lead to qualitatively similar amplitudes. Finally, frequency estimates differ between all measurement functions, as expected given that the sparse measurements contain very little frequency information. Overall, the parameter flex did impact model estimation, and the parameter flex was dependent on the measurement functions.

\begin{figure}
	\centering
	\begin{subfigure}{0.33\textwidth}
		\includegraphics[width=\linewidth]{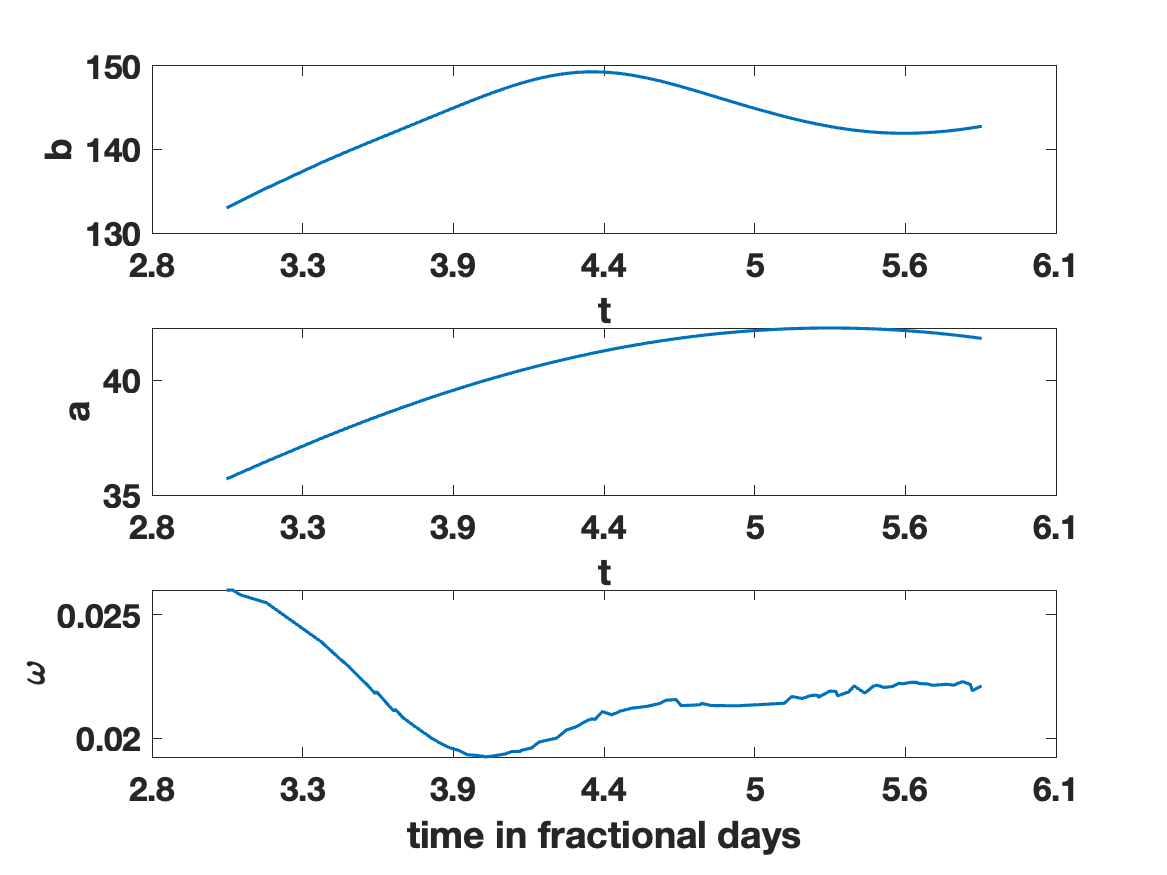}
	\end{subfigure}%
	\begin{subfigure}{0.33\textwidth}
		\includegraphics[width=\linewidth]{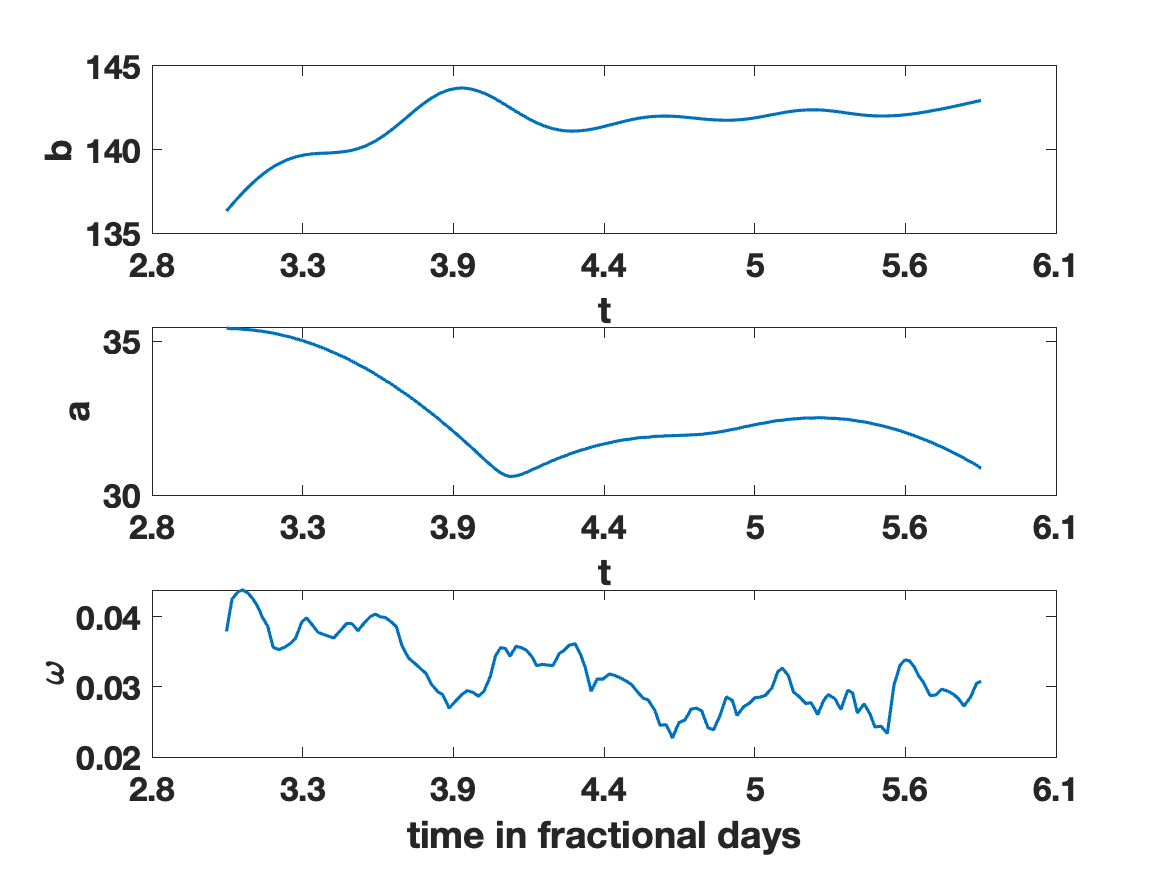}
	\end{subfigure}%
	\begin{subfigure}{0.33\textwidth}
		\includegraphics[width=\linewidth]{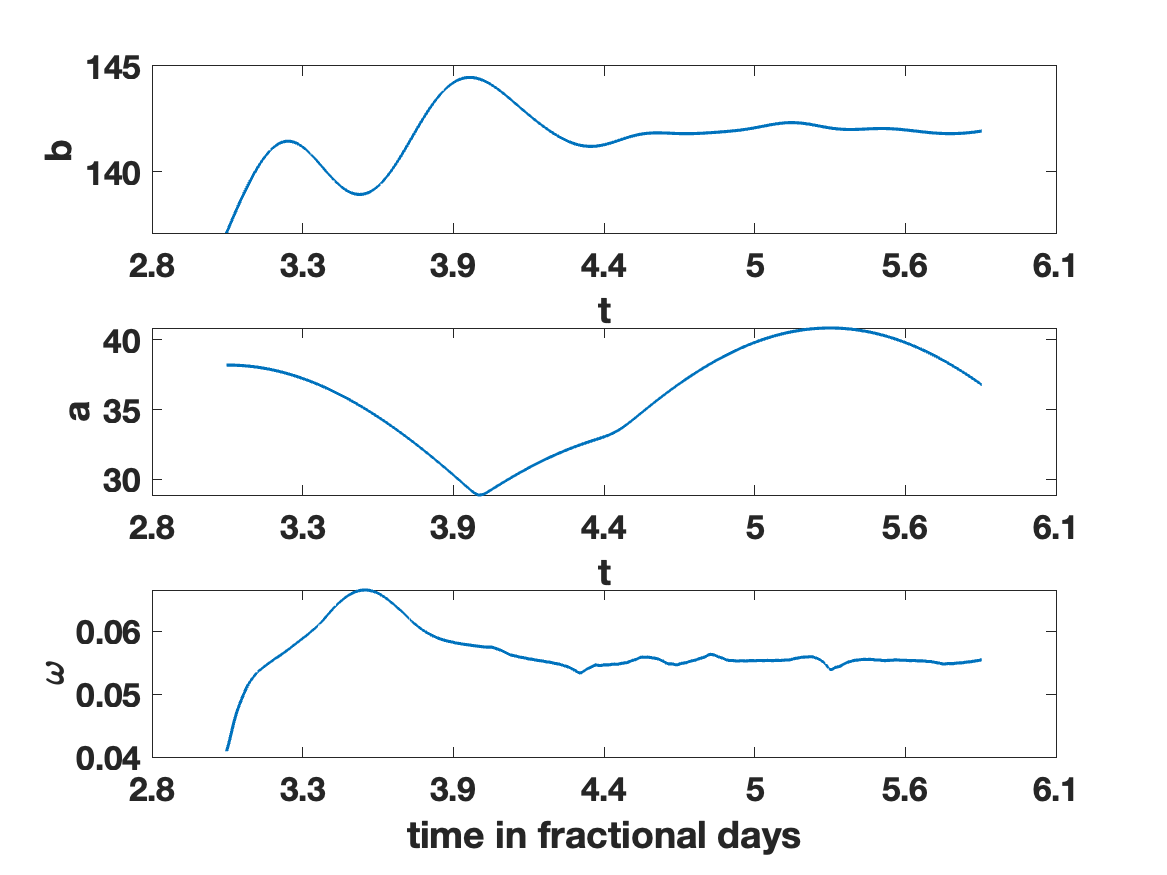}
	\end{subfigure}
	\caption{Given the oscillatory ICU glycemic dynamics and the model basic state assuming oscillatory dynamics, we can see that all measurement functions, $h_1$ (top), $h_2$ (center) and $h_3$ bottom, utilized the temporal flexibility in the parameters to decrease estimation error, but in very different ways. Interestingly, the sparse random and dense measurement functions ($h_2$, $h_3$) had nearly the same flex in \emph{amplitude} estimation, $a$ and $b$, sometimes in direct opposition to flexible parameter trends of the clinician (human) measurement case ($h_1$). All three measurement functions differ in their estimation of  frequency of oscillations, $\omega$. This result has particularly important potential implications for solving inverse problems to estimate parameters because it implies that differently measured data, even when the number of data points are similar, can have a substantial impact on the estimability of model parameters.}
\end{figure}

\begin{figure}
	\centering
	\begin{subfigure}{0.33\textwidth}
		\includegraphics[width=\linewidth]{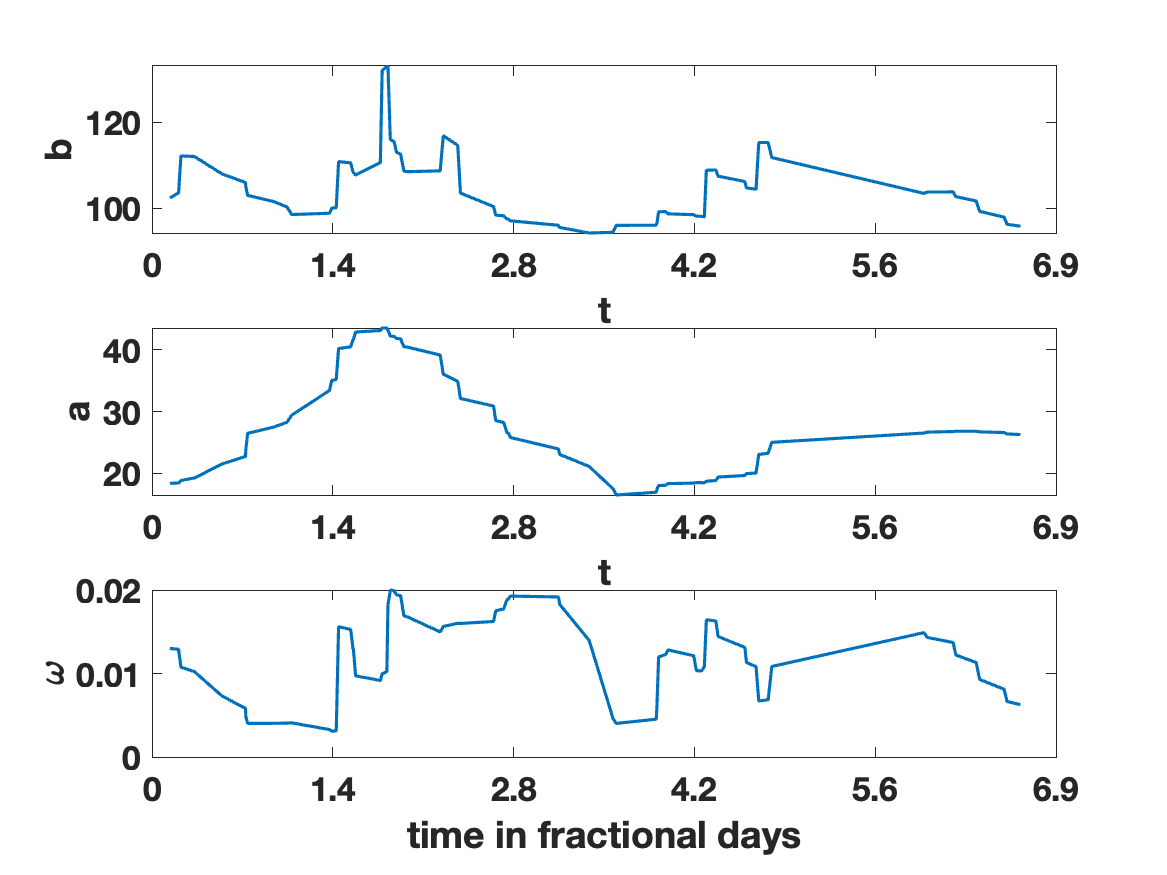}
	\end{subfigure}%
	\begin{subfigure}{0.33\textwidth}
		\includegraphics[width=\linewidth]{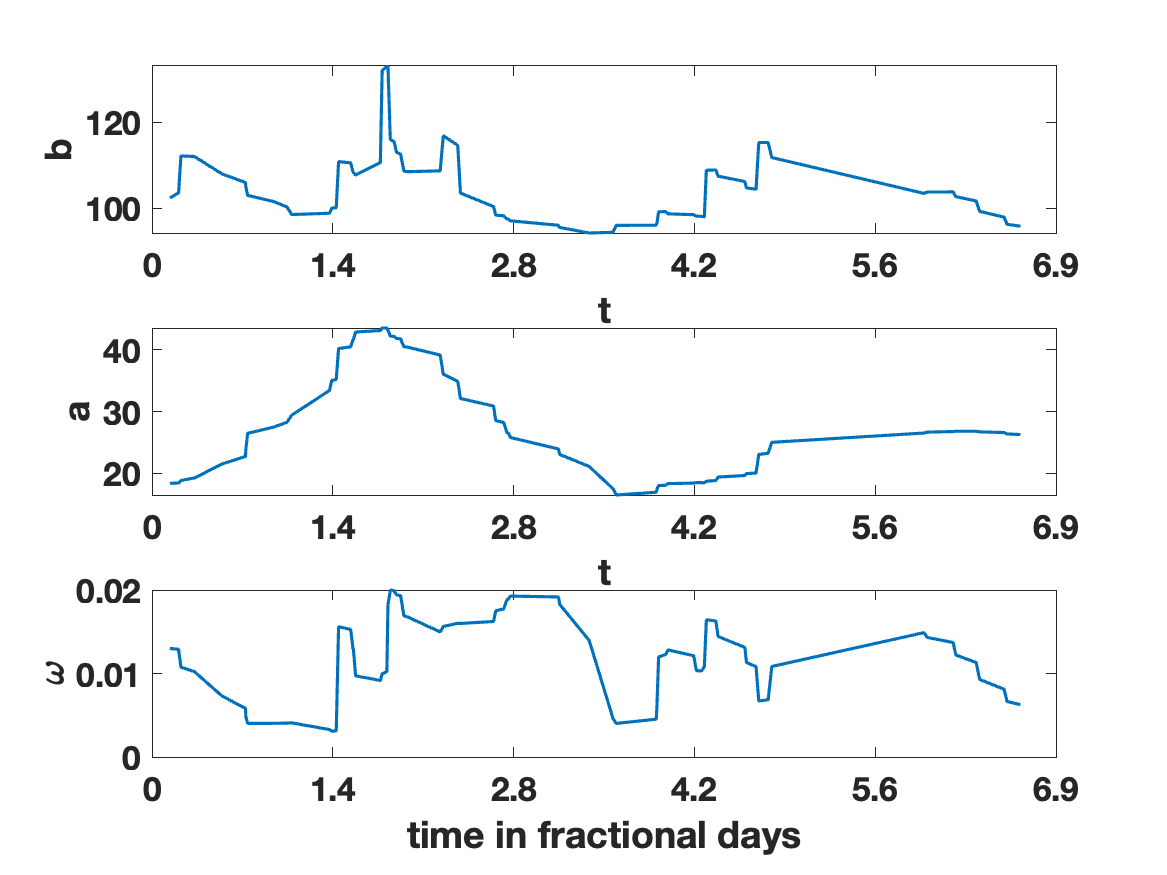}
	\end{subfigure}%
	\begin{subfigure}{0.33\textwidth}
		\includegraphics[width=\linewidth]{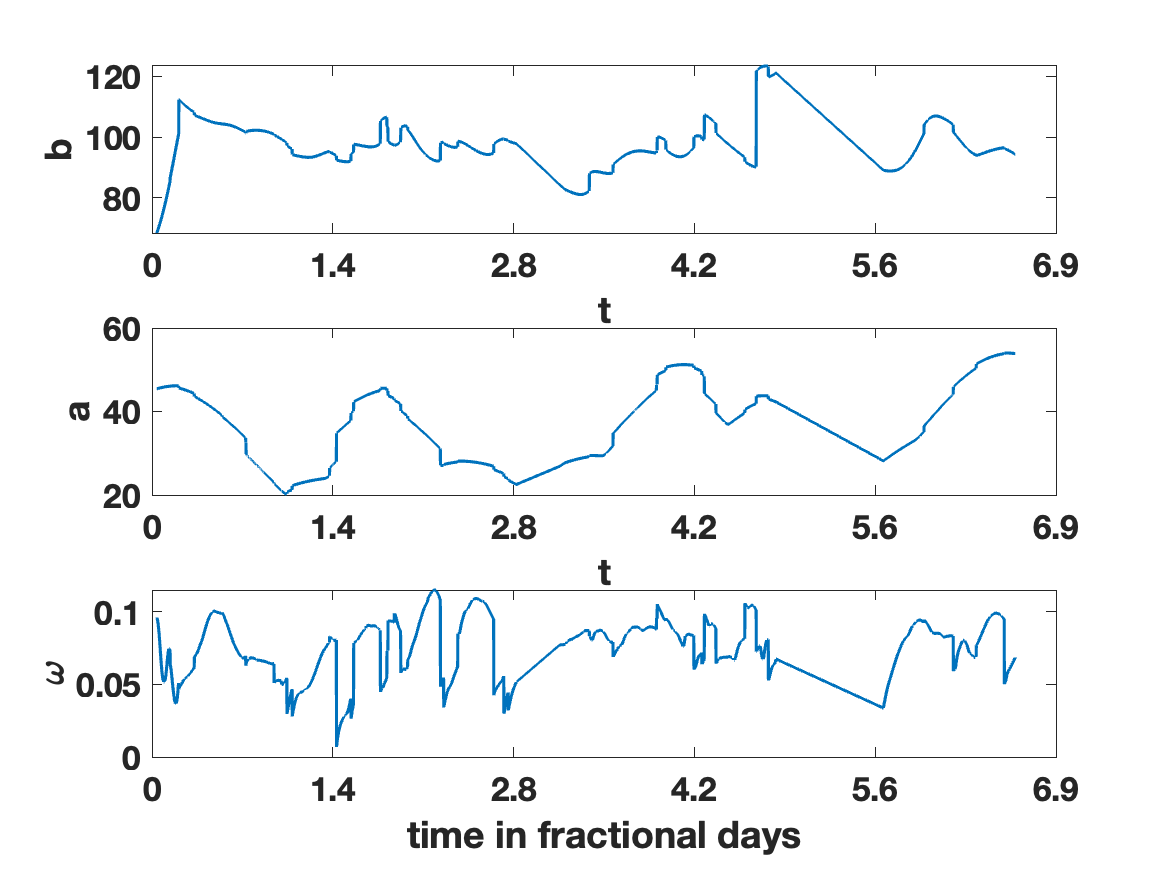}
	\end{subfigure}

	\caption{Given the in-the-wild glycemic damped, driven oscillatory dynamics we can see and the model basic state assuming \emph{oscillatory} dynamics, we can see that all measurement functions, $h_1$ (left), $h_2$ (center) and $h_3$ (right), utilized the temporal flexibility in the parameters to decrease estimation error, but in very different ways. Interestingly, the sparse patient-defined and random measurement functions ($h_1$, $h_2$) had nearly the same flex over the course of the time window, and both were quite different from the dense measurement case ($h_3$). This result has particularly important potential implications for solving inverse problems to estimate parameters because differences between parameters estimated with sparse versus dense data can have a substantial impact on the estimability of model parameters. Understanding the degree of importance and uncertainty in the inverse problems setting will be important.}
	\label{fig:wild_osc_parms}
\end{figure}

\begin{figure}
	\centering
	\begin{subfigure}{0.33\textwidth}
		\includegraphics[width=\linewidth]{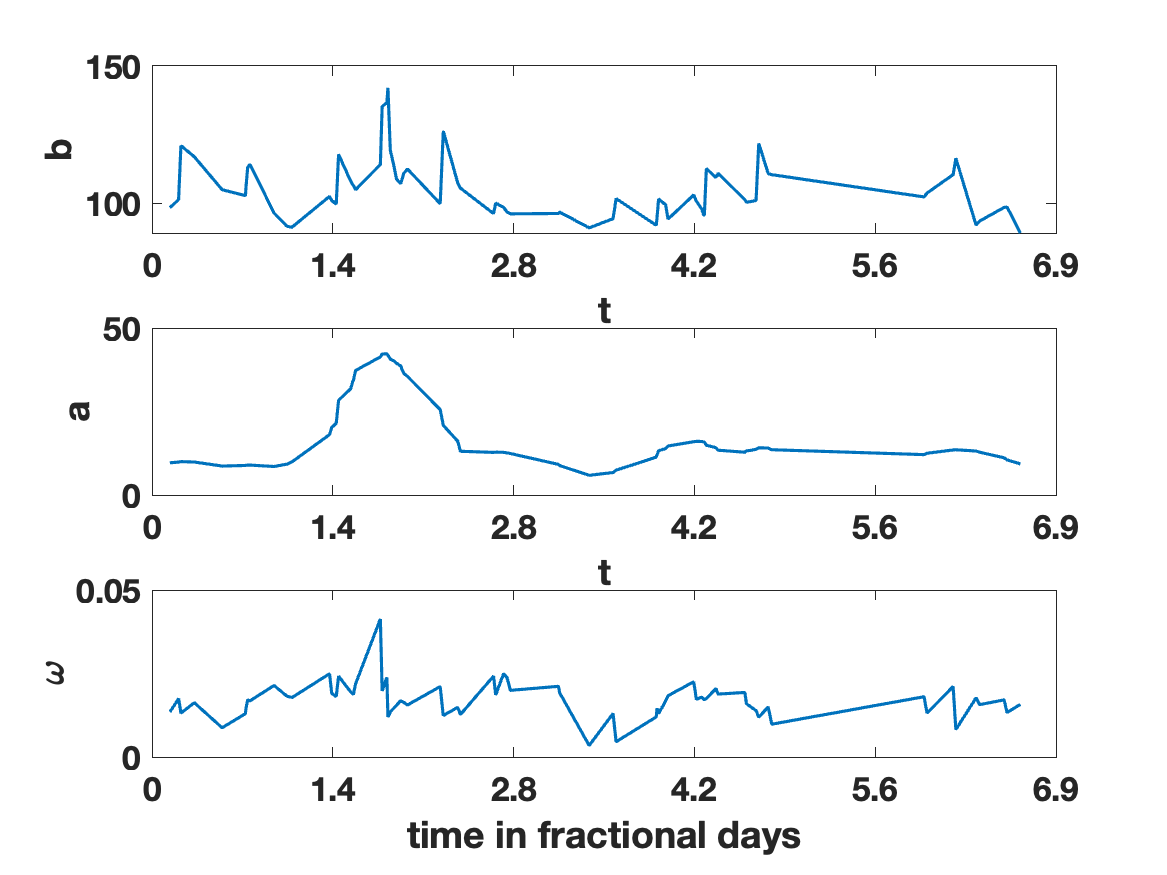}
	\end{subfigure}%
	\begin{subfigure}{0.33\textwidth}
		\includegraphics[width=\linewidth]{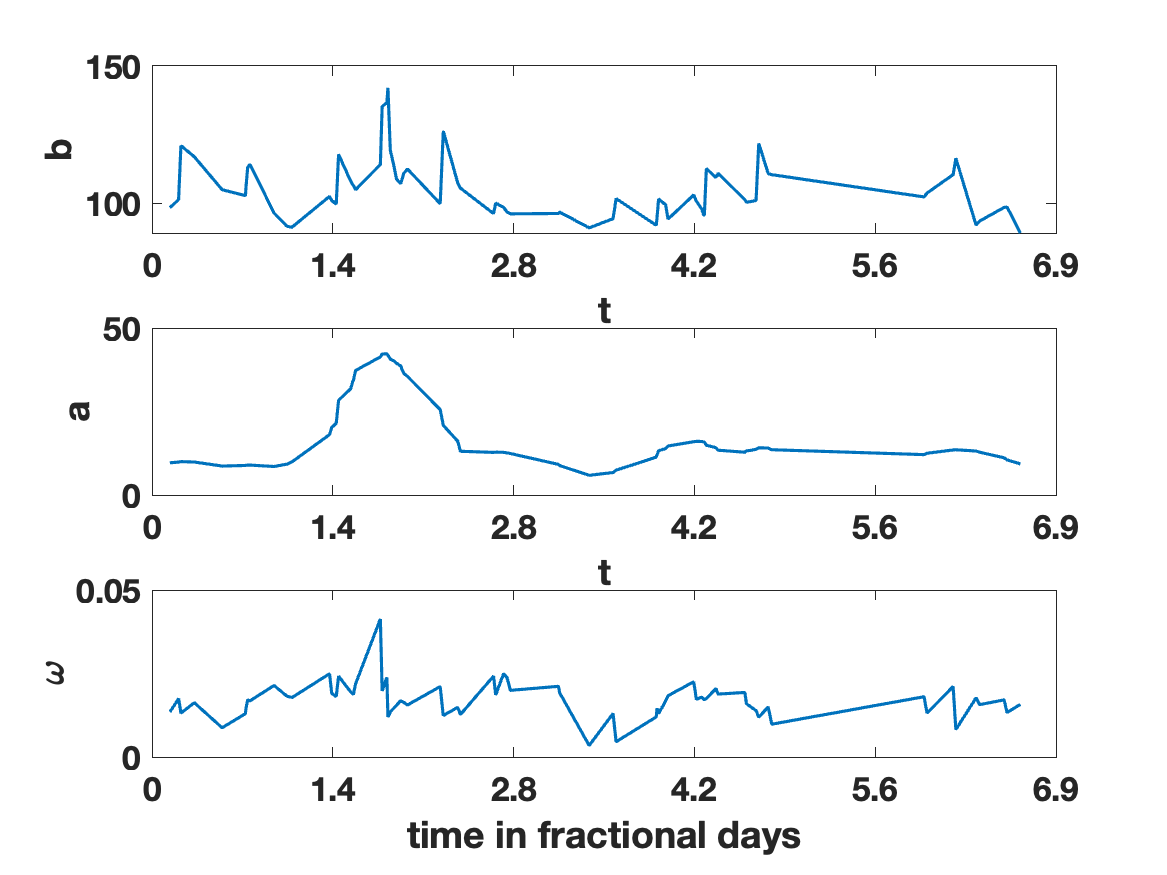}
	\end{subfigure}%
	\begin{subfigure}{0.33\textwidth}
		\includegraphics[width=\linewidth]{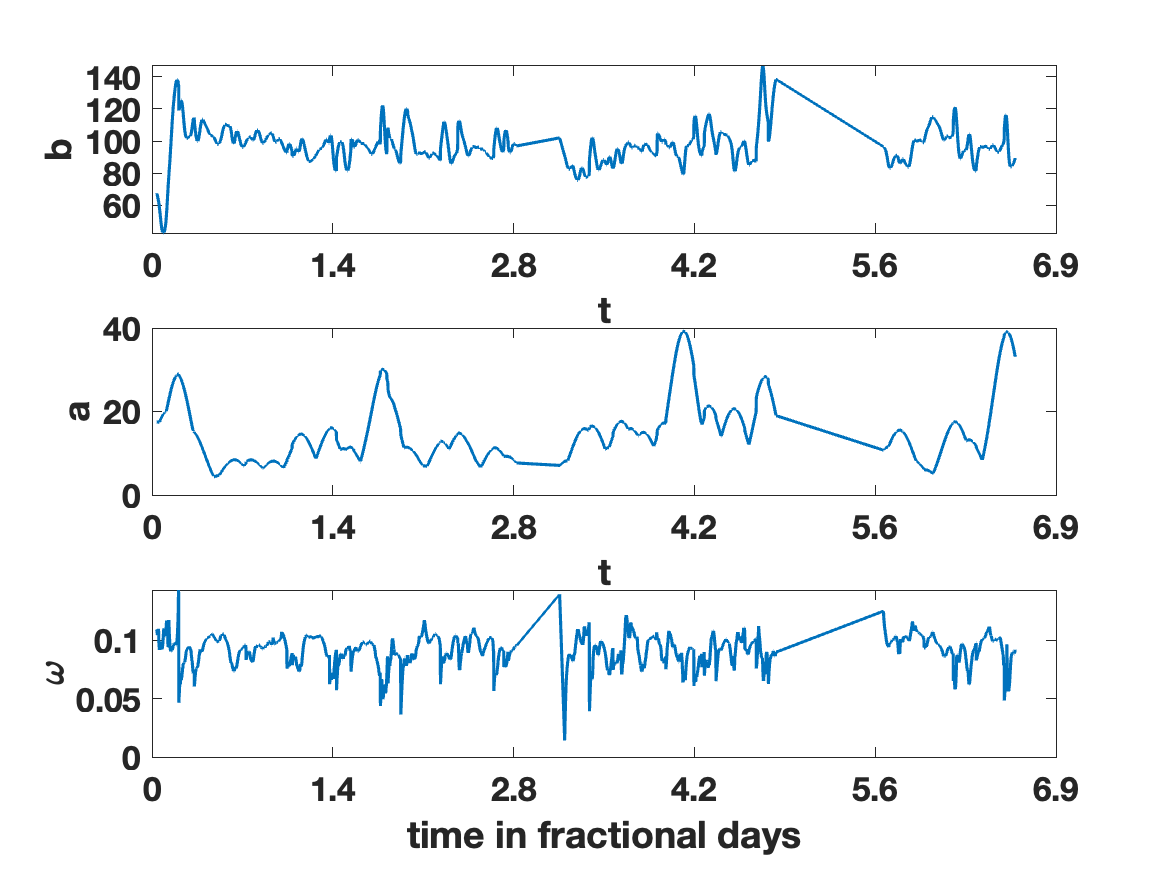}
	\end{subfigure}

	\caption{Given the in-the-wild glycemic damped, driven oscillatory dynamics we can see and the model basic state assuming \emph{non-oscillatory} dynamics, we can see that all measurement functions, $h_1$ (left), $h_2$ (center) and $h_3$ (right), utilized the temporal flexibility in the parameters to decrease estimation error, but in very different ways. Similarly to the case where the dynamics are assumed oscillatory, the sparse patient-defined and random measurement functions ($h_1$, $h_2$) had nearly the same flex over the course of the time window, and both were quite different from the dense measurement case ($h_3$). This result carries the same potential implications for solving inverse problems as shown in Fig. \ref{fig:wild_osc_parms}. Interestingly, the parameter trajectories and exercised parameter flex over the estimate window here is substantially lower than for the case where the underlying dynamics are assumed to be oscillatory in Fig. \ref{fig:wild_osc_parms}, as expected given that the assumed dynamics in that case were less representative of the generating dynamics.}
	\label{fig:wild_no_osc_parms}
\end{figure}

\paragraph{\textbf{In-the-wild}} In the case where the underlying dynamics were assumed to be oscillatory, Fig. \ref{fig:wild_osc_parms}, the sparse measurement functions, $h_1$ and $h_2$ produced nearly identical parameter trajectories, and both differed from the densely measured data produced by the CGM, $h_3$.  Whether the underlying dynamics were assumed to be oscillatory, Fig. \ref{fig:wild_osc_parms}, or not, Fig. \ref{fig:wild_no_osc_parms}, did not substantially impact either the estimation of the amplitude parameters, or the differences between the sparse and dense measurement function dependence on the amplitude parameter estimation. Not surprisingly, the primary dependencies for parameter estimation for a given type of assumed dynamics, oscillatory or not, appeared in the estimation of the frequency parameters.  \emph{The most important distinction} between parameter estimate trajectories for the in-the-wild case was between the assumed dynamics. The parameter variability over the course of the estimate window was substantially higher when the baseline model dynamics were oscillatory because the generating dynamics were not particularly oscillatory, and leading the transition probabilities to vary quite a lot to minimize the point-wise loss function. Overall, the parameter flex did impact model estimation, and the parameter flex was dependent on the measurement functions, but these differences differed from the highly non-stationary ICU case.

\section{Discussion}

\paragraph{\textbf{Overall summary}} Constructing a cost function structure with multiple cost functions supported an estimation balancing the minimization of point-wise errors with global properties such as qualitative dynamics. Estimation accuracy depended on the measurement processes, but for all measurement functions qualitative dynamics were preserved while minimizing point-wise errors. Additionally, the added intra-estimation window parameter flex allowed rigid mechanistic models to be robust to noice in measurement times, model error, and non-stationarity. 


\paragraph{\textbf{Data assimilation initialization}} Accurate initialization for data assimilation is crucial for immediate and accurate forecasting and for DA use with sequential smoothing and filtering. Our method provides a new path for DA initialization when data are particularly sparse and there is some knowledge of underlying model dynamics and state.

\paragraph{\textbf{Measurement functions and intra-window parameter trajectories}} Computed parameter trajectories were dependent on the measurement functions. It is not obvious how to identify whether a parameter trajectory is correct. In the case where the trajectory does not vary over the estimation window, this is not a problem. Similarly, if the parameter flex is only a matter of accounting for noise and the trajectory variability is not scientifically important, identification of intra-estimation-window trajectories is not a problem. But, if the intra-estimation-window parameter trajectories are scientifically or otherwise important, this opens a new problem of trying to understand and quantify the uncertainty in estimated parameter trajectories as they depend on measurement functions. 

\paragraph{\textbf{Data sparsity and non-stationarity}} Data sparsity and non-stationarity are deeply coupled concepts because non-stationarity can induce sparsity. If a system's parameters change substantially and faster than we measure, there will not exist statistically stable data set to estimate the model, and all data sets will be sparse. This induces a balance in estimation error due to use of small data sets against estimation errors due to data sets that are not statistically stable. The intra-estimation window parameter flexibility was designed to manage this complexity along with allowing a rigid model to be more robust to errors in measurement times.

\paragraph{\textbf{Managing data sparsity}} When data are not sparse, the estimation problem is relatively straightforward, and using both standard methods and our new method, minimizing point-wise errors can reproduce many global properties such as the distribution of model states. Our estimation results show that data sparsity, when the generating dynamics are complex, can have a substantial impact on our ability to accurately estimate the underlying dynamics and model states and parameters.  Model estimation with limited data decays as data become increasingly sparse, and it becomes difficult to quantify or identify what information we can reliably determine from sparse data. Our method was developed to help push the boundaries of the information we can extract from data when paired with external knowledge of the system.  And while we were able to estimate global and local properties of the system, explicit new open problems could be articulated related to identifying the limits of what information can be extracted and the uncertainty of that information. In the background, there is surely a \emph{Nyqist-Shannon-type theorem} for data assimilation and mechanistic modeling \cite{ALBERS2023104477} just as there has been for other spaces such as the reproducing kernel Hilbert space \cite{nyqist_shannon_rkhs1,nyqist_shannon_rkhs2}. Even rough numeric guidelines that could help limit the estimation errors given sparse data would be helpful, but such a result does not currently exist to our knowledge.

\paragraph{\textbf{Managing non-stationarity}} The impact of non-stationarity on the measurement functions and parameter estimation was interesting and differed depending on the generating dynamics.  When the system was non-stationary, here the ICU patients, it seems that random sparse measurements may lead to parameter estimation that is more representative to the parameters estimated in the densely measured case.  In contrast, the case where the system is relatively stationary, here the in-the-wild patient, differences in the sparse measurement functions translated to little variability in the intra-estimation-window parameter trajectories. Of course, two cases are not enough to understand and establish general rules of such measurement function dependent variability other than establishing that there can be important estimation dependencies on measurement functions. 
 
\paragraph{\textbf{Next steps}} There are at least five directions for advancement. \emph{First}, we need to construct a scalable and generalizable computational pipeline such that this methodology can be applied to any ODE and SDE models.  For example, because optimizing the $L_3$ and $L_4$ components of the objective function require derivatives of the model, integrating the methodology with auto-differentiation would be a meaningful advancement. \emph{Second,} our framework presents a new set of  variables to optimize over, the weights $\lambda_k$ attached to the components $L_k$ of the objective function.  While in this paper we set these weights externally, this could be further generalized to chose them automatically.  \emph{Third,} while some standard UQ methods apply to our methodology, there are many opportunities for UQ methods development. For example, devising verifiable methods for computing the model parameter estimate uncertainty given sparse data \emph{for the parameter trajectory and a parameter distribution estimated from the parameter trajectory} would be substantial advancements.  \emph{Fourth,} our method allows us to imbue the model with a given underlying dynamic such as oscillatory or fixed point dynamics. In our motivating example, we know from other work, what the baseline underlying dynamics are. However, learning more accurate initial parameters beyond a crude set of parameters that produce, e.g., oscillatory dynamics, could be devised as an iterative application of the methods devised here. And \emph{fifth,} given the motivation for the work in this paper,  adapting and applying our methodology directly to modeling patients in the ICU and other similar settings would be highly valuable.



\section{Conclusions}

Motivated by an inability to robustly estimate model parameters to initialize a DA-filter given sparse, non-stationary, noisy data, we developed a new methodology.  This methodology has three notable features.  \emph{First}, it balances optimizing local point-wise errors with global distributional errors and agreement with the model.  \emph{Second}, it allows users to imbue the model with known, underlying dynamics such as oscillatory dynamics. \emph{Third}, the method allows for the parameters to flex over an optimization window to manage both errors in measurement times and their impact on model rigidity and non-stationarity of the generating processes over the optimization window.  When we applied this method in a few contexts using both simulated (estimated) and real blood glucose data from two contexts including the ICU and a patient in the wild, the method was able in all cases to manage the sparse data issues, balance point-wise errors and global distributional errors while also allowing agreement with the model. This implies that we were able to robustly preserve global dynamics properties while also minimizing point-wise errors, a major goal of this work. We believe that this methodological pathway will be useful in inverse problems applications as well as embedded within sequential smoothing and DA-initialization applications.

\appendix

\section{Computing the derivative for the $L_k$'s}
\label{sec:dlk}

This appendix contains explicit formulas for the derivatives of the various components $L_k$ of the objective function, required for gradient descent. We write them down here for completeness and easy reproducibility, but also to stress a point: even for the relatively simple model that we have adopted for the examples in this article, computing and implementing the required derivatives is a laborious task, prone to errors. In order to add flexibility to the selection of models, one should consider resorting to automatic differentiation.

\smallskip

%
$$ L_1 = \frac{1}{n} \sum_j \log\left[(1-\epsilon) K^y\left(y^j , x^j \right) + \epsilon \frac{1}{n} \sum_{i=1}^n K^y\left(y^j, y^i\right)\right] $$ 

$$ \frac{\partial L_1}{\partial x^j} = \frac{1}{n} \frac{1-\epsilon}{h^2} \frac{\left(y^j - x^j\right) K^y\left(y^j , x^j \right)}{(1-\epsilon) K^y\left(y^j , x^j \right) + \epsilon \rho_0^j} $$



\smallskip


%
\begin{eqnarray*}
	L_2 &=&  -\frac{1}{n}\sum_i \frac{1}{\sum_l K^t\left(t^i, t^l\right)} \sum_j \left[\left(K^y\left(x^i , x^j \right) - K^y\left(y^i , x^j \right)\right) - \left(K^y\left(x^i , y^j \right) - K^y\left(y^i , y^j \right)\right)\right] K^t\left(t^i, t^j\right) \\
	&=&
	-\frac{1}{2 n} \sum_{i,j} \left(\frac{K^t\left(t^i, t^j\right)}{\sum_l K^t\left(t^j , t^l \right)} 
	+ \frac{K^t\left(t^i, t^j\right)}{\sum_l K^t\left(t^i , t^l \right)}\right) \left[K^y\left(x^i , x^j \right) - 2 K^y\left(y^i , x^j \right)
	+ K^y\left(y^i , y^j \right) \right] 
\end{eqnarray*}

$$ \frac{\partial L_2}{\partial x^j} = -\frac{1}{n h^2} \sum_i \left(\frac{K^t\left(t^i, t^j\right)}{\sum_l K^t\left(t^j , t^l \right)} 
+ \frac{K^t\left(t^i, t^j\right)}{\sum_l K^t\left(t^i , t^l \right)}\right) \left[\left(x^i - x^j\right) K^y\left(x^i , x^j \right) - \left(y^i - x^j\right) K^y\left(y^i , x^j \right) \right] $$

\smallskip

For $L_3$ and $L_4$, the components of the objective function that enforce compliance with the model, the pairs $(x^{j-1}, z^{j-1})$ are replaced by $(r^{j-1}, \theta^{j-1})$, and we have, with all variables evaluated at time $t_{j-1}$, 
$$ \frac{\partial}{\partial x} =  \frac{x-b}{r} \frac{\partial}{\partial r} - \frac{z}{r^2}  \frac{\partial}{\partial \theta}, $$
$$ \frac{\partial}{\partial z} =  \frac{z}{r} \frac{\partial}{\partial r} + \frac{x-b}{r^2}  \frac{\partial}{\partial \theta}, $$
$$ \frac{\partial}{\partial b} =  -\frac{x-b}{r} \frac{\partial}{\partial r} + \frac{z}{r^2}  \frac{\partial}{\partial \theta}. $$

\smallskip

$$ \rho^x\left(x^j | r^{j-1}, \theta^{j-1}, \Delta t^j, \alpha^{j} \right) = \frac{1}{\sqrt{2\pi} \sigma} e^{-\frac{\left(x^j - b^j - r_+^j \cos\left(\theta^{j-1} + \omega^{j-1} \Delta t^j\right)  \right)^2}{2 \sigma^2}},$$
$$ \rho^z\left(z^j | r^{j-1}, \theta^{j-1}, \Delta t^j, \alpha^{j} \right) = \frac{1}{\sqrt{2\pi} \sigma} e^{-\frac{\left(z^j - r_+^j \sin\left(\theta^{j-1} + \omega^{j-1} \Delta t^j\right)  \right)^2}{2 \sigma^2}},$$
where
$$ r_+^{j} = \left(1 - d_s^{j}\right) a^{j} + d_s^{j} \ r^{j-1} , \quad d_s^{j} = e^{-\frac{\Delta t^{j}}{T_s}}. $$
Then
$$ \frac{\partial \rho^x\left(x^j | r^{j-1}, \theta^{j-1}, \Delta t^j, \alpha^j \right)}{\partial x^j}  = - \frac{x^j - b^j - r_+^j \cos\left(\theta^{j-1} + \omega^{j-1} \Delta t^j\right)}{\sigma^2} \rho^x,$$
$$ \frac{\partial \rho^x\left(x^j | r^{j-1}, \theta^{j-1}, \Delta t^j, \alpha^j \right)}{\partial r^{j-1}}  = -d_s^j
\cos\left(\theta^{j-1} + \omega^{j-1} \Delta t^j\right)
\frac{\partial \rho^x}{\partial x^j}, $$
$$ \frac{\partial \rho^x\left(x^j | r^{j-1}, \theta^{j-1}, \Delta t^j, \alpha^j \right)}{\partial \theta^{j-1}}  = 
r_+^j
\sin\left(\theta^{j-1} + \omega^{j-1} \Delta t^j\right)
\frac{\partial \rho^x}{\partial x^j}, $$
$$ \frac{\partial \rho^x\left(x^j | r^{j-1}, \theta^{j-1}, \Delta t^j, \alpha^j \right)}{\partial b^j}  = 
-  \frac{\partial \rho^x\left(x^j | r^{j-1}, \theta^{j-1}, \Delta t^j, \alpha^j \right)}{\partial x^j}, $$  
%
$$ \frac{\partial \rho^x\left(x^j | r^{j-1}, \theta^{j-1}, \Delta t^j, \alpha^j \right)}{\partial a^{j}}  = 
- \frac{\partial \rho^x\left(x^j, \theta^j | r^{j-1}, \theta^{j-1}, \Delta t^j, \alpha^j \right)}{\partial r^{j-1}}, $$
$$ \frac{\partial \rho^x\left(x^j | r^{j-1}, \theta^{j-1}, \Delta t^j, \alpha^j \right)}{\partial \omega^{j-1}}  =
\Delta t^j \frac{\partial \rho^x\left(x^j | r^{j-1}, \theta^{j-1}, \Delta t^j, \alpha^j \right)}{\partial \theta^{j-1}} , $$
$$ \frac{\partial \rho^z\left(z^j | r^{j-1}, \theta^{j-1}, \Delta t^j, \alpha^j \right)}{\partial z^j}  = - \frac{z^j - r_+^j \sin\left(\theta^{j-1} + \omega^{j-1} \Delta t^j\right)}{\sigma^2} \ \rho^z,$$
$$ \frac{\partial \rho^z\left(z^j, \theta^j | r^{j-1}, \theta^{j-1}, \Delta t^j, \alpha^j \right)}{\partial r^{j-1}}  = -d_s^j
\sin\left(\theta^{j-1} + \omega^{j-1} \Delta t^j\right)
\frac{\partial \rho^z}{\partial z^j}, $$
$$ \frac{\partial \rho^z\left(z^j | r^{j-1}, \theta^{j-1}, \Delta t^j, \alpha^j \right)}{\partial \theta^{j-1}}  = 
-r_+^j
\cos\left(\theta^{j-1} + \omega^{j-1} \Delta t^j\right)
\frac{\partial \rho^z}{\partial z^j} $$
$$ \frac{\partial \rho^z\left(z^j | r^{j-1}, \theta^{j-1}, \Delta t^j, \alpha^j \right)}{\partial a^{j}}  = 
- \frac{\partial \rho^z\left(z^j, \theta^j | r^{j-1}, \theta^{j-1}, \Delta t^j, \alpha^j \right)}{\partial r^{j-1}}, $$
$$ \frac{\partial \rho^z\left(z^j | r^{j-1}, \theta^{j-1}, \Delta t^j, \alpha^j \right)}{\partial \omega^{j-1}}  =
\Delta t^j \frac{\partial \rho^z\left(z^j | r^{j-1}, \theta^{j-1}, \Delta t^j, \alpha^j \right)}{\partial \theta^{j-1}} . $$

\smallskip


%
$$ L_3 = \frac{1}{n} \sum_j \log\left[\rho^x\left(x^j | x^{j-1}, z^{j-1}, \Delta t^j, \alpha^j \right)\right],$$

$$ \frac{\partial L_3}{\partial x^j} = \frac{1}{n} \left[ \frac{\frac{\partial}{\partial x^j} \rho^x\left(x^j | x^{j-1}, z^{j-1}, \Delta t^j, \alpha^j \right)}{\rho^x\left(x^j | x^{j-1}, z^{j-1}, \Delta t^j, \alpha^j \right)} + \frac{\frac{\partial}{\partial x^j} \rho^x\left(x^{j+1} | x^{j}, z^{j}, \Delta t^{j+1}, \alpha^{j+1} \right)}{\rho^x\left(x^{j+1} | x^{j}, z^{j}, \Delta t^{j+1}, \alpha^{j+1} \right)} \right]$$

$$ \frac{\partial L_3}{\partial z^j} = \frac{1}{n} \frac{\frac{\partial}{\partial z^j} \rho^x\left(x^{j+1} | z^{j}, \Delta t^{j+1}, \alpha^{j+1} \right)}{\rho^x\left(x^{j+1} | x^{j}, z^{j}, \Delta t^{j+1}, \alpha^{j+1} \right)}$$

$$ \frac{\partial L_3}{\partial \alpha_k^j} = \frac{1}{n} \frac{\frac{\partial}{\partial \alpha_k} \rho^x\left(x^{j} | x^{j-1}, z^{j-1}, \Delta t^{j}, \alpha^{j} \right)}{\rho^x\left(x^{j} | x^{j-1}, z^{j-1}, \Delta t^{j}, \alpha^{j} \right)}$$ 

\smallskip


%
$$ L_4 = \frac{1}{n} \sum_j \log\left[\rho^z\left(z^j | x^{j-1}, z^{j-1}, \Delta t^j, \alpha^j \right)\right]$$ 

$$ \frac{\partial L_4}{\partial x^j} =  \frac{1}{n}\frac{\frac{\partial}{\partial x^j} \rho^z\left(z^{j+1} | z^{j}, \Delta t^{j+1}, \alpha^{j+1} \right)}{\rho^z\left(z^{j+1} | x^{j}, z^{j}, \Delta t^{j+1}, \alpha^{j+1} \right)}$$

$$ \frac{\partial L_4}{\partial z^j} =  \frac{1}{n} \left[ \frac{\frac{\partial}{\partial z^j} \rho^z\left(z^j | x^{j-1}, z^{j-1}, \Delta t^j, \alpha^j \right)}{\rho^z\left(z^j | x^{j-1}, z^{j-1}, \Delta t^j, \alpha^j \right)} + \frac{\frac{\partial}{\partial z^j} \rho^z\left(z^{j+1} | x^{j}, z^{j}, \Delta t^{j+1}, \alpha^{j+1} \right)}{\rho^z\left(x^{j+1} | x^{j}, z^{j}, \Delta t^{j+1}, \alpha^{j+1} \right)}\right]$$

$$ \frac{\partial L_4}{\partial \alpha_k^j} = \frac{1}{n} \frac{\frac{\partial}{\partial \alpha_k^j} \rho^z\left(z^{j} | x^{j-1}, z^{j-1}, \Delta t^{j}, \alpha^{j} \right)}{\rho^z\left(z^{j} | x^{j-1}, z^{j-1}, \Delta t^{j}, \alpha^{j} \right)}$$ 

\smallskip


%
$$ L_{4+k} = \frac{1}{n} \sum_j \log\left[\rho^k\left(\alpha^j_k | \alpha^{j-1}_k, \Delta t^j\right) \right], $$
$$ \frac{\partial L_{4+k}}{\partial \alpha_k^j} = \frac{1}{n} \left[\frac{\frac{\partial}{\partial \alpha_k^j} \rho^k\left(\alpha_k^{j} | \alpha_k^{j-1}\right)}{\rho^l\left(\alpha_k^{j} | \alpha_k^{j-1} \right)}
+ \frac{\frac{\partial}{\partial \alpha_k^j} \rho^k\left(\alpha_k^{j+1} | \alpha_k^{j}\right)}{\rho^k\left(\alpha_k^{j+1} | \alpha_k^{j} \right)} \right], $$ 

$$ \rho^k\left(\alpha_k^{j} | \alpha_k^{j-1}\right) = 
N\left[d_l^j \alpha_l^{j-1} + \left(1-d_l^j\right) \tilde{\alpha}_l , \left(1-d_l^j\right) {\sigma_l}^2 \right], $$ 

$$ \frac{\partial \rho^k\left(\alpha_k^{j} | \alpha_k^{j-1}\right)}{\partial \alpha_k^j} =
-\frac{\alpha_k^j - \left(d_l^j \alpha_l^{j-1} + \left(1-d_l^j\right) \tilde{\alpha}_l  \right)}{\left(1-d_l^j\right) {\sigma_l}^2}  \rho^k\left(\alpha_k^{j} | \alpha_k^{j-1}\right),$$

$$ \frac{\partial \rho^k\left(\alpha_k^{j} | \alpha_k^{j-1}\right)}{\partial \alpha_k^{j-1}} = -d_l^j \frac{\partial \rho^k\left(\alpha_k^{j} | \alpha_k^{j-1}\right)}{\partial \alpha_k^j}. $$

\section{Ultradian glucose-insulin model}

The primary state variables include the plasma glucose concentration $G$, the plasma insulin concentration $I_{p}$, and the interstitial insulin concentration $I_{i}$. Additionally there is a delay in the hepatic response of plasma insulin to glucose approximated using the linear chain trick resulting in $(h_1,h_2,h_3)$  \cite{sturis_91}. The system of ordinary differential equations take the form \cite{keenerII}:
\begin{subequations}
	\label{eq:model1}
	\begin{eqnarray}
		\frac{dI_p}{dt} & = &  f_1(G)-E\bigl(\frac{I_{p}}{V_{p}}-\frac{I_i}{V_{i}}\bigr)-\frac{I_{p}}{t_{p}}\\
		\frac{dI_i}{dt} & = & E\bigl(\frac{I_{p}}{V_{p}}-\frac{I_i}{V_{i}}\bigr)-\frac{I_{i}}{t_{i}}\\
		\frac{dG}{dt} & = & f_4(h_3)+I_{G}(t)-f_2(G)-f_3(I_i)G\\
		\frac{dh_1}{dt} & = & \frac{1}{t_d}\bigl(I_p-h_1\bigr) \\
		\frac{dh_2}{dt} & = & \frac{1}{t_d}\bigl(h_1-h_2\bigr) \\
		\frac{dh_3}{dt} & = & \frac{1}{t_d}\bigl(h_2-h_3\bigr)
	\end{eqnarray}
\end{subequations}
This model includes many parameterized processes including the rate of insulin production, $f_1(G)$, insulin-independent glucose utilization $f_2(G)$, insulin-dependent glucose utilization, $f_3(I_i)G$, and  represents insulin-dependent glucose utilization, delayed insulin-dependent glucose utilization,
$f_4(h_3)$, defined by:
\begin{eqnarray}
	f_1(G) & = & \frac{R_m}{1+ \exp(\frac{-G}{V_g c_1} + a_1)} \\
	f_2(G) & = & U_b(1-\exp(\frac{-G}{C_2V_g})) \\
	f_3(I_i) & = & \frac{1}{C_3 V_g}( U_0 + \frac{U_m - U_0}{1 + (\kappa I_i)^{-\beta}}) \\
	f_4(h_3) & = & \frac{R_g}{1 + \exp(\alpha (\frac{h_3}{C_5 V_p}  -1))} \\
	\kappa & = & \frac{1}{C_4} (\frac{1}{V_i} - \frac{1}{E t_i}),
\end{eqnarray}
respectively.

The nutritional driver of the model $I_G(t)$ is defined over $N$ non-overlapping intervals where nutrition is delivered through an enteral tube at constant rates.

\begin{eqnarray}
	I_G(t) = \sum^N_{j=1}m_j\mathbbm{1}_{\{t_j\leq t<t_{j+1}\}}(t),
\end{eqnarray}
where $m_j$ is the carbohydrate rate (mg/min) delivered over the interval $[t_j,t_{j+1})$ and $\mathbbm{1}(\cdot)$ is the characteristic function.

\begin{table}[!ht]
	\centering
	\caption{Full list of parameters for the ultradian glucose-insulin model \cite{keenerII}. Note that IIGU and IDGU denote insulin-independent glucose utilization and insulin-dependent glucose utilization, respectively.}
	\begin{tabular}{|p{1.2cm}|p{2.7cm}|p{8cm}|}
		\hline
		\multicolumn{3}{|p{8cm}|}{\textbf{Ultradian model parameters}} \\ \hline
		\cline{1-3} \hline
		\cline{1-3}
		Name & Nominal Value  & Meaning \\ \hline \hline
		$V_p$  & $3$ l  & plasma volume  \\ \hline \hline
		$V_i$  & $11$ l  & interstitial  volume \\ \hline \hline
		$V_g$ & $10$ l  & glucose space \\ \hline \hline
		$E$  & $0.2$ l min$^{-1}$ &   exchange rate for insulin between remote
		and plasma compartments \\ \hline \hline
		$t_p$  & $6$ min  & time constant for plasma insulin degradation (via
		kidney and liver filtering) \\ \hline \hline
		$t_i$  & $100$ min & time constant for remote insulin degradation (via
		muscle and adipose tissue) \\ \hline \hline
		$t_d$  & $12$ min  & delay between plasma insulin and glucose
		production \\ \hline \hline
		$k$  & $0.5$ min$^{-1}$  & rate of decayed appearance of ingested glucose \\ \hline \hline
		$R_m$  & $209$ mU min$^{-1}$  & linear constant affecting insulin secretion  \\ \hline \hline
		$a_1$  & $6.6$ & exponential constant affecting insulin secretion \\ \hline \hline
		$C_1$  & $300$ mg l$^{-1}$ & exponential constant affecting insulin secretion \\ \hline \hline
		$C_2$  & $144$ mg l$^{-1}$  & exponential constant affecting IIGU \\ \hline \hline
		$C_3$  & $100$ mg l$^{-1}$  & linear constant affecting IDGU \\ \hline \hline
		$C_4$  & $80$ mU l$^{-1}$ & factor affecting IDG \\ \hline \hline
		$C_5$  & $26$ mU l$^{-1}$  & exponential constant affecting IDGU \\ \hline \hline
		$U_b$  & $72$ mg min$^{-1}$  & linear constant affecting IIGU \\ \hline \hline
		$U_0$  & $4$ mg min$^{-1}$ & linear constant affecting IDGU \\ \hline \hline
		$U_m$  & $94$ mg min$^{-1}$  &  linear constant affecting IDGU \\ \hline \hline
		$R_g$  & $180$ mg min$^{-1}$  & linear constant affecting IDGU \\ \hline \hline
		$\alpha$  & $7.5$ & exponential constant affecting IDGU \\ \hline \hline
		$\beta$  & $1.772$ & exponent affecting IDGU \\ \hline \hline
	\end{tabular}
	\label{table:model_parameters}
\end{table}

\end{document}

%% file: ex_shared.tex
\usepackage{lipsum}
\usepackage{amsfonts}
\usepackage{graphicx}
\usepackage{epstopdf}
\usepackage{algorithmic}
\ifpdf
  \DeclareGraphicsExtensions{.eps,.pdf,.png,.jpg}
\else
  \DeclareGraphicsExtensions{.eps}
\fi

\newsiamremark{remark}{Remark}
\newsiamremark{hypothesis}{Hypothesis}
\crefname{hypothesis}{Hypothesis}{Hypotheses}
\newsiamthm{claim}{Claim}

\headers{A new approach to data assimilation initialization}{D. J. Albers, G. Hripcsak, L. Mamykina, M. Sirlanci, E. Tabak}

\title{A new approach to data assimilation initialization problems with sparse data using multiple cost functions\thanks{Funding: G.H. is supported by NIH R01 LM006910. L. M. is suported by NIH R01 DK113189. E. G. T. is supported in part by the Office of Naval Research, through grant N00014-22-1-2192.}}

\author{David J. Albers\thanks{Department of Biomedical Informatics, University of Colorado Anschutz Medical Campus, Aurora, CO
  (\email{david.albers@cuanschutz.edu}, \email{melike.sirlanci@cuanschutz.edu}).}
\and George Hripcsak\thanks{Department of Biomedical Informatics, Columbia University, New York, NY (\email{hripcsak@columbia.edu}, \email{om2196@cumc.columbia.edu}).}
\and Lena Mamykina\footnotemark[3]
\and Melike Sirlanci\footnotemark[2]
\and Esteban G. Tabak\thanks{Courant Institute of Mathematical Sciences, New York University, New York, NY
	(\email{tabak@cims.nyu.edu}).}}

\usepackage{amsopn}

\makeatletter
\newcommand*{\addFileDependency}[1]{
  \typeout{(#1)}
  \@addtofilelist{#1}
  \IfFileExists{#1}{}{\typeout{No file #1.}}
}
\makeatother
